\documentclass[11pt,a4paper]{article}
\usepackage[L7x]{fontenc}
\usepackage[utf8]{inputenc}

\usepackage[pdftex]{graphicx}
\usepackage{amssymb, amsmath}
\usepackage{amsthm}
\usepackage{dsfont}
\usepackage[pdftex,bookmarks=true,unicode]{hyperref}
\usepackage{booktabs}
\usepackage{color}
\usepackage{hyperref,footnote}
\usepackage[bottom]{footmisc}
\hypersetup{hidelinks}
\usepackage{booktabs}
\usepackage{csquotes}
\usepackage[toc,page]{appendix}
\usepackage{microtype}
\DisableLigatures{encoding = *, family = * }
\usepackage[margin=0.2cm]{caption}
\usepackage{enumerate,bbm}

\MakeOuterQuote{"}
\usepackage[left=2cm, top=2cm, right=1.5cm, bottom=2cm] {geometry}
\usepackage{indentfirst}
\theoremstyle{plain}

\theoremstyle{definition}

 \newcommand{\PP}{\mathbb{P} }
 \newcommand{\NN}{\mathbb{N} }
  \newcommand{\RR}{\mathbb{R} }
 \newcommand{\EE}{\mathbb{E} }
  \newcommand{\ZZ}{\mathbb{Z} }
  \newcommand{\al}{\alpha }
  \newcommand{\w}{V }
  \newcommand{\bb}{\beta}
    \newcommand{\g}{\gamma}
\newcommand{\e}{\varepsilon }
\newcommand{\dt}{{\rm d}t }
\newcommand{\ds}{{\rm d}s }
\newcommand{\dd}{{\rm d}}
 \newcommand{\abs}[1]{\left|#1\right| }
   \newcommand{\sv}[1]{\lfloor #1 \rfloor}

  \newcommand{\ind}[1]{\mathbbm{1}_{#1}}
  
\usepackage{lineno,enumitem}

\def\ba{{\bf a}}
\def\bk{{\bf k}}
\def\bi{{\bf i}}
\def\bj{{\bf j}}
\def\bz{{\bf z}}
\def\bx{{\bf x}}
\def\BO{{\bf 0}}

\def\phi{\varphi}

\def\bbr{{\mathbb R}}
\def\bbs{{\mathbb S}}

\def\bX{{\mathbf{X}}}

\def\bz{{\bf Z}}

\def\bz{{\mathbb Z}} 

\def\bj{{\bf j}}
\def\bk{{\bf k}}
\def\bt{{\bf t}}
\def\bs{{\bf s}}

\def\bx{{\bf x}}
\def\bX{{\bf X}}

\def\bi{{\bf i}}
\def\Bl1{{\bf 1}}
\def\B2{{\bf 2}}
\def\B0{{\bf 0}}

\newcommand{\R}{\mathbb{R}}

\def\G{\Gamma}
\def\y{y}

\newtheorem{teo}{Theorem}

\begin{document}

	\title{Asymptotic of spectral covariance for linear
		random fields  with infinite variance
		\footnotemark[0]\footnotetext[0]{ \textit{Short title:}
			spectral covariance for fields with infinite variance }
		\footnotemark[0]\footnotetext[0]{%
			\textit{MSC 2000 subject classifications}. Primary 60G60, secondary
			 60E07 .} \footnotemark[0]
		\footnotemark[0]\footnotetext[0]{ \textit{Corresponding author:}
			Vygantas Paulauskas, Department of Mathematics and
			Informatics, Vilnius university, Naugarduko 24, Vilnius 03225, Lithuania,
			e-mail:vygantas.paulauskas@mif.vu.lt}
	}

	\author{Julius Damarackas$^{\text{\small 1}}$,   Vygantas Paulauskas$^{\text{\small 1}}$ \\
		{\small $^{\text{1}}$ Vilnius University, Department of Mathematics
			and
			Informatics,}\\}

	
	
	\maketitle
	
	\begin{abstract}
		
		In the paper we continue to investigate measures of dependence for random variables with infinite variance. The asymptotic of spectral  covariance $\rho (X_{(0,0)}, X_{(k_1,k_2)})$  for linear random field $X_{k,l}=\sum_{i,j=0}^\infty c_{i,j}\e_{k-i, l-j}, \ (k, l)\in \bz^2,$ with special form of filter $\{c_{i,j}\}$ and with innovations $\{ \e_{i, j}\}$ having infinite second moment is investigated. Different behavior of $\rho (X_{(0,0)}, X_{(k_1,k_2)})$ is obtained in the cases $n\to \infty, \ m\to \infty$ and $n\to \infty, \ m\to -\infty$, the latter case being much more complicated.
	\end{abstract}
	
	\vfill
	\eject
	\section{Introduction }

	The present paper can be considered as continuation of the paper \cite{PaulDam} and further development of the theory of spectral covariances which was outlined in the preprint \cite{Paul15}. In \cite{PaulDam} spectral covariances were calculated for several classes of stochastic processes with infinite second moment. In this paper we investigate spectral covariances for linear random fields without finite second moment.
	
We recall the notion of spectral covariance, mainly repeating definitions from \cite{PaulDam}.  Let $\bbs^{d-1}=\{x\in \R^d:
	\|x\|=1\}$ be the unit sphere in $\R^d$, here $\|x\|$ stands for the
	Euclidean norm in $\R^d$. Random vector $\bX =(X_1, \dots X_d)$ is
	$\al$-stable with parameter $0<\al<2$ if there exists a  finite
	measure $\G$ on  $\bbs^{d-1}$ and a vector $\mu \in \R^d$ such, that the
	characteristic function (ch.f.) of $\bX$, in the case $\al \ne 1,$ is
	given by formula
	\begin{equation}\label{stabdef}
	\EE\exp \left \{i(\bt, \bX) \right \}=\exp \left \{-\int_{\bbs^{d-1}}|(\bt,
	\bs)|^\al \left (1-{\rm sign} (\bt, \bs)\tan \frac{\pi \al}{2} \right
	)\G (d\bs)+i(\bt, \mu) \right \}.
	\end{equation}
	In the case $\al=1$  we shall consider only symmetric measure $\G$.
	In the case where $\G$ is symmetric, we get the
	so-called symmetric $\al$-stable ($S\al S$) distributions, and ch.f.
	takes very simple form
	\begin{equation}\label{stabsymdef}
	\EE\exp \left \{i(\bt, \bX) \right \}=\exp \left \{-\int_{\bbs^{d-1}}|(\bt,
	\bs)|^\al \G (d\bs) \right \}.
	\end{equation}
	$\G$ in (\ref{stabdef}) is called the spectral measure of an $\al$-stable random vector
	$\bX$, and the pair $(\G, \mu)$ is unique. The Gaussian case $\al=2$
	is excluded from this definition, since in the Gaussian case there
	is no uniqueness of the spectral measure $\G$: many different
	measures $\G$ will give the same ch.f.. In what follows we shall
	assume that $\mu=0$. Taking $d=2$ and $\mu=0$ in (\ref{stabdef}) we
	have $\al $-stable random vector $\bX =(X_1, X_2)$ with spectral
	measure $\G$ on $\bbs^1=\{x\in \R^2: x_1^2+x_2^2=1\}$ and we define
	spectral covariance of random vector $\bX$ (or spectral covariance between
	coordinates $X_1$ and $X_2$) as follows
	\begin{equation} \label{acovdef}
	{\rho}={\rho} (X_1, X_2) = \int_{\bbs^{1}}s_1s_2{\G}(ds).
	\end{equation}
	Also in analogy with the usual correlation coefficient we define
	spectral correlation coefficient (s.c.c.) for an $\al $-stable random vector
	$\bX$
	\begin{eqnarray} \label{gapdef}
	{\tilde \rho}={\tilde \rho} (X_1, X_2)
	=\int_{\bbs^1}s_1s_2{\G}(ds)\left
	(\int_{\bbs^{1}}s_1^2{\G}(ds)\int_{\bbs^{1}}s_2^2{\G}(ds) \right
	)^{-1/2}.
	\end{eqnarray}

	Let us note that this notion of spectral covariance was introduced in  \cite{Paul3} (only for $S\al S$ random vectors and under the name of "generalized association parameter"), in \cite{Paul15} it was called $\al$-covariance and only in \cite{PaulDam} it was named "spectral covariance".
	In  \cite{Paul3} (see Proposition 1 in \cite{PaulDam})
	it was shown that the introduced s.c.c. has the same properties as the usual correlation coefficient in the case $\al=2$.
	Also in \cite{PaulDam} and \cite{Paul3}  it was shown that the notion of spectral covariance can be extended naturally to $\al $-stable random vector with values in $\R^d$ or even in a separable Banach space, while two other measures of dependence - covariation and codifference - were defined only for bivariate $S\al S$ random vectors (we do not provide here the definitions of these two notions, a reader can find them in \cite{Samorod} or in \cite{PaulDam}).
	One of the main advantages of the spectral covariance is the fact that it can be extended to all random vectors belonging to the  domain of  attraction of  an $\al $-stable random vector in the following way.  Let $\xi=(\xi_1,\xi_2)$ be a random vector satisfying the following condition:
	there exists a number $0<\al\le 2,$ a slowly varying function L, and a finite  measure $\G$ on $\bbs^{1}$ such that
	\begin{equation}\label{domofattr}
	\lim_{x\to \infty} \frac{x^\al}{L(x)} \PP\left (\|\xi\|>x, \ \xi\|\xi\|^{-1}\in A \right )=\G (A),
	\end{equation}
	for any  Borel set $A$ on $\bbs^{1}$ with $\G(\partial A)=0.$ It is well-known that this condition is necessary and sufficient for  the random vector $(\xi_1,
	\xi_2)$ to belong to the domain of  attraction of  an $\al $-stable random
	vector $\bX=(X_1, X_2)$ with the exponent $\al$ and the spectral measure $\G$. We suggest
	to define spectral covariance and
	s.c.c. of $(\xi_1, \xi_2)$ by means of the measure $\G$, as these quantities
	are defined for $\al$-stable random vector $\bX=(X_1, X_2)$:
	\begin{equation}\label{acovardef2}
	\rho (\xi_1, \xi_2)=\rho (X_1, X_2)=\int_{\bbs^{1}}s_1s_2{\G}(ds)
	\end{equation}
	and similarly for ${\tilde \rho}(\xi_1, \xi_2)$, see \cite{PaulDam} for motivation of such definition.
	
	In the paper we consider linear random fields
	\begin{equation}\label{linfield}
	X_{\bk}=\sum_{\bi \in (\bz^+)^d} c_{\bi}\e_{\bk-\bi}, \ \bk\in \bz^d,
	\end{equation}
	where $(\bz^+)^d=\{\bi: i_j \ge 0, j=1, \dots ,d \}$ and $\e_{\bi}, \ \bi \in \bz^d, $ are
	i.i.d.  random variables with $S\al S$ distribution  (as in the case of random processes this assumption is made only for the simplicity of writing), and $c_{\bi}, \ \bi\in (\bz^+)^d$, are real numbers satisfying
	condition
	\begin{equation}\label{cond2}
	\sum_{\bi\in (\bz^+)^d} |c_{\bi}|^\al <\infty.
	\end{equation}
	Here and in what follows letters in bald stand for vectors in $\bbr^d$, and inequalities or equalities
	between them is componentwise. We are interested in the expression of $\rho (X_{\B0}, X_{\bk})$ via the filter $\{c_{\bi}\}$ and the asymptotic behavior of this quantity as $\bk \to \infty$ (we shall understand this relation as $\min_{1\le i\le d}k_i \to \infty $). For linear random processes ($d=1, \ X_{k}=\sum_{i \in \bz^+} c_{i}\e_{k-i}, \ k\in \bz$) we easily got the expression of $\rho (X_0, X_n)$  and the asymptotic decay of this quantity (see (2.2) and  Theorem 2 in \cite{PaulDam}), but it turned out that the generalization to the case $d>1$ is not so trivial. Since the main difficulties in passing from the case $d=1$ to the case $d>1$ can be easier understandable in the case $d=2$, for a moment we shall change the notation in this case, but in the future we shall use both notations in the case $d=2$, letters in bald mainly will be used in the cases where there is no essential difference in formulae between the cases $d=2$ and $d>2$. But at once we must say that for $d>2$ the complexity of asymptotic is growing with $d$ very rapidly, therefore we restricted ourselves with the case $d=2.$ Let  us denote
	\begin{equation}\label{field2}
	X_{k,l}=\sum_{i,j=0}^\infty c_{i,j}\e_{k-i, l-j}, \ (k, l)\in \bz^2,
	\end{equation}
	where $\e_{k, l}, \ (k, l)\in \bz^2,$ are i.i.d.  random $S\al S$  variables and $c_{i, j}, \ i\ge 1,  j\ge 0,$ are real numbers satisfying
	condition
	\begin{equation}\label{cond2a}
	A:=\sum_{i, j=0}^\infty |c_{i,j}|^\al <\infty.
	\end{equation}
	Considering $\rho_{n,m}:=\rho (X_{0,0}, X_{n,m})$ we found that the expression of $\rho_{n,m}$ is different in the cases $n>0, m>0$ and $n>0, m<0$, (due to the stationarity the remaining two cases can be transformed into the first ones). Then, trying to find the asymptotic of $\rho_{n,m}$, as $\min (m, n)\to \infty,$ we faced with much more complicated picture, comparing with the asymptotic for processes. For an arbitrary filter $\{c_{i,j}\}$ there is no hope to get asymptotic of $\rho_{n,m}$, therefore it is necessary to assume some regular behavior  of this filter. Usually one can consider two types of such behavior:
	$$
	c_{i,j} \sim \frac{1}{(i^{\bb_1}+j^{\bb_2})^{\bb_3}}, \quad {\rm  in \ \ particular,} \quad c_{i,j} \sim \frac{1}{||(i, j)||^{\bb_3}},
	$$
	with some positive $\bb_i, \ i=1,2,3,$ such, that (\ref{cond2a}) is satisfied, or
	\begin{equation}\label{cond3}
	c_{i, j}\sim (1+i)^{-\bb_1}(1+j)^{-\bb_2},
	\end{equation}
	with some $\bb_k>1/\al, \ k=1,2$. We had chosen the condition (\ref{cond3}), and this was done for two reasons: firstly, such behavior of the filter is easier to investigate, secondly, linear fields with such filter motivated to define the so-called directional memory (see \cite{Paul20}) and there is a hope to relate the behavior of spectral covariance with the memory properties. But even under condition (\ref{cond3}) the investigation of the asymptotic of $\rho_{n,m}$ was complicated and different in different areas, defined by various combinations of parameters $0<\al <2, \ \bb_k>1/\al, \ k=1,2.$

	Even the number of such different areas was unusually big in the case  $n>0, m<0$,   in some areas the asymptotic behavior of $\rho_{n,m}$ was dependent on the way how $m, n$ tends to infinity, more precisely, on the limit for  the ratio $m^{-\bb_2}/n^{-\bb_1}$.

	The rest of the paper is organized as follows. In the next section we formulate two theorems giving  the asymptotics of $\rho_{k,l}$, separately in the cases $k>0, l>0$ and $k>0, l<0$. In the last two sections proofs of these theorems are presented.

\section{Asymptotics of spectral covariance $\rho_{k,l}$}	

We consider linear random fields (\ref{field2}) satisfying (\ref{cond2a}).
For any $\ba=(a_1,a_2)\in\ZZ^2$  we will denote $Q_{\ba}=\left\lbrace \bx=(x_1,x_2)\in\ZZ^2: x_i\ge -a_i \right\rbrace$. Also, we will use the notations $\ba_+=((a_1)_+,(a_2)_+ )$ and $\ba_-=((a_1)_-,(a_2)_- )$, where $(\cdot)_+=\max(\cdot,0)$ and $(\cdot)_-=-\min(\cdot,0)$.
We are interested in the dependence between the r.v. $X_{\textbf{\BO}}$ and $X_{\bk}$, which are defined by \eqref{field2}. Due to the stationarity of the field, vectors $(X_{\BO}, X_{\bk})$ and $(X_{\BO+\bk_-}, X_{\bk+\bk_-})=(X_{\bk_-}, X_{\bk_+})$ have the same distribution.

If $\ba\in Q_\BO$ we have $Q_\BO\subset Q_\ba$,
\begin{equation*}
X_{\ba}=\sum_{\bi\in Q_\BO}c_\bi\e_{\ba-\bi}=
\sum_{\bi-\ba\in Q_\ba}c_{\bi-\ba+\ba}\e_{\ba-\bi}=
\sum_{\bj\in Q_\ba}c_{\bj+\ba}\e_{-\bj}=
\sum_{\bj\in Q_\BO}c_{\bj+\ba}\e_{-\bj}+
\sum_{\bj\in Q_\ba\backslash Q_\BO}c_{\bj+\ba}\e_{-\bj},
\end{equation*}
and we easily obtain the ch.f. of the vector $(X_{\bk_-}, X_{\bk_+})$:
\begin{equation*}
\EE\exp\left( {\rm i}(\theta_1X_{\bk_-}+\theta_2X_{\bk_+}) \right)
\end{equation*}
\begin{equation*}
=\exp\left(
-\abs{\theta_1}^\al
\sum_{\bj\in Q_{\bk_-}\backslash Q_\BO}\abs{c_{\bj+\bk_-}}^\al-
\sum_{\bj\in Q_\BO}\abs{ \theta_1c_{\bj+\bk_-}+\theta_2
	c_{\bj+\bk_+}}^\al
-\abs{\theta_2}^\al
\sum_{\bj\in Q_{\bk_+}\backslash Q_\BO}\abs{c_{\bj+\bk_+}}^\al
\right).
\end{equation*}
The obtained ch.f. allows us to see the structure of the spectral measure $\Gamma$ and obtain the expression for the spectral covariance:

\begin{equation}\label{eq:spectralCov}
\rho (X_{\BO}, X_{\bk})=\rho (X_{\bk_-}, X_{\bk_+})=\sum_{\bj\in Q_\BO}\frac{c_{\bj+\bk_-}c_{\bj+\bk_+}}{\left( c_{\bj+\bk_-}^2+c_{\bj+\bk_+}^2 \right)^{\frac{2-\al}{2}}}.
\end{equation}

As it was mentioned in the Introduction, we will consider one particular case - we will assume that the coefficients $c_{i,j}$ in \eqref{field2} have the form
\begin{equation}\label{eq:cijForm}
c_{i,j}=w_{(i,j)}(1+i)^{-\bb_1}(1+j)^{-\bb_2},\quad i,j\ge 0,
\end{equation}
where $\bb_k>\frac{1}{\al}, k=1,2$, and coefficients $w_{(i,j)}$ satisfy the following conditions:
\begin{enumerate}[label=(A{\arabic*})]
	\item\label{c:A1} there exists limit $\lim_{i,j\rightarrow\infty}w_{(i,j)}=1$,
	\item\label{c:A2} for every $i\geq 0$ there exists limit $w_{(i,\infty)}:=\lim_{j\rightarrow\infty}w_{(i,j)}>0$,
	\item\label{c:A3} for every $j\geq 0$ there exists limit $w_{(\infty,j)}:=\lim_{i\rightarrow\infty}w_{(i,j)}>0$.
\end{enumerate}
For convenience of writing we assume that $w_{(i,j)}\neq 0,\ (i,j)\in(\ZZ^+)^2$.

It follows from these assumptions that there exist numbers $d,e>0$ such that
\begin{equation}\label{ineq:wBounds}
d<|w_{(i,j)}|<e,\ \forall i,j\geq 0.
\end{equation}

We are interested in the asymptotic behavior of $\rho (k_1,k_2):=\rho (X_{(0,0)}, X_{(k_1,k_2)})$ as $k_1,k_2 \to \infty $ or $k_1\rightarrow\infty,\ k_2\rightarrow-\infty$.

It is convenient to introduce notation
 \begin{equation*}
 V(x,y)=\left\lbrace
 \begin{aligned}
 \frac{xy}{\left( x^2+y^2 \right)^{\frac{2-\al}{2}}}\ &\text{ if } x^2+y^2>0,\\
 0 &\text{ if } x=y=0.
 \end{aligned}\right.
 \end{equation*}
 Function $V$ has the following properties:
 \begin{itemize}
 	\item $V(x,y)=V(y,x)$,
 	\item $V(\lambda x,\lambda y)=\abs{\lambda}^\alpha V(x,y)$,
 	\item when $\al\in(0,2]$, the function $V$ is continuous,
 	\item The following inequalities hold
 	\begin{equation*}
 	\abs{V(x,y)}\leqslant 2^{\frac{\al-2}{2}}\abs{xy}^{\frac{\al}{2}},
 	\end{equation*}
 	\begin{equation*}
 	\abs{V(x,y)}\leqslant \abs{x}^{\al-1}\abs{y}.
 	\end{equation*}
 \end{itemize}
 This notation gives a more compact form of \eqref{eq:spectralCov}:
 \begin{equation}\label{eq:spectralCovA}
 \rho (X_{\BO}, X_{\bk})=\rho (X_{\bk_-}, X_{\bk_+})=
 \sum_{\bj\in Q_\BO}V\left( c_{\bj+\bk_-},c_{\bj+\bk_+} \right)
 .
 \end{equation}

Say $n,m>0$. From \eqref{eq:spectralCovA} we obtain the expressions of $\rho(n,m)$ and $\rho(n,-m)$:
\begin{equation}\label{eq:mainsum}
\rho(n,m)=\sum_{i=0}^{\infty}\sum_{j=0}^{\infty} V(w_{(i,j)}(1+i)^{-\bb_1}(1+j)^{-\bb_2}, \ w_{(i+n,j+m)}(1+i+n)^{-\bb_1}(1+j+m)^{-\bb_2}),
\end{equation}
\begin{equation}\label{spectcorneg}
\rho(n,-m)=\sum_{i=0}^{\infty}\sum_{j=0}^{\infty}V(w_{(i+n,j)}(1+i+n)^{-\bb_1}(1+j)^{-\bb_2}, \ w_{(i,j+m)}(1+i)^{-\bb_1}(1+j+m)^{-\bb_2}).
\end{equation}

 \begin{teo}\label{thm1}
 	Suppose that a linear field \eqref{field2} with coefficients $c_{i,j}$ having form \eqref{eq:cijForm}, satisfies conditions (A1)-(A3). Then the asymptotic behavior of spectral covariance $\rho(n,m)$ is as follows.
 	\begin{enumerate}
 		\item If $1<\al\leq 2$ and $\bb_i>\frac{1}{\al-1},\ i=1,2$,
 		\begin{equation}\label{asrel1}
 		\lim_{n,m\rightarrow\infty}\frac{\rho({n,m})}{n^{-\bb_1}m^{-\bb_2}}= \sum_{i=0}^{\infty}\sum_{j=0}^{\infty} w_{(i,j)}^{\langle\al-1\rangle}(1+i)^{-\bb_1(\al-1)}(1+j)^{-\bb_2(\al-1)}.
 		\end{equation}
 		\item If $1<\al\leq 2$ and $\frac{1}{\al}<\bb_i<\frac{1}{\al-1},\ i=1,2$ or $0<\al\leq 1$ and $\bb_i>\frac{1}{\al},\ i=1,2$,
 		\begin{equation}\label{asrel2}
 		\lim_{n,m\rightarrow\infty}
 		\frac{\rho(n,m)}{n^{1-\al\bb_1}m^{1-\al\bb_2}}= \int_{0}^{\infty}\int_{0}^{\infty} V(t^{-\bb_1}s^{-\bb_2}, \ (t+1)^{-\bb_1}(s+1)^{-\bb_2})\dt \ds.
 		\end{equation}
 		\item If $1<\al\leq 2$ and $\bb_1>\frac{1}{\al-1},\ \frac{1}{\al}<\bb_2<\frac{1}{\al-1},$
 		\begin{equation*}
 		\lim_{n,m\rightarrow\infty}\frac{\rho(n,m)}{n^{-\bb_1}m^{1-\bb_2\al}}= \sum_{i=0}^{\infty}w_{(i,\infty)}^{\al-1}(1+i)^{-\bb_1(\al-1)}\int\limits_{0}^{\infty} u^{-\bb_2(\al-1)}(1+u)^{-\bb_2}\dd u.
 		\end{equation*}
 		\item If $1<\al\leq 2$ and $\bb_1=\frac{1}{\al-1},\ \frac{1}{\al}<\bb_2<\frac{1}{\al-1},$
 		\begin{equation*}
 		\lim_{n,m\rightarrow\infty}\frac{\rho(n,m)}{n^{-\bb_1}m^{1-\bb_2\al}\ln (n)}= \int_{0}^{\infty}
 		v^{-\bb_2(\al-1)}
 		\left(1+v\right)^{-\bb_2}  \dd v.
 		\end{equation*}
 		\item If $1<\al\leq 2$ and $\bb_1=\frac{1}{\al-1},\ \bb_2>\frac{1}{\al-1},$
 		\begin{equation*}
 		\lim_{n,m\rightarrow\infty}\frac{\rho(n,m)}{n^{-\bb_1}m^{-\bb_2}\ln (n)} \sum_{j=0}^{\infty}w_{(\infty,j)}^{\al-1}(1+j)^{-\bb_2(\al-1)}.
 		\end{equation*}
 		\item If $1<\al\leq 2$ and $\bb_1=\frac{1}{\al-1},\ \bb_2=\frac{1}{\al-1},$
 		\begin{equation*}
 		\lim_{n,m\rightarrow\infty}\frac{\rho(n,m)}{n^{-\bb_1}m^{-\bb_2}\ln (n)\ln (m)}=1.
 		\end{equation*}
 	\end{enumerate}
 	
 \end{teo}

\begin{teo}\label{thm2}
	Suppose that a linear field \eqref{field2} with coefficients $c_{i,j}$ having form \eqref{eq:cijForm}, satisfies conditions (A1)-(A3). Then the asymptotic behavior of spectral covariance $\rho(n,-m)$ is as follows.
	\begin{enumerate}
		\item If $1<\al\le 2,\ \frac{1}{\bb_1}+\frac{1}{\bb_2}>\al$ and $\frac{1}{\al}<\bb_i<\frac{1}{\al-1 },\ i=1,2,$ or $0<\al\le 1$ and $\frac{1}{\bb_1}+\frac{1}{\bb_2}>\al$:
		\begin{equation*}
		\lim_{n,m\rightarrow\infty}\frac{\rho(n,-m)}{n^{1-\bb_1\al}m^{1-\bb_2\al}}=
		\int\limits_{0}^{\infty}
		\int\limits_{0}^{\infty} V((1+u)^{-\bb_1}v^{-\bb_2}, \ u^{-\bb_1}(1+v)^{-\bb_2} )\dd u\dd v.
		\end{equation*}
		\item If $1<\al\le 2,\ \frac{1}{\al-1}<\bb_1$ and $\frac{1}{\al}<\bb_2<1$:
		\begin{equation*}
		\lim_{n,m\rightarrow\infty}\frac{\rho(n,-m)}{n^{-\bb_1}m^{1-\bb_2\al}}=
		\sum_{i=0}^{\infty}
		w_{(i,\infty)}^{\al-1}
		(1+i)^{-\bb_1(\al-1)}
		\int\limits_{0}^{\infty}v^{-\bb_2}(1+v)^{-\bb_2(\al-1)}\dd v.
		\end{equation*}
		\item If $ 1<\al\le 2,\ \bb_1=\frac{1}{\al-1}$ and $\frac{1}{\al}<\bb_2<1 $:
		\begin{equation*}
		\lim_{n,m\rightarrow\infty}\frac{\rho(n,-m)}{n^{1-\bb_1\al}m^{1-\bb_2\al}\ln(n)}
		\int_{0}^{\infty}v^{-\bb_2}(1+v)^{-\bb_2(\al-1)}\dd v
		\ds.
		\end{equation*}
		\item If  $1<\al\le 2$ and $\frac{1}{\al-1}<\bb_i,\ i=1,2$:
		\begin{enumerate}
			\item If $m_n$ is a sequence such that $\frac{m_n^{-\bb_2}}{n^{-\bb_1}}\rightarrow 0$:
			\begin{equation*}
			\lim_{n\rightarrow\infty}\frac{\rho(n,-m_n)}{n^{-\bb_1(\al-1)}m_n^{-\bb_2}}=
			\sum_{i=0}^{\infty}\sum_{j=0}^{\infty}
			w_{(\infty,j)}^{\al-1}w_{(i,\infty)}(1+j)^{-\bb_2(\al-1)}(1+i)^{-\bb_1}
			.
			\end{equation*}
			\item If $m_n$ is a sequence such that $\frac{m_n^{-\bb_2}}{n^{-\bb_1}}\rightarrow c\in(0;\infty)$:
			\begin{equation*}
			\lim_{n\rightarrow\infty}\frac{\rho(n,-m_n)}{n^{-\bb_1(\al-1)}m_n^{-\bb_2}}=\frac{1}{c}
			\sum_{i=0}^{\infty}\sum_{j=0}^{\infty} V(w_{(\infty,j)}(1+j)^{-\bb_2}, \ cw_{(i,\infty)}(1+i)^{-\bb_1}).
			\end{equation*}
			\item If $m_n\rightarrow\infty$ is a sequence such that $\frac{m_n^{-\bb_2}}{n^{-\bb_1}}\rightarrow \infty$:
			\begin{equation*}
			\lim_{n\rightarrow\infty}\frac{\rho(n,-m_n)}{n^{-\bb_1}m_n^{-\bb_2(\al-1)}}=
			\sum_{i=0}^{\infty}\sum_{j=0}^{\infty}{w_{(\infty,j)}w_{(i,\infty)}^{\al-1}(1+j)^{-\bb_2}(1+i)^{-\bb_1(\al-1)}}.
			\end{equation*}
		\end{enumerate}
		\item If $1<\al< 2,\ \frac{1}{\al-1}<\bb_1$ and $1<\bb_2<\frac{1}{\al-1}$:
		\begin{enumerate}
			\item If $m_n$ is a sequence such that $\frac{m_n^{-\bb_2}}{n^{-\bb_1}}\rightarrow 0 $:
			\begin{equation*}
			\lim_{n\rightarrow\infty}\frac{\rho(n,-m_n)}{n^{-\frac{\bb_1}{\bb_2}}m_n^{1-\bb_2\al}}=
			\sum_{i=0}^{\infty}\int_{0}^{\infty}V(v^{-\bb_2}, \ w_{(i,\infty)}(1+i)^{-\bb_1})\dd v.
			\end{equation*}
			\item If $m_n$ is a sequence such that $\frac{m_n^{-\bb_2}}{n^{-\bb_1}}\rightarrow c\in(0;\infty) $:
			\begin{equation*}
			\lim_{n\rightarrow\infty}\frac{\rho(n,-m_n)}{n^{-\bb_1}m_n^{-\bb_2(\al-1)}}=
			c\sum_{i=0}^{\infty}\sum_{j=0}^{\infty} V(w_{(\infty,j)}c^{-1}(1+j)^{-\bb_2}, \ w_{(i,\infty)}(1+i)^{-\bb_1}).
			\end{equation*}
			\item If $m_n\rightarrow\infty$ is a sequence such that $\frac{m_n^{-\bb_2}}{n^{-\bb_1}}\rightarrow \infty $:
			\begin{equation*}
			\lim_{n\rightarrow\infty}\frac{\rho(n,-m_n)}{n^{-\bb_1}m_n^{-\bb_2(\al-1)}}=
			\sum_{i=0}^{\infty}\sum_{j=0}^{\infty}
			w_{(\infty,j)}w_{(i,\infty)}^{\al-1}(1+j)^{-\bb_2}(1+i)^{-\bb_1(\al-1)}
			.
			\end{equation*}
		\end{enumerate}
		\item If $1<\al< 2,\ 1<\bb_i<\frac{1}{\al-1},\ i=1,2$ and $ \frac{1}{\bb_1}+\frac{1}{\bb_2}<\al$ or $0<\al\le 1$ and $\frac{1}{\bb_1}+\frac{1}{\bb_2}<\al$:
		\begin{enumerate}
			\item If $m_n$ is a sequence such that $\frac{m_n^{-\bb_2}}{n^{-\bb_1}}\rightarrow 0$:
			\begin{equation*}
			\lim_{n\rightarrow\infty}\frac{\rho(n,-m_n)}{n^{-\frac{\bb_1}{\bb_2}}m_n^{1-\bb_2\al}}=
			\sum_{i=0}^{\infty}\int_{0}^{\infty} V(v^{-\bb_2}, \ w_{(i,\infty)}(1+i)^{-\bb_1})\dd v.
			\end{equation*}
			\item If $m_n$ is a sequence such that $\frac{m_n^{-\bb_2}}{n^{-\bb_1}}\rightarrow c\in(0;\infty) $:
			\begin{equation*}
			\lim_{n\rightarrow\infty}\frac{\rho(n,-m_n)}{n^{-\bb_1(\al-1)}m_n^{-\bb_2}}=\frac{1}{c}
			\sum_{i=0}^{\infty}\sum_{j=0}^{\infty}V(w_{(\infty,j)}(1+j)^{-\bb_2}, \ w_{(i,\infty)}c^2(1+i)^{-\bb_1}).
			\end{equation*}
			\item If $m_n\rightarrow\infty$ is a sequence such that $\frac{m_n^{-\bb_2}}{n^{-\bb_1}}\rightarrow \infty$:
			\begin{equation*}
			\lim_{n\rightarrow\infty}\frac{\rho(n,-m_n)}{n^{1-\bb_1\al}m_n^{-\frac{\bb_2}{\bb_1}}}=
			\sum_{j=0}^{\infty}\int_{0}^{\infty}V(u^{-\bb_1}, \ w_{(\infty,j)}(1+j)^{-\bb_2})\dd u.
			\end{equation*}
		\end{enumerate}
		\item If $1<\al< 2,\ \frac{1}{\al-1}<\bb_1$ and $\bb_2=\frac{1}{\al-1}$:
		\begin{enumerate}
			\item If $m_n$ is a sequence such that $\frac{m_n^{-\bb_2}}{n^{-{\bb_1}}}\rightarrow 0$:
			\begin{equation*}
			\lim_{n\rightarrow\infty}\frac{\rho(n,-m_n)}{n^{-\frac{\bb_1}{\bb_2}}m_n^{1-\bb_2\al}\ln \left(  n^{-\frac{\bb_1}{\bb_2}}m_n \right)}=
			\sum_{i=0}^{\infty}w_{(i,\infty)}(1+i)^{-\bb_1}.
			\end{equation*}
			\item If $m_n$ is a sequence such that $\frac{m_n^{-\bb_2}}{n^{-{\bb_1}}}\rightarrow c\in(0;\infty)$:
			\begin{equation*}
			\lim_{n\rightarrow\infty}\frac{\rho(n,-m_n)}{n^{-\bb_1}m_n^{-\bb_2(\al-1)}}=
			c\sum_{i=0}^{\infty}\sum_{j=0}^{\infty}V(w_{(\infty,j)}c^{-1}(1+j)^{-\bb_2}, \ w_{(i,\infty)}(1+i)^{-\bb_1}).
			\end{equation*}
			\item If $m_n\rightarrow\infty$ is a sequence such that $\frac{m_n^{-\bb_2}}{n^{-{\bb_1}}}\rightarrow \infty$:
			\begin{equation*}
			\lim_{n\rightarrow\infty}\frac{\rho(n,-m_n)}{n^{-\bb_1}m_n^{-\bb_2(\al-1)}}=
			\sum_{i=0}^{\infty}\sum_{j=0}^{\infty}
			w_{(\infty,j)}w_{(i,\infty)}^{\al-1}(1+j)^{-\bb_2}(1+i)^{-\bb_1(\al-1)}.
			\end{equation*}
			\end{enumerate}
		\item If $1<\al< 2,\ \bb_1=\frac{1}{\al-1}$ and $1<\bb_2<\frac{1}{\al-1} $
		\begin{enumerate}
			\item If $m_n$ is a sequence such that $\frac{m_n^{-\bb_2}}{n^{-\bb_1}}\rightarrow 0$:
			\begin{equation*}
			\lim_{n\rightarrow\infty}\frac{\rho(n,-m_n)}{n^{-\frac{\bb_1}{\bb_2}}m_n^{1-\bb_2\al}}=
			\sum_{i=0}^{\infty}\int_{0}^{\infty}V(w_{(i,\infty)}\left(1+i\right)^{-\bb_1}, \ s^{-\bb_2})\dd s.
			\end{equation*}
			\item If $m_n$ is a sequence such that $\frac{m_n^{-\bb_2}}{n^{-\bb_1}}\rightarrow c\in(0;\infty) $:
			\begin{equation*}
			\lim_{n\rightarrow\infty}\frac{\rho(n,-m_n)}{n^{-\bb_1}m_n^{-\bb_2(\al-1)}}=
			c\sum_{i=0}^{\infty}\sum_{j=0}^{\infty}V(w_{(\infty,j)}c^{-1}(1+j)^{-\bb_2}, \ w_{(i,\infty)}(1+i)^{-\bb_1}).
			\end{equation*}
			\item If $m_n\rightarrow\infty$ is a sequence such that $\frac{m_n^{-\bb_2}}{n^{-\bb_1}}\rightarrow \infty $:
			\begin{equation*}
			\lim_{n\rightarrow\infty}\frac{\rho(n,-m_n)}{n^{1-\bb_1\al}m_n^{-\frac{\bb_2}{\bb_1}}\ln\left(nm_n^{-\frac{\bb_2}{\bb_1}}\right)}=
			\sum_{j=0}^{\infty}w_{(\infty,j)}\left(1+j\right)^{-\bb_2}.
			\end{equation*}
		\end{enumerate}
		\item If $1<\al<2$ and $\bb_i=\frac{1}{\al-1},\ i=1,2$:
		\begin{enumerate}
			\item If $m_n$ is a sequence such that $\frac{m_n^{-\bb_2}}{n^{-\bb_1}}\rightarrow 0 $:
			\begin{equation*}
			\lim_{n\rightarrow\infty}\frac{\rho(n,-m_n)}{n^{-\bb_1/\bb_2}m_n^{1-\bb_2\al}\ln\left( n^{-\bb_1/\bb_2}m_n \right)}=
			\sum_{i=0}^{\infty}w_{(i,\infty)}(1+i)^{-\bb_1}.
			\end{equation*}
			\item If $m_n$ is a sequence such that $\frac{m_n^{-\bb_2}}{n^{-\bb_1}}\rightarrow c\in(0;\infty) $:
			\begin{equation*}
			\lim_{n\rightarrow\infty}\frac{\rho(n,-m_n)}{n^{-\bb_1}m_n^{-\bb_2(\al-1)}}=
			c\sum_{i=0}^{\infty}\sum_{j=0}^{\infty}V(w_{(\infty,j)}c^{-1}(1+j)^{-\bb_2}, \ w_{(i,\infty)}(1+i)^{-\bb_1}).
			\end{equation*}
			\item If $m_n\rightarrow\infty$ is a sequence such that $\frac{m_n^{-\bb_2}}{n^{-\bb_1}}\rightarrow \infty $:
			\begin{equation*}
			\lim_{n\rightarrow\infty}\frac{\rho(n,-m_n)}{n^{1-\bb_1\al}m_n^{-\bb_2/\bb_1}\ln\left( nm_n^{-\bb_2/\bb_1} \right)}=
			\sum_{j=0}^{\infty}w_{(\infty,j)}(1+j)^{-\bb_2}.
			\end{equation*}
		\end{enumerate}
	\end{enumerate}
\end{teo}

\section{Proof of Theorem \ref{thm1}}

%
%
%
%

In order not to make the paper too long, we will not provide all of the proofs. We feel that the provided proofs illustrate the ideas used, the other sets of parameters are dealt with in a similar way.
Here we will investigate the following sets of parameters:
\begin{enumerate}
	\item $1<\al\leq 2,\ \bb_i>\frac{1}{\al-1},\ i=1,2;$
	\item $1<\al\leq 2,\ \frac{1}{\al}<\bb_i<\frac{1}{\al-1},\ i=1,2$;
	\item $1<\al\leq 2,\ \bb_1>\frac{1}{\al-1},\ \frac{1}{\al}<\bb_2<\frac{1}{\al-1};$
	\item $1<\al\leq 2,\ \bb_1=\frac{1}{\al-1},\ \bb_2=\frac{1}{\al-1};$
\end{enumerate}


\subsection{The case $1<\al\leq 2,\ \bb_i>\frac{1}{\al-1},\ i=1,2$}
In this case we have
\begin{equation*}
\sum_{i=0}^{\infty}\sum_{j=0}^{\infty} |c_{i,j}|^{\al-1}=
\sum_{i=0}^{\infty}\sum_{j=0}^{\infty} |w_{(i,j)}(1+i)^{-\bb_1}(1+j)^{-\bb_2}|^{\al-1}
\end{equation*}
\begin{equation*}
\leq
e^{\al-1}
\sum_{i=0}^{\infty}(1+i)^{-\bb_1(\al-1)}
\sum_{j=0}^{\infty}(1+j)^{-\bb_2(\al-1)}<\infty.
\end{equation*}
We will show that
\begin{equation}\label{eq:rez1}
\frac{\rho({n,m})}{n^{-\bb_1}m^{-\bb_2}}\rightarrow \sum_{i=0}^{\infty}\sum_{j=0}^{\infty} c_{i,j}^{<\al-1>},\ n,m\rightarrow\infty,
\end{equation}
here $b^{<a>}=|b|^a{\rm sgn}(b)$.

 Fix $\e>0$ and choose $N_1\in\NN$ such that
 \begin{equation*}
 \mathop{\sum\sum}_{\max(i,j)> N_1} |c_{i,j}|^{\al-1}<\e.
 \end{equation*}
 Such a number exists since $\sum_{i=0}^{\infty}\sum_{j=0}^{\infty} |c_{i,j}|^{\al-1}<\infty$.

Let us consider
 \begin{equation*}
 \Delta_{n,m}:=\left| \frac{\rho(n,m)}{n^{-\bb_1}m^{-\bb_2}}-\sum_{i=0}^{\infty}\sum_{j=0}^{\infty} c_{i,j}^{<\al-1>} \right|.
 \end{equation*}

 From the triangle inequality we obtain
 \begin{equation*}
  \Delta_{n,m}\leq S_1+S_2+S_3,
 \end{equation*}
where
 \begin{equation*}
 S_1=\left|
 \mathop{\sum\sum}_{i,j\leq N_1}\frac{c_{i,j}c_{i+n,j+m}n^{\bb_1}m^{\bb_2}}{\left( c_{i,j}^2+c_{i+n,j+m}^2 \right)^{\frac{2-\al}{2}}}  -\mathop{\sum\sum}_{i,j\leq N_1} c_{i,j}^{<\al-1>}
 \right|,
 \end{equation*}
 \begin{equation*}
 S_2=\left|
 \mathop{\sum\sum}_{\max(i,j)> N_1}\frac{c_{i,j}c_{i+n,j+m}n^{\bb_1}m^{\bb_2}}{\left( c_{i,j}^2+c_{i+n,j+m}^2 \right)^{\frac{2-\al}{2}}}   \right|,
 \end{equation*}
 \begin{equation*}
 S_3= \left|
 \mathop{\sum\sum}_{\max(i,j)> N_1} c_{i,j}^{<\al-1>}
 \right|.
 \end{equation*}

 The choice of $N_1$ implies
 \begin{equation}\label{ineq:21S3}
 S_3<\e.
 \end{equation}
  We also have
 \begin{equation}\label{ineq:21S2}\begin{gathered}
 S_2\leq \mathop{\sum\sum}_{\max(i,j)> N_1}|c_{i,j}|^{\al-1}|c_{i+n,j+m}|n^{\bb_1}m^{\bb_2} \\
 \leq
 e \mathop{\sum\sum}_{\max(i,j)> N_1}|c_{i,j}|^{\al-1}\left( \frac{n}{1+i+n} \right)^{\bb_1}\left( \frac{m}{1+j+m} \right)^{\bb_2}\leq e\mathop{\sum\sum}_{\max(i,j)> N_1}|c_{i,j}|^{\al-1}\leq e\e.
 \end{gathered}\end{equation}
 For fixed values of $i,j$ we have $c_{i+n,j+m}\rightarrow 0$ and $c_{i+n,j+m}n^{\bb_1}m^{\bb_2}\rightarrow 1$, as $n,m\rightarrow \infty$, thus
 \begin{equation*}
 \frac{c_{i,j}c_{i+n,j+m}n^{\bb_1}m^{\bb_2}}{\left( c_{i,j}^2+c_{i+n,j+m}^2 \right)^{\frac{2-\al}{2}}} \rightarrow c_{i,j}^{<\al-1>},\ \text{ as } n,m\rightarrow\infty.
 \end{equation*}
 Since
 \begin{equation*}
 S_1\leq  \mathop{\sum\sum}_{i,j\leq N_1} \left| \frac{c_{i,j}c_{i+n,j+m}n^{\bb_1}m^{\bb_2}}{\left( c_{i,j}^2+c_{i+n,j+m}^2 \right)^{\frac{2-\al}{2}}}  -c_{i,j}^{<\al-1>}\right|,
 \end{equation*}
 and the sum on the right side is finite, we obtain $\lim_{n,m\rightarrow\infty}S_1=0$. Together with \eqref{ineq:21S3} and \eqref{ineq:21S2} this implies $\limsup_{n,m\rightarrow\infty}\Delta_{n,m}<(1+e)\e$. Since $\e>0$ is arbitrary, we obtain $\lim_{n,m\rightarrow\infty}\Delta_{n,m}=0$.


\subsection{The case $1<\al\leq 2,\ \frac{1}{\al}<\bb_i<\frac{1}{\al-1},\ i=1,2$}
 We will show that
 \begin{equation}\label{IIrez}
 \frac{\rho(n,m)}{n^{1-\al\bb_1}m^{1-\al\bb_2}}\rightarrow \int_{0}^{\infty}\int_{0}^{\infty}\frac{t^{-\bb_1}(t+1)^{-\bb_1}s^{-\bb_2}(s+1)^{-\bb_2}}{\left( t^{-2\bb_1}s^{-2\bb_2}+(t+1)^{-2\bb_1}(s+1)^{-2\bb_2} \right)^{\frac{2-\al}{2}}}\dt \ds.
 \end{equation}

  Denote $\bb=\max(\bb_1,\bb_2)$ and define $\g=\frac{1}{2}\left( 1-\bb(\al-1) \right)$. We split \eqref{eq:mainsum} into four terms
 \begin{equation}\label{eq:rhoSplit}
 \rho({n,m})= I_{1,1}+I_{1,2}+I_{2,1}+I_{2,2},
 \end{equation}
 where
 \begin{equation*}\begin{gathered}
 I_{1,1}= \sum_{i\leq n^\g}\sum_{j\leq m^\g}V\left( c_{i,j}, c_{i+n,j+m} \right), \quad I_{1,2}= \sum_{i\leq n^\g}\sum_{j> m^\g}V\left( c_{i,j}, c_{i+n,j+m} \right),\\
 I_{2,1}= \sum_{i> n^\g}\sum_{j\leq m^\g}V\left( c_{i,j}, c_{i+n,j+m} \right),\quad I_{2,2}= \sum_{i> n^\g}\sum_{j> m^\g}V\left( c_{i,j}, c_{i+n,j+m} \right).
 \end{gathered}
 \end{equation*}

 We start by estimating $I_{1,1}$.
 \begin{equation}\label{bound:2I11}
 \frac{|I_{1,1}|}{n^{1-\al\bb_1}m^{1-\al\bb_2}}\leq
 n^{\al\bb_1-1}m^{\al\bb_2-1}\sum_{i\leq n^\g}\sum_{j\leq m^\g} |c_{i,j}|^{\al-1}|c_{i+n,j+m}|.
 \end{equation}
 Inequalities \eqref{ineq:wBounds} imply that $|c_{i,j}|<e$. Thus, from \eqref{bound:2I11} we obtain
 \begin{equation*}
 \frac{|I_{1,1}|}{n^{1-\al\bb_1}m^{1-\al\bb_2}}\leq
  e^{\al-1}n^{\al\bb_1-1}m^{\al\bb_2-1} \sum_{i\leq n^\g}\sum_{j\leq m^\g}|c_{i+n,j+m}|
 \end{equation*}
 \begin{equation*}
 \leq  e^{\al}n^{\al\bb_1-1}m^{\al\bb_2-1}\sum_{i\leq n^\g}\sum_{j\leq m^\g}(1+i+n)^{-\bb_1}(1+j+m)^{-\bb_2}
 \end{equation*}
 \begin{equation*}
 \leq
 e^{\al}n^{\al\bb_1-1}m^{\al\bb_2-1}\sum_{i\leq n^\g}\sum_{j\leq m^\g}n^{-\bb_1}m^{-\bb_2}
 \end{equation*}
 \begin{equation*}
 \leq    e^{\al} n^{\al\bb_1-1}m^{\al\bb_2-1} n^\g n^{-\bb_1}m^\g m^{-\bb_2}
 \end{equation*}
 \begin{equation*}
 =
 e^{\al} n^{\bb_1(\al-1)-1+\g}m^{\bb_2(\al-1)-1+\g}.
 \end{equation*}

 The choice of $\gamma$ implies $\frac{|I_{1,1}|}{n^{1-\al\bb_1}m^{1-\al\bb_2}}\rightarrow 0$ as $n,m\rightarrow\infty$.
Let us investigate $I_{1,2}$ and $I_{2,1}$ next. Due to the symmetry it is enough to deal with  one of these terms.
 \begin{equation*}\begin{gathered}
 {|I_{2,1}|}\leq
 \sum_{i> n^\g}\sum_{j\leq m^\g}  |c_{i,j}|^{\al-1} |c_{i+n,j+m}| \\
 \leq
 e^\al\sum_{i> n^\g}\sum_{j\leq m^\g}
 (1+i)^{-\bb_1(\al-1)}(1+j)^{-\bb_2(\al-1)} (1+n+i)^{-\bb_1}(1+m+j)^{-\bb_2}\\
 =e^\al
 \sum_{i> n^\g}
 (1+i)^{-\bb_1(\al-1)} (1+n+i)^{-\bb_1}
 \sum_{j\leq m^\g}
 (1+j)^{-\bb_2(\al-1)}(1+m+j)^{-\bb_2}\\
 \leq e^\al m^{\g-\bb_2}
 \sum_{i> n^\g}
 (1+i)^{-\bb_1(\al-1)} (1+n+i)^{-\bb_1}\\
 \leq e^\al m^{\g-\bb_2} \int\limits_{n^\g}^{\infty} (1+\lfloor t\rfloor)^{-\bb_1(\al-1)} (1+n+\lfloor t\rfloor)^{-\bb_1}\dt
 \end{gathered}
 \end{equation*}
  \begin{equation*}
  \leq e^\al m^{\g-\bb_2} \int\limits_{n^\g}^{\infty}  t^{-\bb_1(\al-1)} (t+n)^{-\bb_1}\dt
  \end{equation*}
  \begin{equation*}
  = e^\al m^{\g-\bb_2} \int\limits_{n^{\g-1}}^{\infty} (nt )^{-\bb_1(\al-1)} (nt+n)^{-\bb_1}\dd nt
  \end{equation*}
 \begin{equation*}
 = e^\al m^{\g-\bb_2}n^{1-\bb_1\al} \int\limits_{n^{\g-1}}^{\infty} t^{-\bb_1(\al-1)} (t+1)^{-\bb_1}\dd t=:
 e^\al m^{\g-\bb_2}n^{1-\bb_1\al}B(\al,\bb_1).
 \end{equation*}
  $B(\al,\bb_1)$ is finite, since $\frac{1}{\al}<\bb_i<\frac{1}{\al-1},\ i=1,2$ in the case under consideration.
 We obtain
 \begin{equation*}
 \frac{|I_{2,1}|}{n^{1-\al\bb_1}m^{1-\al\bb_2}}\leq
 e^\al m^{\g+\bb_2(\al-1)-1}B(\al,\bb_1) \rightarrow 0,\text{ as } n,m\rightarrow\infty,
 \end{equation*}
 and symmetry implies
 \begin{equation*}
 \frac{|I_{1,2}|}{n^{1-\al\bb_1}m^{1-\al\bb_2}}\rightarrow 0,\text{ as } n,m\rightarrow\infty.
 \end{equation*}

 It remains to deal with $I_{2,2}$.
 \begin{equation*}\begin{gathered}
 \frac{I_{2,2}}{n^{1-\al\bb_1}m^{1-\al\bb_2}}=\sum_{i> n^\g}\sum_{j> m^\g}
 \frac{c_{i,j}c_{i+n,j+m}n^{\al\bb_1-1}m^{\al\bb_2-1}}
 {\left( c_{i,j}^2+c_{i+n,j+m}^2 \right)^{\frac{2-\al}{2}}}\\
 \end{gathered} \end{equation*}

 Choose arbitrary $0<\tau<1$. There exists $N_2\in\NN$, such that
  \begin{equation*}
 (1-\tau)(1+i)^{-\bb_1}(1+j)^{-\bb_2}\leq c_{i,j}\leq(1+\tau)(1+i)^{-\bb_1}(1+j)^{-\bb_2},\ i,j>N_2.
 \end{equation*}
 These inequalities imply that for large $n,m$ (large enough, that the inequalities $n^\g>N_2$ and  $m^\g>N_2$ hold)
 \begin{equation}\label{I22Anm}
 \frac{(1-\tau)^2}{(1+\tau)^{2-\al}}A_{n,m}\leq
 I_{2,2}
 \leq  \frac{(1+\tau)^2}{(1-\tau)^{2-\al}}A_{n,m},
 \end{equation}
 where
 \begin{equation*}
 A_{n,m}=\sum_{i> n^\g}\sum_{j> m^\g}
 \frac{(1+i)^{-\bb_1}(1+j)^{-\bb_2}(1+i+n)^{-\bb_1}(1+j+m)^{-\bb_2}}
 {\left( (1+i)^{-2\bb_1}(1+j)^{-2\bb_2}+(1+i+n)^{-2\bb_1}(1+j+m)^{-2\bb_2} \right)^{\frac{2-\al}{2}}}.
 \end{equation*}

 We will investigate asymptotic behavior of $A_{n,m}$ as $n,m\rightarrow\infty$.
 {\small
 	\begin{equation*}\begin{gathered}
 	A_{n,m}=\int\limits_{\lfloor n^\g\rfloor+1}^{\infty}
 	\int\limits_{\lfloor m^\g\rfloor+1}^{\infty}
 	V\left((1+\lfloor t\rfloor)^{-\bb_1}(1+\lfloor s\rfloor)^{-\bb_2},(1+\lfloor t\rfloor+n)^{-\bb_1}(1+\lfloor s\rfloor+m)^{-\bb_2}\right)\dt\ds\\
 	=
 	\int\limits_{\frac{\lfloor n^\g\rfloor+1}{n}}^{\infty}
 	\int\limits_{\frac{\lfloor m^\g\rfloor+1}{m}}^{\infty}
 	V\left((1+\lfloor nt\rfloor)^{-\bb_1}(1+\lfloor ms\rfloor)^{-\bb_2},(1+\lfloor nt\rfloor+n)^{-\bb_1}(1+\lfloor ms\rfloor+m)^{-\bb_2}\right)\dt n\ds m \\
 	=n^{1-\al\bb_1}m^{1-\al\bb_2}\times\\\times
 	\int\limits_{\frac{\lfloor n^\g\rfloor+1}{n}}^{\infty}
 	\int\limits_{\frac{\lfloor m^\g\rfloor+1}{m}}^{\infty}
 	V\left(\left(\frac{1+\lfloor nt\rfloor}{n}\right)^{-\bb_1}\left(\frac{1+\lfloor ms\rfloor}{m}\right)^{-\bb_2},\left(\frac{1+\lfloor nt\rfloor}{n}+1\right)^{-\bb_1}\left(\frac{1+\lfloor ms\rfloor}{m}+1\right)^{-\bb_2}\right)\dt\ds .
 	\end{gathered}   \end{equation*}}
 Let us denote the function under integral as $V_{n,m,t,s}$, then we have
 \begin{equation}\label{AnmIntegr}
 \frac{A_{n,m}}{n^{1-\al\bb_1}m^{1-\al\bb_2}}=
 \int\limits_{\frac{\lfloor n^\g\rfloor+1}{n}}^{\infty}
 \int\limits_{\frac{\lfloor m^\g\rfloor+1}{m}}^{\infty}V_{n,m,t,s}\dt\ds
  \end{equation}

 It is easy to see that $\w_{n,m,t,s}$ converges pointwise to $\w\left(t^{-\bb_1}s^{-\bb_2},\left(t+1\right)^{-\bb_1}\left(s+1\right)^{-\bb_2}\right)$, $(s,t)\in(0,\infty)^2$.
Since $1+\lfloor nt \rfloor>nt$ and $1+\lfloor ms \rfloor>ms$, it is easy to see that
 \begin{equation}\label{wrezis}
 \w_{n,m,t,s}\leq t^{-\bb_1(\al-1)}s^{-\bb_2(\al-1)}\left(t+1\right)^{-\bb_1}\left(s+1\right)^{-\bb_2},\ t,s\in(0,\infty).
 \end{equation}
 The function $t^{-\bb_1(\al-1)}s^{-\bb_2(\al-1)}\left(t+1\right)^{-\bb_1}\left(s+1\right)^{-\bb_2}$ is integrable in $(0,\infty)^2$ since $1<\al\leq 2,\ \frac{1}{\al}<\bb_i<\frac{1}{\al-1},\ i=1,2,$
 therefore, applying the Dominated Convergence Theorem, we obtain
%
\begin{equation*}
\frac{A_{n,m}}{n^{1-\al\bb_1}m^{1-\al\bb_2}}\rightarrow A:= \int_{0}^{\infty}\int_{0}^{\infty}\frac{t^{-\bb_1}(t+1)^{-\bb_1}s^{-\bb_2}(s+1)^{-\bb_2}}{\left( t^{-2\bb_1}s^{-2\bb_2}+(t+1)^{-2\bb_1}(s+1)^{-2\bb_2} \right)^{\frac{2-\al}{2}}}\dt \ds.
\end{equation*}
\eqref{I22Anm} implies
\begin{equation*}
\frac{(1-\tau)^2}{(1+\tau)^{2-\al}}A\leq
\liminf_{n,m\rightarrow\infty}\frac{I_{2,2}}{n^{1-\al\bb_1}m^{1-\al\bb_2}}
\leq
\limsup_{n,m\rightarrow\infty}
\frac{I_{2,2}}{n^{1-\al\bb_1}m^{1-\al\bb_2}}
\leq  \frac{(1+\tau)^2}{(1-\tau)^{2-\al}}A
\end{equation*}
for arbitrary $\tau\in(0,1)$. Passing to the limit as $\tau\rightarrow0$, we get
\begin{equation*}
\lim_{n,m\rightarrow\infty}\frac{I_{2,2}}{n^{1-\al\bb_1}m^{1-\al\bb_2}}=A.
\end{equation*}

In conclusion
 \begin{equation*}
 \frac{\rho(n,m)}{n^{1-\al\bb_1}m^{1-\al\bb_2}}
 \end{equation*}
 \begin{equation*}=
 \frac{I_{1,1}+I_{1,2}+I_{2,1}+I_{2,2}}{n^{1-\al\bb_1}m^{1-\al\bb_2}}\rightarrow \int_{0}^{\infty}\int_{0}^{\infty}\frac{t^{-\bb_1}(t+1)^{-\bb_1}s^{-\bb_2}(s+1)^{-\bb_2}}{\left( t^{-2\bb_1}s^{-2\bb_2}+(t+1)^{-2\bb_1}(s+1)^{-2\bb_2} \right)^{\frac{2-\al}{2}}}\dt \ds.
 \end{equation*}

    \subsection{The case $1<\al\leq 2,\ \bb_1>\frac{1}{\al-1},\ \frac{1}{\al}<\bb_2<\frac{1}{\al-1}$;}
We will show that
\begin{equation*}
\frac{\rho(n,m)}{n^{-\bb_1}m^{1-\bb_2\al}} \rightarrow \sum_{i=0}^{\infty}w_{(i,\infty)}^{\al-1}(1+i)^{-\bb_1(\al-1)}\int\limits_{0}^{\infty} t^{-\bb_2(\al-1)}(1+t)^{-\bb_2}\dt,\ n,m\rightarrow\infty.
\end{equation*}

Fix $\gamma\in(0,1)$. Let us split the main sum \eqref{eq:mainsum} into two terms
\begin{equation}\label{4rhoSkaid}
\rho({n,m})=
\sum_{j\leq m^\g}\sum_{i=0}^{\infty}\w\left( c_{i,j}, c_{i+n,j+m} \right)+
\sum_{j> m^\g}\sum_{i=0}^{\infty}\w\left( c_{i,j}, c_{i+n,j+m} \right).
\end{equation}
The first term on the right side of equality can be estimated as follows
\begin{equation*}\begin{gathered}
\left|\sum_{j\leq m^\g}\sum_{i=0}^{\infty}\w\left( c_{i,j}, c_{i+n,j+m} \right)\right| \leq
\sum_{j\leq m^\g}\sum_{i=0}^{\infty}|\w\left( c_{i,j}, c_{i+n,j+m} \right)|\\
\leq
\sum_{j\leq m^\g}\sum_{i=0}^{\infty} |c_{i,j}|^{\al-1}|c_{i+n,j+m}|\\
\leq
e^{\al} \sum_{j\leq m^\g}\sum_{i=0}^{\infty} (1+i)^{-\bb_1(\al-1)}(1+j)^{-\bb_2(\al-1)}(1+i+n)^{-\bb_1}(1+j+m)^{-\bb_2}\\
= e^{\al}
\sum_{i=0}^{\infty}(1+i)^{-\bb_1(\al-1)}(1+i+n)^{-\bb_1}
\sum_{j\leq m^\g} (1+j)^{-\bb_2(\al-1)}(1+j+m)^{-\bb_2}.
\end{gathered} \end{equation*}

Notice that
\begin{equation}\label{ineq:4pirmasdaug}
\sum_{i=0}^{\infty}(1+i)^{-\bb_1(\al-1)}(1+i+n)^{-\bb_1}\leq
n^{-\bb_1} \sum_{i=0}^{\infty}(1+i)^{-\bb_1(\al-1)},
\end{equation}
and the series $\sum_{i=0}^{\infty}(1+i)^{-\bb_1(\al-1)}$ converge, since $\bb_1>\frac{1}{\al-1}$.
Simple inequalities give us
\begin{equation*}\begin{gathered}
\sum_{j\leq m^\g} (1+j)^{-\bb_2(\al-1)}(1+j+m)^{-\bb_2}=
\int\limits_{0}^{\lfloor m^\g \rfloor +1} (1+\lfloor s \rfloor)^{-\bb_2(\al-1)}(1+\lfloor s \rfloor+m)^{-\bb_2} \dd s\\
=
\int\limits_{0}^{\frac{\lfloor m^\g \rfloor +1}{m}} (1+\lfloor ms \rfloor)^{-\bb_2(\al-1)}(1+\lfloor ms \rfloor+m)^{-\bb_2} \dd ms \\
\leq
\int\limits_{0}^{\frac{\lfloor m^\g \rfloor +1}{m}}
(ms)^{-\bb_2(\al-1)}(ms+m)^{-\bb_2} \dd ms=
m^{1-\bb_2\al}
\int\limits_{0}^{\frac{\lfloor m^\g \rfloor +1}{m}}
s^{-\bb_2(\al-1)}(s+1)^{-\bb_2} \dd s\\
\leq
m^{1-\bb_2\al}
\int\limits_{0}^{\frac{\lfloor m^\g \rfloor +1}{m}}
s^{-\bb_2(\al-1)} \dd s
=
m^{1-\bb_2\al}\frac{\left( \frac{\lfloor m^\g \rfloor +1}{m} \right)^{1-\bb_2(\al-1)}}{1-\bb_2(\al-1)}.
\end{gathered}\end{equation*}
This bound together with \eqref{ineq:4pirmasdaug} implies that
\begin{equation*}
\frac{1}{n^{-\bb_1}m^{1-\bb_2\al}}\sum_{j\leq m^\g}\sum_{i=0}^{\infty}\w\left( c_{i,j}, c_{i+n,j+m} \right)\rightarrow 0,\ \text{as } n,m\rightarrow \infty,
\end{equation*}
since $1<\al< 2$ and $\frac{1}{\al}<\bb_2<\frac{1}{\al-1}$.

Let us investigate the second term in \eqref{4rhoSkaid}.
\begin{equation*}\begin{gathered}
\sum_{j> m^\g}\sum_{i=0}^{\infty}\w\left( c_{i,j}, c_{i+n,j+m} \right)\\
=
\sum_{i=0}^{\infty} \sum_{j> m^\g}
\w\left( w_{(i,j)}(1+i)^{-\bb_1}(1+j)^{-\bb_2}, w_{(i+n,j+m)}(1+i+n)^{-\bb_1}(1+j+m)^{-\bb_2} \right)\\
=
\sum_{i=0}^{\infty} \int\limits_{\lfloor m^\g \rfloor+1}^{\infty}
\w\left( w_{(i,\lfloor s \rfloor)}(1+i)^{-\bb_1}(1+\lfloor s \rfloor)^{-\bb_2}, w_{(i+n,\lfloor s \rfloor+m)}(1+i+n)^{-\bb_1}(1+\lfloor s \rfloor+m)^{-\bb_2} \right) \dd s
\end{gathered}\end{equation*}
\begin{equation*}   \begin{gathered}
=
\sum_{i=0}^{\infty} \int\limits_{\frac{\lfloor m^\g \rfloor+1}{m}}^{\infty}
\w\left( w_{(i,\lfloor ms \rfloor)}(1+i)^{-\bb_1}(1+\lfloor ms \rfloor)^{-\bb_2}, w_{(i+n,\lfloor ms \rfloor+m)}(1+i+n)^{-\bb_1}(1+\lfloor ms \rfloor+m)^{-\bb_2} \right) \dd ms
\end{gathered}\end{equation*}
\begin{equation*} \footnotesize  \begin{gathered}
=
m^{1-\bb_2\al} \sum_{i=0}^{\infty} \int\limits_{\frac{\lfloor m^\g \rfloor+1}{m}}^{\infty}
\w\left(
w_{(i,\lfloor ms \rfloor)}
(1+i)^{-\bb_1}
\left(\frac{1+\lfloor ms\rfloor}{m} \right)^{-\bb_2},
w_{(i+n,\lfloor ms \rfloor+m)}
(1+i+n)^{-\bb_1}
\left(\frac{1+\lfloor ms \rfloor}{m}+1\right)^{-\bb_2} \right) \dd s.
\end{gathered}\end{equation*}

Denoting the  integrand as $\y_{i,s,n,m}$, we have
\begin{equation*}\begin{gathered}
\frac{\y_{i,s,n,m}}{n^{-\bb_1}}\\
=\frac{n^{\bb_1}
	w_{(i,\lfloor ms \rfloor)}
	(1+i)^{-\bb_1}
	\left(\frac{1+\lfloor ms\rfloor}{m} \right)^{-\bb_2}
	w_{(i+n,\lfloor ms \rfloor+m)}
	(1+i+n)^{-\bb_1}
	\left(\frac{1+\lfloor ms \rfloor}{m}+1\right)^{-\bb_2}
}{
\left(  w_{(i,\lfloor ms \rfloor)}^2
(1+i)^{-2\bb_1}
\left(\frac{1+\lfloor ms\rfloor}{m} \right)^{-2\bb_2}+
w_{(i+n,\lfloor ms \rfloor+m)}^2
(1+i+n)^{-2\bb_1}
\left(\frac{1+\lfloor ms \rfloor}{m}+1\right)^{-2\bb_2} \right)^{\frac{2-\al}{2}} }\\
=\frac{
	w_{(i,\lfloor ms \rfloor)}
	(1+i)^{-\bb_1}
	\left(\frac{1+\lfloor ms\rfloor}{m} \right)^{-\bb_2}
	w_{(i+n,\lfloor ms \rfloor+m)}
	\left(\frac{1+i}{n}+1\right)^{-\bb_1}
	\left(\frac{1+\lfloor ms \rfloor}{m}+1\right)^{-\bb_2}
}{
\left(  w_{(i,\lfloor ms \rfloor)}^2
(1+i)^{-2\bb_1}
\left(\frac{1+\lfloor ms\rfloor}{m} \right)^{-2\bb_2}+
w_{(i+n,\lfloor ms \rfloor+m)}^2
(1+i+n)^{-2\bb_1}
\left(\frac{1+\lfloor ms \rfloor}{m}+1\right)^{-2\bb_2} \right)^{\frac{2-\al}{2}} }
\end{gathered}\end{equation*}
It is evident from the obtained expression that, for fixed $i,s$,
\begin{equation*}
\frac{\y_{i,s,n,m}}{n^{-\bb_1}}\rightarrow
\frac{
	w_{(i,\infty)}
	(1+i)^{-\bb_1}
	s^{-\bb_2}
	\left(s+1\right)^{-\bb_2}
}{
\left(  w_{(i,\infty)}^2
(1+i)^{-2\bb_1}
s^{-2\bb_2} \right)^{\frac{2-\al}{2}} }=
w_{(i,\infty)}^{\al-1} (1+i)^{-\bb_1(\al-1)} s^{-\bb_2(\al-1)} \left(s+1\right)^{-\bb_2}.
\end{equation*}

Properties of function $\w$ imply that
\begin{equation*}
\frac{|\y_{i,s,n,m}|}{n^{-\bb_1}}= n^{\bb_1}|\y_{i,s,n,m}|
\end{equation*}
\begin{equation*}\small
\leq
|w_{(i,\lfloor ms \rfloor)}|^{\al-1}
(1+i)^{-\bb_1(\al-1)}
\left(\frac{1+\lfloor ms\rfloor}{m} \right)^{-\bb_2(\al-1)}
|w_{(i+n,\lfloor ms \rfloor+m)}|
\left(\frac{1+i}{n}+1\right)^{-\bb_1}
\left(\frac{1+\lfloor ms \rfloor}{m}+1\right)^{-\bb_2}
\end{equation*}
\begin{equation*}
\leq
e^{\al}
(1+i)^{-\bb_1(\al-1)}
\left(\frac{1+\lfloor ms\rfloor}{m} \right)^{-\bb_2(\al-1)}
\left(\frac{1+\lfloor ms \rfloor}{m}+1\right)^{-\bb_2}
\end{equation*}
\begin{equation*}
\leq
e^{\al}
(1+i)^{-\bb_1(\al-1)}
s^{-\bb_2(\al-1)}
\left(s+1\right)^{-\bb_2}.
\end{equation*}

Notice that
\begin{equation*}
\sum_{i=0}^{\infty}\int\limits_{0}^{\infty}
(1+i)^{-\bb_1(\al-1)}
s^{-\bb_2(\al-1)}
\left(s+1\right)^{-\bb_2}\ds=
\sum_{i=0}^{\infty}(1+i)^{-\bb_1(\al-1)}
\int\limits_{0}^{\infty}            	    			
s^{-\bb_2(\al-1)}
\left(s+1\right)^{-\bb_2}\ds
<\infty,
\end{equation*}
since $-\bb_1(\al-1)<-1,\ -\bb_2(\al-1)>-1$ and $-\bb_2\al<-1$.

Thus, $\frac{|\y_{i,s,n,m}|}{n^{-\bb_1}}$ converges pointwise and is dominated by an integrable function. The Dominated Convergence Theorem can be applied to obtain
\begin{equation*}
\frac{\sum_{j> m^\g}\sum_{i=0}^{\infty}\w\left( c_{i,j}, c_{i+n,j+m} \right)}{n^{-\bb_1}m^{1-\bb_2\al}}
\end{equation*}
\begin{equation*}
=
\sum_{i=0}^{\infty} \int\limits_{\frac{\lfloor m^\g \rfloor+1}{m}}^{\infty}\frac{|\y_{i,s,n,m}|} {n^{-\bb_1}}
 \dd s
\rightarrow \sum_{i=0}^{\infty}w_{(i,\infty)}^{\al-1}(1+i)^{-\bb_1(\al-1)}\int\limits_{0}^{\infty} s^{-\bb_2(\al-1)}(1+s)^{-\bb_2}\dd s,\ n,m\rightarrow\infty.
\end{equation*}

This implies
\begin{equation}\label{IVRezult}
\frac{\rho({n,m})}{n^{-\bb_1}m^{1-\bb_2\al}} \rightarrow \sum_{i=0}^{\infty}w_{(i,\infty)}^{\al-1}(1+i)^{-\bb_1(\al-1)}\int\limits_{0}^{\infty} s^{-\bb_2(\al-1)}(1+s)^{-\bb_2}\dd s,\ n,m\rightarrow\infty.
\end{equation}

\subsection{The case $1<\al\leq 2,\ \bb_1=\frac{1}{\al-1},\ \bb_2=\frac{1}{\al-1};$}

  We will show that
  \begin{equation*}
  \frac{\rho(n,m)}{n^{-\bb_1}m^{-\bb_2}\ln n\ln m}\rightarrow 1,\ n,m\rightarrow\infty.
  \end{equation*}

  Denote
  \begin{equation*}\begin{gathered}
  q_n^{0}=0,\ q_n^{1}=\lfloor \ln n \rfloor+1,\ \ q_n^{2}=n+1,\ \  q_n^{3}=\infty,\\
  \hat{S}_{n,m}^{k,l}=\sum_{q_n^{k-1}\leq i<q_n^{k}}\sum_{q_m^{l-1}\leq j<q_m^{l}}\w\left( c_{i,j}, c_{i+n,j+m} \right),\\
  \hat{Z}_{n}^{k}(\bb_1)=\sum_{q_n^{k-1}\leq i<q_n^{k}}(1+i)^{-\bb_1(\al-1)}(1+i+n)^{-\bb_1}.
  \end{gathered}\end{equation*}

  The introduced notation enables us to write
  \begin{equation}\label{eq:2.5rhoSplit}
  \rho(n,m)=\sum_{k=1}^{3}\sum_{l=1}^{3}\hat{S}_{n,m}^{k,l}.
  \end{equation}

  Let us now investigate the asymptotic behavior of $\hat{Z}_{n}^{k}(\bb_1)$, $k=1,2,3$.
  \begin{equation*}
  \hat{Z}_{n}^{1}(\bb_1)=\sum_{q_n^{0}\leq i<q_n^{1}}(1+i)^{-\bb_1(\al-1)}(1+i+n)^{-\bb_1}=
  \sum_{0\leq i<\sv{\ln(n)}+1}(1+i)^{-\bb_1(\al-1)}(1+i+n)^{-\bb_1}
  \end{equation*}
  \begin{equation*}
  \leq
  n^{-\bb_1}\sum_{0\leq i<\sv{\ln(n)}+1}(1+i)^{-\bb_1(\al-1)}=
  n^{-\bb_1}\left(1+\sum_{1\leq i<\sv{\ln(n)}+1}(1+i)^{-1}\right)
  \end{equation*}
  \begin{equation*}
  =n^{-\bb_1}\left(1+\int_{1}^{\sv{\ln(n)}+1}(1+\sv{s})^{-1}\ds\right)\leq
  n^{-\bb_1}\left(1+\int_{1}^{\sv{\ln(n)}+1}s^{-1}\ds\right)
  \end{equation*}
  \begin{equation*}
  =
  n^{-\bb_1}\left(1+\ln\left({\sv{\ln(n)}+1}\right)\right).
  \end{equation*}
  The obtained upper bound implies
  \begin{equation*}
  \frac{\hat{Z}_{n}^{1}(\bb_1)}{n^{-\bb_1}\ln (n)}\rightarrow 0,\ n\rightarrow\infty.
  \end{equation*}

  We continue with $\hat{Z}_{n}^{2}(\bb_1)$:
  \begin{equation*}
  \hat{Z}_{n}^{2}(\bb_1)=\sum_{q_n^{1}\leq i<q_n^{2}}(1+i)^{-\bb_1(\al-1)}(1+i+n)^{-\bb_1}=
  \sum_{\sv{\ln(n)}+1\leq i<n+1}(1+i)^{-\bb_1(\al-1)}(1+i+n)^{-\bb_1}
  \end{equation*}
  \begin{equation*}
  =
  \int_{\sv{\ln(n)}+1}^{n+1}(1+\sv{s})^{-\bb_1(\al-1)}(1+\sv{s}+n)^{-\bb_1}\dd s.
  \end{equation*}
  Since
  \begin{equation*}
   \hat{Z}_{n}^{2}(\bb_1)=\int_{\sv{\ln(n)}+1}^{n+1}(1+\sv{s})^{-\bb_1(\al-1)}(1+\sv{s}+n)^{-\bb_1}\dd s\leq
  n^{-\bb_1}\int_{1}^{n+1}s^{-\bb_1(\al-1)}\dd s=n^{-\bb_1}\ln(n+1),
  \end{equation*}
  we have
  \begin{equation}\label{ineq:2.5Z2limsup}
  \limsup_{n\rightarrow\infty}\frac{\hat{Z}_{n}^{2}(\bb_1)}{n^{-\bb_1}\ln(n)}\leq 1.
  \end{equation}
  Next, choose arbitrary $\varepsilon\in(0,1)$. There exists $N\in\NN$ such that $n\varepsilon>\sv{\ln(n)}+2$ for all $n\geq N$. Then, for $n\geq N$ we have
  \begin{equation*}
  \hat{Z}_{n}^{2}(\bb_1)=\int_{\sv{\ln(n)}+1}^{n+1}(1+\sv{s})^{-\bb_1(\al-1)}(1+\sv{s}+n)^{-\bb_1}\dd s
  \end{equation*}
  \begin{equation*}
  \geq
  \int_{\sv{\ln(n)}+1}^{n+1}(1+s)^{-\bb_1(\al-1)}(1+s+n)^{-\bb_1}\dd s
  \end{equation*}
  \begin{equation*}
  =
  \int_{\sv{\ln(n)}+2}^{n+2}s^{-\bb_1(\al-1)}(s+n)^{-\bb_1}\dd s\geq
  \int_{\sv{\ln(n)}+2}^{n\varepsilon}s^{-\bb_1(\al-1)}(s+n)^{-\bb_1}\dd s
  \end{equation*}
  \begin{equation*}
 \geq
  \int_{\sv{\ln(n)}+2}^{n\varepsilon}s^{-\bb_1(\al-1)}(n\varepsilon+n)^{-\bb_1}\dd s=
  n^{-\bb_1}(1+\varepsilon)^{-\bb_1}\int_{\sv{\ln(n)}+2}^{n\varepsilon}s^{-1}\dd s
  \end{equation*}
  \begin{equation*}
  =
  n^{-\bb_1}(1+\varepsilon)^{-\bb_1}
  \left( \ln(n\varepsilon)-\ln\left( \sv{\ln(n)}+2 \right) \right)
  \end{equation*}
  \begin{equation*}
  =
  n^{-\bb_1}(1+\varepsilon)^{-\bb_1}
  \left( \ln(n)+\ln(\varepsilon)-\ln\left( \sv{\ln(n)}+2 \right) \right).
  \end{equation*}
  Thus,
  \begin{equation*}
  \liminf_{n\rightarrow\infty}\frac{\hat{Z}_{n}^{2}(\bb_1)}{n^{-\bb_1}\ln(n)}\geq \left( 1+\varepsilon \right)^{-\bb_1}
  \end{equation*}
   for arbitrary $\varepsilon\in(0,1)$. Passing to the limit as $\varepsilon\rightarrow0$, we obtain
   \begin{equation*}
   \liminf_{n\rightarrow\infty}\frac{\hat{Z}_{n}^{2}(\bb_1)}{n^{-\bb_1}\ln(n)}\geq 1.
   \end{equation*}
   This, together with \eqref{ineq:2.5Z2limsup}, gives
   \begin{equation}\label{lim:2.5z2}
   \lim_{n\rightarrow\infty}\frac{\hat{Z}_{n}^{2}(\bb_1)}{n^{-\bb_1}\ln(n)}=1.
   \end{equation}

   Finally, we estimate $\hat{Z}_{n}^{3}(\bb_1)$:
   \begin{equation*}
   \hat{Z}_{n}^{3}(\bb_1)=\sum_{q_n^{2}\leq i<q_n^{3}}(1+i)^{-\bb_1(\al-1)}(1+i+n)^{-\bb_1}=
   \sum_{i=n+1}^{\infty}(1+i)^{-\bb_1(\al-1)}(1+i+n)^{-\bb_1}
   \end{equation*}
   \begin{equation*}
   =
   \int_{n+1}^{\infty}(1+\sv{s})^{-\bb_1(\al-1)}(1+\sv{s}+n)^{-\bb_1}\dd s\leq
   \int_{n}^{\infty}s^{-\bb_1(\al-1)}(s+n)^{-\bb_1}\dd s
   \end{equation*}
   \begin{equation*}
   =
   \int_{1}^{\infty}(ns)^{-\bb_1(\al-1)}(ns+n)^{-\bb_1}\dd ns=
   n^{1-\bb_1\al}\int_{1}^{\infty}s^{-\bb_1(\al-1)}(s+1)^{-\bb_1}\dd s
   \end{equation*}
   \begin{equation*}
   =
   n^{-\bb_1}\int_{1}^{\infty}s^{-\bb_1(\al-1)}(s+1)^{-\bb_1}\dd s.
   \end{equation*}
  The obtained upper bound implies
  \begin{equation*}
  \frac{\hat{Z}_{n}^{3}(\bb_1)}{n^{-\bb_1}\ln (n)}\rightarrow 0,\ n\rightarrow\infty.
  \end{equation*}

  If $(k,l)\neq (2,2)$ we can use the inequality
  \begin{equation*}
  |\hat{S}_{n,m}^{k,l}|\leq e^\al \hat{Z}_{n}^{k}(\bb_1)\hat{Z}_{m}^{l}(\bb_2)
  \end{equation*}
  to obtain
  \begin{equation}\label{lim:2.5Snm}
  \frac{\hat{S}_{n,m}^{k,l}}{n^{-\bb_1}m^{-\bb_2}\ln (n)\ln (m)}\rightarrow 0,\ n,m\rightarrow\infty.
  \end{equation}

  It remains to investigate the asymptotic behavior of $\hat{S}_{n,m}^{2,2}$.
  Choose arbitrary $\tau\in(0,1)$. Since $w_{(i,j)}\rightarrow 1$ as $i,j\rightarrow\infty$, there exists $N\in\NN$ such that $1-\tau<w_{(i,j)}<1+\tau$ for $i,j\geq N$. Assume that $n$ and $m$ are large enough so that inequalities $\ln(n)>N$ and $\ln(m)>N$ hold.


  If $i\geq \ln(n)$ and $j\geq \ln(m)$ we have
  \begin{equation*}
  (1-\tau)(1+i)^{-\bb_1}(1+j)^{-\bb_2}\leq c_{i,j}\leq (1+\tau)(1+i)^{-\bb_1}(1+j)^{-\bb_2}.
  \end{equation*}
  This implies
  \begin{equation*}
  \w\left( c_{i,j}, c_{i+n,j+m} \right)\leq
  \frac{(1+\tau)^2}{(1-\tau)^{2-\al}}
  \frac{(1+i)^{-\bb_1}(1+j)^{-\bb_2}(1+i+n)^{-\bb_1}(1+j+m)^{-\bb_2}}{\left( (1+i)^{-2\bb_1}(1+j)^{-2\bb_2}+(1+i+n)^{-2\bb_1}(1+j+m)^{-2\bb_2} \right)^{\frac{2-\al}{2}}}
  \end{equation*}
  and
  \begin{equation*}
  \w\left( c_{i,j}, c_{i+n,j+m} \right)\geq
  \frac{(1-\tau)^2}{(1+\tau)^{2-\al}}
  \frac{(1+i)^{-\bb_1}(1+j)^{-\bb_2}(1+i+n)^{-\bb_1}(1+j+m)^{-\bb_2}}{\left( (1+i)^{-2\bb_1}(1+j)^{-2\bb_2}+(1+i+n)^{-2\bb_1}(1+j+m)^{-2\bb_2} \right)^{\frac{2-\al}{2}}}.
  \end{equation*}
  Thus,
  \begin{equation}\label{ineq:2.5AnmS22nm}
  \frac{(1-\tau)^2}{(1+\tau)^{2-\al}}A_{n,m}\leq \hat{S}_{n,m}^{2,2}\leq \frac{(1+\tau)^2}{(1-\tau)^{2-\al}}A_{n,m},
  \end{equation}
  where
  \begin{equation*}
  A_{n,m}=\sum_{i=\sv{\ln (n)}+1}^{n}\sum_{j=\sv{\ln (m)}+1}^{m}
  \frac{(1+i)^{-\bb_1}(1+j)^{-\bb_2}(1+i+n)^{-\bb_1}(1+j+m)^{-\bb_2}}{\left( (1+i)^{-2\bb_1}(1+j)^{-2\bb_2}+(1+i+n)^{-2\bb_1}(1+j+m)^{-2\bb_2} \right)^{\frac{2-\al}{2}}}.
  \end{equation*}

  Notice that
  \begin{equation*}
  A_{n,m}\leq \sum_{i=\sv{\ln (n)}+1}^{n}\sum_{j=\sv{\ln (m)}+1}^{m}
  (1+i)^{-\bb_1(\al-1)}(1+j)^{-\bb_2(\al-1)}(1+i+n)^{-\bb_1}(1+j+m)^{-\bb_2}
  \end{equation*}
  \begin{equation*}
  =\hat{Z}_{n}^{2}(\bb_1)\hat{Z}_{m}^{2}(\bb_2).
  \end{equation*}

  Previously obtained result \eqref{lim:2.5z2} implies
  \begin{equation*}
  \limsup_{n,m\rightarrow\infty}\frac{A_{n,m}}{n^{-\bb_1}m^{-\bb_2}\ln (n)\ln (m)}\leq 1.
  \end{equation*}

  Choose arbitrary $\varepsilon\in(0,1)$ and assume that $n,m$ are large enough so that $n\varepsilon>\sv{\ln(n)}+1$ and $m\varepsilon>\sv{\ln(m)}+1$. If $i\leq n\varepsilon$ we have
  \begin{equation*}
  (1+i+n)^{-\bb_1}=\frac{(1+i)^{\bb_1}}{(1+i+n)^{\bb_1}}(1+i)^{-\bb_1}\leq
  \frac{(1+n\varepsilon)^{\bb_1}}{(1+n\varepsilon+n)^{\bb_1}}(1+i)^{-\bb_1}\leq
  {\left(\frac{1}{n}+\varepsilon\right)^{\bb_1}}(1+i)^{-\bb_1}.
  \end{equation*}
  We also have
  \begin{equation*}
  (1+j+m)^{-\bb_2}\leq
  {\left(\frac{1}{m}+\varepsilon\right)^{\bb_2}}(1+j)^{-\bb_2}
  \end{equation*}
  for $j\leq m\varepsilon$.
   Then
  \begin{equation*}
  A_{n,m}\geq \sum_{i=\sv{\ln (n)}+1}^{\sv{n\varepsilon}}\sum_{j=\sv{\ln (m)}+1}^{\sv{m\varepsilon}}
  \frac{(1+i)^{-\bb_1}(1+j)^{-\bb_2}(1+i+n)^{-\bb_1}(1+j+m)^{-\bb_2}}{\left( (1+i)^{-2\bb_1}(1+j)^{-2\bb_2}+(1+i+n)^{-2\bb_1}(1+j+m)^{-2\bb_2} \right)^{\frac{2-\al}{2}}}
  \end{equation*}
  \begin{equation*}
  \geq \sum_{i=\sv{\ln (n)}+1}^{\sv{n\varepsilon}}\sum_{j=\sv{\ln (m)}+1}^{\sv{m\varepsilon}}
  \frac{(1+i)^{-\bb_1}(1+j)^{-\bb_2}(1+i+n)^{-\bb_1}(1+j+m)^{-\bb_2}}
  {\left( (1+i)^{-2\bb_1}(1+j)^{-2\bb_2}+{\left(\frac{1}{n}+\varepsilon\right)^{2\bb_1}}(1+i)^{-2\bb_1}{\left(\frac{1}{m}+\varepsilon\right)^{2\bb_2}}(1+j)^{-2\bb_2} \right)^{\frac{2-\al}{2}}}
  \end{equation*}
  \begin{equation*}
  = \sum_{i=\sv{\ln (n)}+1}^{\sv{n\varepsilon}}\sum_{j=\sv{\ln (m)}+1}^{\sv{m\varepsilon}}
  \frac{(1+i)^{-\bb_1(\al-1)}(1+j)^{-\bb_2(\al-1)}(1+i+n)^{-\bb_1}(1+j+m)^{-\bb_2}}
  {\left( 1+{\left(\frac{1}{n}+\varepsilon\right)^{2\bb_1}}{\left(\frac{1}{m}+\varepsilon\right)^{2\bb_2}}\right)^{\frac{2-\al}{2}}}
  \end{equation*}
  \begin{equation*}
  \geq
  \frac{(1+n\varepsilon+n)^{-\bb_1}(1+m\varepsilon+m)^{-\bb_2}}
  {\left( 1+{\left(\frac{1}{n}+\varepsilon\right)^{2\bb_1}}{\left(\frac{1}{m}+\varepsilon\right)^{2\bb_2}}\right)^{\frac{2-\al}{2}}}
  \sum_{i=\sv{\ln (n)}+1}^{\sv{n\varepsilon}}(1+i)^{-\bb_1(\al-1)}
  \sum_{j=\sv{\ln (m)}+1}^{\sv{m\varepsilon}}(1+j)^{-\bb_2(\al-1)}.
  \end{equation*}

  Let us deal with the sums on the right side.
  \begin{equation*}
  \sum_{i=\sv{\ln (n)}+1}^{\sv{n\varepsilon}}(1+i)^{-\bb_1(\al-1)}=
  \int_{\sv{\ln (n)}+1}^{\sv{n\varepsilon}+1}(1+\sv{s})^{-1}\dd s\geq
  \int_{\sv{\ln (n)}+1}^{\sv{n\varepsilon}+1}(1+s)^{-1}\dd s=
  \int_{\sv{\ln (n)}+2}^{\sv{n\varepsilon}+2}s^{-1}\dd s
  \end{equation*}
  \begin{equation*}
  \geq
  \int_{\sv{\ln (n)}+2}^{n\varepsilon}s^{-1}\dd s=\ln(n\varepsilon)-\ln\left( \sv{\ln (n)}+2 \right)=\ln(n)+\ln(\varepsilon)-\ln\left( \sv{\ln (n)}+2 \right).
  \end{equation*}
  Hence
  \begin{equation*}
  \liminf_{n\rightarrow\infty}\frac{\sum_{i=\sv{\ln (n)}+1}^{\sv{n\varepsilon}}(1+i)^{-\bb_1(\al-1)}}{\ln(n)}\geq 1.
  \end{equation*}
  Similarly
  \begin{equation*}
  \liminf_{m\rightarrow\infty}\frac{\sum_{j=\sv{\ln (m)}+1}^{\sv{m\varepsilon}}(1+j)^{-\bb_2(\al-1)}}{\ln(m)}\geq 1.
  \end{equation*}
  These results imply
  \begin{equation*}
  \liminf_{n,m\rightarrow\infty}\frac{A_{n,m}}{n^{-\bb_1}m^{-\bb_2}\ln (n)\ln (m)}\geq
   \frac{(\varepsilon+1)^{-\bb_1}(\varepsilon+1)^{-\bb_2}}
   {\left( 1+{\left(\varepsilon\right)^{2\bb_1}}{\left(\varepsilon\right)^{2\bb_2}}\right)^{\frac{2-\al}{2}}}=
   \frac{(\varepsilon+1)^{-\bb_1-\bb_2}}
   {\left( 1+\varepsilon^{2\bb_1+2\bb_2}\right)^{\frac{2-\al}{2}}}.
  \end{equation*}
  Since $\varepsilon$ was chosen arbitrarily, we have
  \begin{equation*}
  \liminf_{n,m\rightarrow\infty}\frac{A_{n,m}}{n^{-\bb_1}m^{-\bb_2}\ln (n)\ln (m)}\geq
 1.
  \end{equation*}

  Returning to \eqref{ineq:2.5AnmS22nm} we obtain
  \begin{equation*}
  \frac{(1-\tau)^2}{(1+\tau)^{2-\al}}
  \leq \liminf_{n,m\rightarrow\infty}\frac{\hat{S}_{n,m}^{2,2}}{n^{-\bb_1}m^{-\bb_2}\ln (n)\ln (m)}
  \leq \limsup_{n,m\rightarrow\infty}\frac{\hat{S}_{n,m}^{2,2}}{n^{-\bb_1}m^{-\bb_2}\ln (n)\ln (m)}\leq
  \frac{(1+\tau)^2}{(1-\tau)^{2-\al}}
  \end{equation*}
  for arbitrary $\tau$. Passing to the limit as $\tau\rightarrow0$ gives us
  \begin{equation*}
  \lim_{n,m\rightarrow\infty}\frac{\hat{S}_{n,m}^{2,2}}{n^{-\bb_1}m^{-\bb_2}\ln (n)\ln (m)}=1.
  \end{equation*}
 Together with \eqref{eq:2.5rhoSplit} and \eqref{lim:2.5Snm} this implies
 \begin{equation*}
 \lim_{n,m\rightarrow\infty}\frac{\rho(n,m)}{n^{-\bb_1}m^{-\bb_2}\ln (n)\ln (m)}=1.
 \end{equation*}

\section{Proof of Theorem \ref{thm2}}

Similarly to the proof of Theorem \ref{thm1}, we  consider only part of the sets of parameters, since  the other cases are investigated in a similar way.
Here we will investigate  the following sets of parameters:
\begin{enumerate}
	\item $ 1<\al\le 2,\ \frac{1}{\al}<\bb_i<\frac{1}{\al-1 },\ i=1,2 $ and $\frac{1}{\bb_1}+\frac{1}{\bb_2}>\al$;
	\item $ 1<\al\le 2,\ \frac{1}{\al-1}<\bb_1$ and $\frac{1}{\al}<\bb_2<1$;
	\item $ 1<\al< 2,\ \bb_1=\frac{1}{\al-1}$ and $1<\bb_2<\frac{1}{\al-1}$;
\end{enumerate}

\subsection{The case $ 1<\al\le 2,\ \frac{1}{\al}<\bb_i<\frac{1}{\al-1 },\ i=1,2,\ \frac{1}{\bb_1}+\frac{1}{\bb_2}>\al$;}

The spectral covariance equals
 \begin{equation*}\footnotesize\begin{gathered}
 		\rho(n,-m)=
 		\sum_{i=0}^{\infty}\sum_{j=0}^{\infty}
 		\frac
 		{w_{(i+n,j)}(1+i+n)^{-\bb_1}(1+j)^{-\bb_2}w_{(i,j+m)}(1+i)^{-\bb_1}(1+j+m)^{-\bb_2}}
 		{\left( w^2_{(i+n,j)}(1+i+n)^{-2\bb_1}(1+j)^{-2\bb_2}+w^2_{(i,j+m)}(1+i)^{-2\bb_1}(1+j+m)^{-2\bb_2} \right)^{\frac{2-\al}{2}}}
 		\\=
 		\int\limits_{0}^{\infty}
 		\int\limits_{0}^{\infty}
 		\frac
 		{w_{(\sv{t}+n,\sv{s})}w_{(\sv{t},\sv{s}+m)}(1+\sv{t}+n)^{-\bb_1}(1+\sv{s})^{-\bb_2}(1+\sv{t})^{-\bb_1}(1+\sv{s}+m)^{-\bb_2}\dd t \dd s}
 		{\left( w^2_{(\sv{t}+n,\sv{s})}(1+\sv{t}+n)^{-2\bb_1}(1+\sv{s})^{-2\bb_2}+w^2_{(\sv{t},\sv{s}+m)}(1+\sv{t})^{-2\bb_1}(1+\sv{s}+m)^{-2\bb_2} \right)^{\frac{2-\al}{2}}}
 		\\=
 		\int\limits_{0}^{\infty}
 		\int\limits_{0}^{\infty}
 		\frac
 		{w_{(\sv{nt}+n,\sv{ms})}w_{(\sv{nt},\sv{ms}+m)}(1+\sv{nt}+n)^{-\bb_1}(1+\sv{ms})^{-\bb_2}(1+\sv{nt})^{-\bb_1}(1+\sv{ms}+m)^{-\bb_2}\dd nt \dd ms}
 		{\left( w^2_{(\sv{nt}+n,\sv{ms})}(1+\sv{nt}+n)^{-2\bb_1}(1+\sv{ms})^{-2\bb_2}+w^2_{(\sv{nt},\sv{ms}+m)}(1+\sv{nt})^{-2\bb_1}(1+\sv{ms}+m)^{-2\bb_2} \right)^{\frac{2-\al}{2}}}
 		\\=
 		\int\limits_{0}^{\infty}
 		\int\limits_{0}^{\infty}
 		\frac
 		{n^{1-\bb_1\al}m^{1-\bb_2\al}w_{(\sv{nt}+n,\sv{ms})}w_{(\sv{nt},\sv{ms}+m)}
 			(1+\frac{\sv{nt}+1}{n})^{-\bb_1}(\frac{\sv{ms}+1}{m})^{-\bb_2}(\frac{\sv{nt}+1}{n})^{-\bb_1}(1+\frac{\sv{ms}+1}{m})^{-\bb_2}\dd t \dd s}
 		{\left( w^2_{(\sv{nt}+n,\sv{ms})}(1+\frac{\sv{nt}+1}{n})^{-2\bb_1}(\frac{\sv{ms}+1}{m})^{-2\bb_2}+w^2_{(\sv{nt},\sv{ms}+m)}(\frac{\sv{nt}+1}{n})^{-2\bb_1}(1+\frac{\sv{ms}+1}{m})^{-2\bb_2} \right)^{\frac{2-\al}{2}}}
 	\end{gathered}\end{equation*}
 	
 	Thus,
 	\begin{equation}\label{IRhoIsr}\footnotesize\begin{gathered}
 			\frac{\rho(n,-m)}{n^{1-\bb_1\al}m^{1-\bb_2\al}}=
 			\int\limits_{0}^{\infty}
 			\int\limits_{0}^{\infty}
 			y_{t,s,n,m}
 			\dd t \dd s,
 		\end{gathered}\end{equation}
 		where
 		\begin{equation*}
 		y_{t,s,n,m}=\frac
 		{w_{(\sv{nt}+n,\sv{ms})}w_{(\sv{nt},\sv{ms}+m)}
 			(1+\frac{\sv{nt}+1}{n})^{-\bb_1}(\frac{\sv{ms}+1}{m})^{-\bb_2}(\frac{\sv{nt}+1}{n})^{-\bb_1}(1+\frac{\sv{ms}+1}{m})^{-\bb_2}}
 		{\left( w^2_{(\sv{nt}+n,\sv{ms})}(1+\frac{\sv{nt}+1}{n})^{-2\bb_1}(\frac{\sv{ms}+1}{m})^{-2\bb_2}+w^2_{(\sv{nt},\sv{ms}+m)}(\frac{\sv{nt}+1}{n})^{-2\bb_1}(1+\frac{\sv{ms}+1}{m})^{-2\bb_2} \right)^{\frac{2-\al}{2}}}.
 		\end{equation*}
 		
	 	For fixed values of $t,s$ we have
 		\begin{equation}\label{IRibine}
 			y_{t,s,n,m}\rightarrow
 			\frac
 			{      (1+t)^{-\bb_1}s^{-\bb_2}t^{-\bb_1}(1+s)^{-\bb_2}}
 			{\left( (1+t)^{-2\bb_1}s^{-2\bb_2}+
 				t^{-2\bb_1}(1+s)^{-2\bb_2} \right)^{\frac{2-\al}{2}}},\ n,m\rightarrow\infty.
 		\end{equation}
 		
 		Next we show that $y_{t,s,n,m}$ is dominated by an integrable function.
 		If $s\geq1$ we have
 		\begin{equation*}\begin{gathered}
 		\abs{y_{t,s,n,m}} \leq
 		e^\al
 		\left(1+\frac{\sv{nt}+1}{n}\right)^{-\bb_1}
 		\left(\frac{\sv{ms}+1}{m}\right)^{-\bb_2}
 		\left(\frac{\sv{nt}+1}{n}\right)^{-\bb_1(\al-1)}
 		\left(1+\frac{\sv{ms}+1}{m}\right)^{-\bb_2(\al-1)}\\ \leq
 		e^\al
 		\left(1+t\right)^{-\bb_1}
 		s^{-\bb_2}
 		t^{-\bb_1(\al-1)}
 		\left(1+s\right)^{-\bb_2(\al-1)}    ,
 		\end{gathered}\end{equation*}
 		and the function on the right hand side of inequality is integrable in the set $\{(t,s):t>0,\ s\geq1\}$.
 		
 		If $0<s<1,t\geq 1$ we have
 		\begin{equation*}\begin{gathered}
 		\abs{y_{t,s,n,m}} \leq
 		e^\al
 		\left(1+\frac{\sv{nt}+1}{n}\right)^{-\bb_1(\al-1)}
 		\left(\frac{\sv{ms}+1}{m}\right)^{-\bb_2(\al-1)}
 		\left(\frac{\sv{nt}+1}{n}\right)^{-\bb_1}
 		\left(1+\frac{\sv{ms}+1}{m}\right)^{-\bb_2}\\ \leq
 		e^\al
 		\left(1+t\right)^{-\bb_1(\al-1)}
 		s^{-\bb_2(\al-1)}
 		t^{-\bb_1}
 		\left(1+s\right)^{-\bb_2}     ,
 		\end{gathered}\end{equation*}
 		and the function on the right hand side of inequality is integrable in the set $\lbrace (s,t):0<s<1,t\geq1\rbrace$.
 		
 		It remains to investigate the set $\lbrace (s,t):0<s<1,0<t<1\rbrace$. If $0<s<1,0<t<1$ we have
 		\begin{equation*}
 		\abs{y_{t,s,n,m}}\leq
 		\frac{e^2}{d^{2-\al}}\frac
 		{
 			(1+\frac{\sv{nt}+1}{n})^{-\bb_1}(\frac{\sv{ms}+1}{m})^{-\bb_2}(\frac{\sv{nt}+1}{n})^{-\bb_1}(1+\frac{\sv{ms}+1}{m})^{-\bb_2}}
 		{\left( (1+\frac{\sv{nt}+1}{n})^{-2\bb_1}(\frac{\sv{ms}+1}{m})^{-2\bb_2}+
 			(\frac{\sv{nt}+1}{n})^{-2\bb_1}(1+\frac{\sv{ms}+1}{m})^{-2\bb_2} \right)^{\frac{2-\al}{2}}}
 		\end{equation*}
 		\begin{equation*}
 		\leq
 		\frac{e^2}{d^{2-\al}}\frac
 		{
 			(\frac{\sv{ms}+1}{m})^{-\bb_2}(\frac{\sv{nt}+1}{n})^{-\bb_1}}
 		{\left( 2^{-2\bb_1}(\frac{\sv{ms}+1}{m})^{-2\bb_2}+
 			(\frac{\sv{nt}+1}{n})^{-2\bb_1}2^{-2\bb_2} \right)^{\frac{2-\al}{2}}}
 		\end{equation*}
 		\begin{equation*}
 		\leq
 		\frac{e^2}{d^{2-\al}}2^{(\bb_1+\bb_2)(2-\al)}\frac
 		{
 			(\frac{\sv{ms}+1}{m})^{-\bb_2}(\frac{\sv{nt}+1}{n})^{-\bb_1}}
 		{\left( (\frac{\sv{ms}+1}{m})^{-2\bb_2}+
 			(\frac{\sv{nt}+1}{n})^{-2\bb_1}
 			 \right)^{\frac{2-\al}{2}}}
 		\end{equation*}
 		\begin{equation*}
 		\leq
 		\frac{e^2}{d^{2-\al}}2^{(\bb_1+\bb_2)(2-\al)}\frac
 		{
 			(\frac{\sv{ms}+1}{m})^{-\bb_2}(\frac{\sv{nt}+1}{n})^{-\bb_1}}
 		{\left(
 			\max\left((\frac{\sv{ms}+1}{m})^{-2\bb_2},
 			(\frac{\sv{nt}+1}{n})^{-2\bb_1}\right)
 			\right)^{\frac{2-\al}{2}}}
 		\end{equation*}
 		\begin{equation*}
 		=
 		\frac{e^2}{d^{2-\al}}2^{(\bb_1+\bb_2)(2-\al)}\frac
 		{
 			(\frac{\sv{ms}+1}{m})^{-\bb_2}(\frac{\sv{nt}+1}{n})^{-\bb_1}}{
 			\max
 			\left(
 			(\frac{\sv{ms}+1}{m})^{-\bb_2(2-\al)},
 			(\frac{\sv{nt}+1}{n})^{-\bb_1(2-\al)}
 			\right)}
 		\end{equation*}
 		\begin{equation*}\small
 		=
 		\frac{e^2}{d^{2-\al}}2^{(\bb_1+\bb_2)(2-\al)}
 		\left( \frac{\sv{ms}+1}{m} \right)^{-\bb_2}
 		\left( \frac{\sv{nt}+1}{n} \right)^{-\bb_1}
 			\min
 			\left(
 			\left( \frac{\sv{ms}+1}{m} \right)^{\bb_2(2-\al)},
 			\left( \frac{\sv{nt}+1}{n} \right)^{\bb_1(2-\al)}
 			\right)
 		\end{equation*}
 		\begin{equation*}\small
 		=
 		\frac{e^2}{d^{2-\al}}2^{(\bb_1+\bb_2)(2-\al)}
 		\min
 		\left(
 		\left( \frac{\sv{ms}+1}{m} \right)^{-\bb_2(\al-1)}
 		\left( \frac{\sv{nt}+1}{n} \right)^{-\bb_1},
 		\left( \frac{\sv{ms}+1}{m} \right)^{-\bb_2}
 		\left( \frac{\sv{nt}+1}{n} \right)^{-\bb_1(\al-1)}
 		\right)
 		\end{equation*}
 		\begin{equation*}
 		\leq
 		\frac{e^2}{d^{2-\al}}2^{(\bb_1+\bb_2)(2-\al)}
 		\min
 		\left(
 		s^{-\bb_2(\al-1)}
 		t^{-\bb_1},
 		s^{-\bb_2}
 		t^{-\bb_1(\al-1)}
 		\right)
 		\end{equation*}
 		\begin{equation*}
 		=
 		\frac{e^2}{d^{2-\al}}2^{(\bb_1+\bb_2)(2-\al)}
 		s^{-\bb_2}
 		t^{-\bb_1}
 		\min
 		\left(
 		s^{\bb_2(2-\al)}
 		,
 		t^{\bb_1(2-\al)}
 		\right)
 		\end{equation*}
 		\begin{equation*}
 		=
 		\frac{e^2}{d^{2-\al}}2^{(\bb_1+\bb_2)(2-\al)}
 		s^{-\bb_2}
 		t^{-\bb_1}
 		\left( \min
 		\left(
 		s^{2\bb_2}
 		,
 		t^{2\bb_1}
 		\right) \right)^{\frac{2-\al}{2}}
 		\end{equation*}
 		\begin{equation*}
 		=
 		\frac{e^2}{d^{2-\al}}2^{(\bb_1+\bb_2)(2-\al)}
 		\frac{s^{-\bb_2}
 			t^{-\bb_1}}{
 		\left( \max
 		\left(
 		s^{-2\bb_2}
 		,
 		t^{-2\bb_1}
 		\right) \right)^{\frac{2-\al}{2}}}
 		\end{equation*}
 		\begin{equation*}
 		\leq
 		\frac{e^2}{d^{2-\al}}2^{(\bb_1+\bb_2)(2-\al)}
 		2^{\frac{2-\al}{2}}
 		\frac{s^{-\bb_2}
 			t^{-\bb_1}}{
 			\left(
 			s^{-2\bb_2}
 			+
 			t^{-2\bb_1}
 			 \right)^{\frac{2-\al}{2}}}.
 		\end{equation*}
 		
 		We continue by showing that
 		\begin{equation*}
 		\int_{0}^{1}\int_{0}^{1}
 		\frac{s^{-\bb_2}
 			t^{-\bb_1}}{
 			\left(
 			s^{-2\bb_2}
 			+
 			t^{-2\bb_1}
 			\right)^{\frac{2-\al}{2}}}
 		\dd t \dd s<\infty.
 		\end{equation*}
 		
 		Change of variables $t=x^{-\frac{1}{\bb_1}},\ s=y^{-\frac{1}{\bb_2}}$ gives us
 					\begin{equation*}\begin{gathered}
 							\int\limits_{0}^{1}
 							\int\limits_{0}^{1} \frac{t^{-\bb_1}s^{-\bb_2}}{\left( t^{-2\bb_1}+s^{-2\bb_2} \right)^{\frac{2-\al}{2}}}\dd t \dd s\\=
 							\frac{1}{\bb_1\bb_2}
 							\int\limits_{1}^{\infty}
 							\int\limits_{1}^{\infty}
 							\frac{xy}{\left( x^2+y^2 \right)^{\frac{2-\al}{2}}} x^{-\frac{1}{\bb_1}-1}y^{-\frac{1}{\bb_2}-1}
 							\dd x \dd y=
 							\frac{1}{\bb_1\bb_2}
 							\int\limits_{1}^{\infty}
 							\int\limits_{1}^{\infty}
 							\frac{x^{-\frac{1}{\bb_1}}y^{-\frac{1}{\bb_2}}}{\left( x^2+y^2 \right)^{\frac{2-\al}{2}}}
 							\dd x \dd y,
 						\end{gathered}\end{equation*}
	 	We continue by applying the change of variables $x=r\cos\phi,\ y=r\sin\phi$, which gives us
 						\begin{equation*}\begin{gathered}
 								\int\limits_{1}^{\infty}
 								\int\limits_{1}^{\infty}
 								\frac{x^{-\frac{1}{\bb_1}}y^{-\frac{1}{\bb_2}}}{\left( x^2+y^2 \right)^{\frac{2-\al}{2}}}
 								\dd x \dd y\\=
 								\int\limits_{1}^{\infty}
 								\int\limits_{\arcsin\frac{1}{r}}^{\frac{\pi}{2}-\arcsin\frac{1}{r}}
 								\frac{
 									r^{1-\frac{1}{\bb_2}-\frac{1}{\bb_1}}
 									(\cos\phi)^{-\frac{1}{\bb_1}}(\sin\phi)^{-\frac{1}{\bb_2}}
 								}{r^{2-\al}}
 								\dd \phi \dd r\\=
 								\int\limits_{1}^{\infty}
 								r^{\al-1-\frac{1}{\bb_2}-\frac{1}{\bb_1}}
 								\int\limits_{\arcsin\frac{1}{r}}^{\frac{\pi}{2}-\arcsin\frac{1}{r}}
 								(\cos\phi)^{-\frac{1}{\bb_1}}(\sin\phi)^{-\frac{1}{\bb_2}}
 								\dd \phi \dd r \\
 								=\int\limits_{1}^{\infty}
 								r^{\al-1-\frac{1}{\bb_2}-\frac{1}{\bb_1}}\left(
 								\int\limits_{\arcsin\frac{1}{r}}^{\frac{\pi}{4}}
 								(\cos\phi)^{-\frac{1}{\bb_1}}(\sin\phi)^{-\frac{1}{\bb_2}}
 								\dd \phi +
 								\int\limits_{\frac{\pi}{4}}^{\frac{\pi}{2}-\arcsin\frac{1}{r}}
 								(\cos\phi)^{-\frac{1}{\bb_1}}(\sin\phi)^{-\frac{1}{\bb_2}}
 								\dd \phi\right)\dd r
 							\end{gathered}\end{equation*}
 			The first integral in parentheses can be bounded above as
 			\begin{equation*}
 			\int\limits_{\arcsin\frac{1}{r}}^{\frac{\pi}{4}}
 			(\cos\phi)^{-\frac{1}{\bb_1}}(\sin\phi)^{-\frac{1}{\bb_2}}
 			\dd \phi\leq
 			2^{\frac{1}{2\bb_1}}2^{\frac{1}{\bb_2}}
 			\int\limits_{\arcsin\frac{1}{r}}^{\frac{\pi}{4}}
 			\phi^{-\frac{1}{\bb_2}}
 			\dd \phi,
 			\end{equation*}
 			since $\cos\phi\geq2^{-1/2}$ and $\sin\phi\geq \phi/2$ when $\phi\in[0,\frac{\pi}{4}]$. Let us denote
 			\begin{equation*}
 			I_1(r)=\int\limits_{\arcsin\frac{1}{r}}^{\frac{\pi}{4}}
 			\phi^{-\frac{1}{\bb_2}}
 			\dd \phi=
 			\begin{cases}
 			\frac{\left( \frac{\pi}{4} \right)^{1-\frac{1}{\bb_2}}-
 				\left( \arcsin\frac{1}{r} \right)^{1-\frac{1}{\bb_2}}}{1-\frac{1}{\bb_2}},\ \bb_2\neq1,\\
 			\ln\left( \frac{\pi}{4} \right)-\ln\left( \arcsin\frac{1}{r} \right),\ \bb_2=1.
 			\end{cases}
 			\end{equation*}
 			In a similar way we can bound the second integral:
 			\begin{equation*}
 			\int\limits_{\frac{\pi}{4}}^{\frac{\pi}{2}-\arcsin\frac{1}{r}}
 			(\cos\phi)^{-\frac{1}{\bb_1}}(\sin\phi)^{-\frac{1}{\bb_2}}
 			\dd \phi=
 			\int\limits_{\arcsin\frac{1}{r}}^{\frac{\pi}{4}}
 			\left( \cos\left( \frac{\pi}{2}-\phi \right) \right)^{-\frac{1}{\bb_1}}\left( \sin\left( \frac{\pi}{2}-\phi \right) \right)^{-\frac{1}{\bb_2}}
 			\dd \phi
 			\end{equation*}
 			\begin{equation*}
 			=
 			\int\limits_{\arcsin\frac{1}{r}}^{\frac{\pi}{4}}
 			\left( \sin\phi \right)^{-\frac{1}{\bb_1}}
 			\left( \cos\phi \right)^{-\frac{1}{\bb_2}}
 			\dd \phi\leq 2^{\frac{1}{2\bb_2}}2^{\frac{1}{\bb_1}}I_2(r),
 			\end{equation*}
 			where
 			\begin{equation*}
 			I_2(r)=\int\limits_{\arcsin\frac{1}{r}}^{\frac{\pi}{4}}
 			\phi^{-\frac{1}{\bb_1}}
 			\dd \phi=
 			\begin{cases}
 			\frac{\left( \frac{\pi}{4} \right)^{1-\frac{1}{\bb_1}}-
 				\left( \arcsin\frac{1}{r} \right)^{1-\frac{1}{\bb_1}}}{1-\frac{1}{\bb_1}},\ \bb_1\neq1,\\
 			\ln\left( \frac{\pi}{4} \right)-\ln\left( \arcsin\frac{1}{r} \right),\ \bb_1=1.
 			\end{cases}
 			\end{equation*}		
 	We have shown that for some constant $K_1$
 	\begin{equation*}
 		\int\limits_{0}^{1}
 		\int\limits_{0}^{1} \frac{t^{-\bb_1}s^{-\bb_2}}{\left( t^{-2\bb_1}+s^{-2\bb_2} \right)^{\frac{2-\al}{2}}}\dd t \dd s\leq
 		K_1
 		\int\limits_{1}^{\infty}
 		r^{\al-1-\frac{1}{\bb_2}-\frac{1}{\bb_1}}
 		\left(
 		I_1(r)+I_2(r)
 		\right)\dd r.
 	\end{equation*}

 	To show that the integral is finite, it is enough to investigate the asymptotic behavior of $r^{\al-1-\frac{1}{\bb_2}-\frac{1}{\bb_1}}I_i(r)$ as $r\rightarrow\infty$ for $i=1,2$.
 	Since $\arcsin x\sim x,\ x\rightarrow 0$, we have
 	\begin{equation}\label{II1}
 	r^{\al-1-\frac{1}{\bb_2}-\frac{1}{\bb_1}}I_1(r)\sim
 	\begin{cases}
 	r^{\al-1-\frac{1}{\bb_2}-\frac{1}{\bb_1}}\frac{\left( \frac{\pi}{4} \right)^{1-\frac{1}{\bb_2}}}{1-\frac{1}{\bb_2}},\ \bb_2>1,\\
 	r^{\al-1-\frac{1}{\bb_2}-\frac{1}{\bb_1}}\ln r,\ \bb_2=1,\\
 	\frac{1}{\frac{1}{\bb_2}-1}r^{\al-2-\frac{1}{\bb_1}},\ \bb_2<1,
 	\end{cases}
 	\end{equation}
 	as $r\rightarrow\infty$. In the case under consideration we have $ 1<\al\le 2,\ \frac{1}{\al}<\bb_i<\frac{1}{\al-1 },\ i=1,2,\ \frac{1}{\bb_1}+\frac{1}{\bb_2}>\al$, so
 	\begin{equation*}
 	\int_{1}^{\infty}r^{\al-1-\frac{1}{\bb_2}-\frac{1}{\bb_1}}I_1(r)\dd r<\infty.
 	\end{equation*}
 	In the same way we show that
 	\begin{equation*}
 	\int_{1}^{\infty}r^{\al-1-\frac{1}{\bb_2}-\frac{1}{\bb_1}}I_2(r)\dd r<\infty.
 	\end{equation*}
 	
 	In conclusion, we have shown that $y_{t,s,n,m}$ is bounded above by an integrable function. Application of The Dominated Convergence Theorem gives us
 		\begin{equation*}\begin{gathered}
 		\frac{\rho(n,-m)}{n^{1-\bb_1\al}m^{1-\bb_2\al}}=
 		\int\limits_{0}^{\infty}
 		\int\limits_{0}^{\infty}
 		y_{t,s,n,m}
 		\dd t \dd s\rightarrow
 		\int\limits_{0}^{\infty}
 		\int\limits_{0}^{\infty}
 		\frac
 		{      (1+t)^{-\bb_1}s^{-\bb_2}t^{-\bb_1}(1+s)^{-\bb_2}}
 		{\left( (1+t)^{-2\bb_1}s^{-2\bb_2}+
 			t^{-2\bb_1}(1+s)^{-2\bb_2} \right)^{\frac{2-\al}{2}}}
 		\dd t \dd s,\ n,m\rightarrow\infty.
 		\end{gathered}\end{equation*}


\subsection{The case $ 1<\al\le 2,\ \frac{1}{\al-1}<\bb_1$ and $\frac{1}{\al}<\bb_2<1$;}
In this case we have
\begin{equation*}\small\begin{gathered}
\rho(n,-m)\\=
\sum_{i=0}^{\infty}\sum_{j=0}^{\infty}
\w\left( w_{(i+n,j)}(1+i+n)^{-\bb_1}(1+j)^{-\bb_2}, w_{(i,j+m)}(1+i)^{-\bb_1}(1+j+m)^{-\bb_2} \right)
\\=
\sum_{i=0}^{\infty}
\int\limits_{0}^{\infty}
\w\left( w_{(i+n,\sv{s})}(1+i+n)^{-\bb_1}(1+\sv{s})^{-\bb_2}, w_{(i,\sv{s}+m)}(1+i)^{-\bb_1}(1+\sv{s}+m)^{-\bb_2} \right) \dd s
\\=
\sum_{i=0}^{\infty}
\int\limits_{0}^{\infty}
\w\left( w_{(i+n,\sv{ms})}(1+i+n)^{-\bb_1}(1+\sv{ms})^{-\bb_2}, w_{(i,\sv{ms}+m)}(1+i)^{-\bb_1}(1+\sv{ms}+m)^{-\bb_2} \right) \dd ms
\\=
\sum_{j=0}^{\infty}\int\limits_{0}^{\infty}
\frac
{w_{(i+n,\sv{ms})}w_{(i,\sv{ms}+m)}
	(1+\frac{i+1}{n})^{-\bb_1}(\frac{1+\sv{ms}}{m})^{-\bb_2} (1+i)^{-\bb_1}(1+\frac{\sv{ms}+1}{m})^{-\bb_2} n^{-\bb_1}m^{1-\bb_2\al}\dd s}
{\left(
	w_{(i+n,\sv{ms})}^2
	(1+i+n)^{-2\bb_1}(\frac{1+\sv{ms}}{m})^{-2\bb_2}+
	w_{(i,\sv{ms}+m)}^2(1+i)^{-2\bb_1}(1+\frac{\sv{ms}+1}{m})^{-2\bb_2} \right)^{\frac{2-\al}{2}}}.
\end{gathered}\end{equation*}

Let us denote
\begin{equation*}\small
y_{i,s,n,m}:=\frac
{w_{(i+n,\sv{ms})}w_{(i,\sv{ms}+m)}
	(1+\frac{i+1}{n})^{-\bb_1}(\frac{1+\sv{ms}}{m})^{-\bb_2} (1+i)^{-\bb_1}(1+\frac{\sv{ms}+1}{m})^{-\bb_2} n^{-\bb_1}m^{1-\bb_2\al}\dd s}
{\left(
	w_{(i+n,\sv{ms})}^2
	(1+i+n)^{-2\bb_1}(\frac{1+\sv{ms}}{m})^{-2\bb_2}+
	w_{(i,\sv{ms}+m)}^2(1+i)^{-2\bb_1}(1+\frac{\sv{ms}+1}{m})^{-2\bb_2} \right)^{\frac{2-\al}{2}}},
\end{equation*}
then
\begin{equation*}
\rho(n,-m)=
\sum_{j=0}^{\infty}\int\limits_{0}^{\infty}
y_{i,s,n,m}
n^{-\bb_1}m^{1-\bb_2\al}\dd s.
\end{equation*}

For fixed values of $i,s$ we have, as $n,m\rightarrow\infty$,
\begin{equation*}
y_{i,s,n,m}\rightarrow
\frac
{w_{(i,\infty)}
	 (1+i)^{-\bb_1}s^{-\bb_2}(1+s)^{-\bb_2}}
{\left(
	w_{(i,\infty)}^2(1+i)^{-2\bb_1}(1+s)^{-2\bb_2} \right)^{\frac{2-\al}{2}}}=
w_{(i,\infty)}^{\al-1}(1+i)^{-\bb_1(\al-1)}s^{-\bb_2}(1+s)^{-\bb_2(\al-1)}.
\end{equation*}

Also
\begin{equation*}
y_{i,s,n,m}\leq
\frac
{w_{(i+n,\sv{ms})}w_{(i,\sv{ms}+m)}
	(1+\frac{i+1}{n})^{-\bb_1}(\frac{1+\sv{ms}}{m})^{-\bb_2} (1+i)^{-\bb_1}(1+\frac{\sv{ms}+1}{m})^{-\bb_2} }
{\left(
	w_{(i,\sv{ms}+m)}^2(1+i)^{-2\bb_1}(1+\frac{\sv{ms}+1}{m})^{-2\bb_2} \right)^{\frac{2-\al}{2}}}
\end{equation*}
\begin{equation*}
\leq
w_{(i+n,\sv{ms})}w_{(i,\sv{ms}+m)}^{\al-1}
	\left(\frac{1+\sv{ms}}{m}\right)^{-\bb_2} (1+i)^{-\bb_1(\al-1)}\left(1+\frac{\sv{ms}+1}{m}\right)^{-\bb_2(\al-1)}
\end{equation*}
\begin{equation*}
\leq
e^\al
 (1+i)^{-\bb_1(\al-1)}s^{-\bb_2}\left(1+s\right)^{-\bb_2(\al-1)}
\end{equation*}

In the case under investigation we have $\bb_2<1,\  \bb_2>\frac{1}{\al},\ \bb_1>\frac{1}{\al-1}$, so
\begin{equation*}
\sum_{i=0}^{\infty}
\int\limits_{0}^{\infty}(1+i)^{-\bb_1(\al-1)}s^{-\bb_2}\left(1+s\right)^{-\bb_2(\al-1)} \dd s<\infty.
\end{equation*}

The function $y_{i,s,n,m}$ is dominated by an integrable function and converges pointwise. Application of The Dominated Convergence Theorem gives
\begin{equation*}
\frac{\rho(n,-m)}{n^{-\bb_1}m^{1-\bb_2\al}}=
\sum_{j=0}^{\infty}\int\limits_{0}^{\infty}
y_{i,s,n,m}
\dd s\rightarrow
\sum_{j=0}^{\infty}\int\limits_{0}^{\infty}
w_{(i,\infty)}^{\al-1}(1+i)^{-\bb_1(\al-1)}s^{-\bb_2}(1+s)^{-\bb_2(\al-1)}\dd s.
\end{equation*}

\subsection{The case $ 1<\al< 2,\ \bb_1=\frac{1}{\al-1}$ and $1<\bb_2<\frac{1}{\al-1}$;}

 $\bullet$ Say $m_n$ is a sequence such that $h_{n}=\frac{m_n^{-\bb_2}}{n^{-\bb_1}}\rightarrow c \in(0;\infty)$.

 \begin{equation*}
 \rho(n,-m_n)=\sum_{i=0}^{\infty}\sum_{j=0}^{\infty}
 \frac{
 	w_{(i+n,j)}w_{(i,j+m_n)}(1+i)^{-\bb_1}(1+j)^{-\bb_2}(1+i+n)^{-\bb_1}(1+j+m_n)^{-\bb_2}
 }{\left(
 w_{(i+n,j)}^2(1+i+n)^{-2\bb_1}(1+j)^{-2\bb_2}
 +
 w_{(i,j+m_n)}^2(1+i)^{-2\bb_1}(1+j+m_n)^{-2\bb_2}
 \right)^{\frac{2-\al}{2}}}
 \end{equation*}
  \begin{equation*}
 =\sum_{i=0}^{\infty}\sum_{j=0}^{\infty}
  \frac{
  	w_{(i+n,j)}w_{(i,j+m_n)}(1+i)^{-\bb_1}(1+j)^{-\bb_2}(1+\frac{i+1}{n})^{-\bb_1}(1+\frac{j+1}{m_n})^{-\bb_2}n^{-\bb_1}m_n^{-\bb_2(\al-1)}
  }{
  \left(
  w_{(i+n,j)}^2(1+\frac{i+1}{n})^{-2\bb_1}(1+j)^{-2\bb_2}\frac{n^{-2\bb_1}}{m_n^{-2\bb_2}}
  +
  w_{(i,j+m_n)}^2(1+i)^{-2\bb_1}(1+\frac{j+1}{m_n})^{-2\bb_2}
  \right)^{\frac{2-\al}{2}}}
\end{equation*}
\begin{equation*}
=\sum_{i=0}^{\infty}\sum_{j=0}^{\infty}n^{-\bb_1}m_n^{-\bb_2(\al-1)}q_{i,j,n},
\end{equation*}
 where
 \begin{equation*}
q_{i,j,n}= \frac{
 	w_{(i+n,j)}w_{(i,j+m_n)}(1+i)^{-\bb_1}(1+j)^{-\bb_2}(1+\frac{i+1}{n})^{-\bb_1}(1+\frac{j+1}{m_n})^{-\bb_2}
 }{
 \left(
 w_{(i+n,j)}^2(1+\frac{i+1}{n})^{-2\bb_1}(1+j)^{-2\bb_2}\frac{n^{-2\bb_1}}{m_n^{-2\bb_2}}
 +
 w_{(i,j+m_n)}^2(1+i)^{-2\bb_1}(1+\frac{j+1}{m_n})^{-2\bb_2}
 \right)^{\frac{2-\al}{2}}}.
 \end{equation*}

 Let us split $\sum_{i=0}^{\infty}\sum_{j=0}^{\infty}q_{i,j,n}$ into three sums:
 \begin{equation}\label{3.4rhosplit}\begin{gathered}
 \frac{\rho(n,-m_n)}{n^{-\bb_1}m_n^{-\bb_2(\al-1)}}
 =
 \sum_{i=0}^{n-1}\sum_{j=0}^{m_n-1}q_{i,j,n}+
 \sum_{i=n}^{\infty}\sum_{j=0}^{\infty}q_{i,j,n}+
 \sum_{i=0}^{n-1}\sum_{j=m_n}^{\infty}q_{i,j,n}=:Z_1+Z_2+Z_3.
 \end{gathered}\end{equation}

 At first we will show that ${Z_i}\rightarrow 0,\ n\rightarrow\infty,\ i=2,3.$

 \begin{equation*}
 |Z_2|\leq \sum_{i=n}^{\infty}\sum_{j=0}^{\infty}\abs{q_{i,j,n}}
 \end{equation*}
 \begin{equation*}\small
 \leq
 \sum_{i=n}^{\infty}\sum_{j=0}^{\infty}
 e^\al(1+i)^{-\bb_1}(1+j)^{-\bb_2(\al-1)}\left(1+\frac{i+1}{n}\right)^{-\bb_1(\al-1)}\left(1+\frac{j+1}{m_n}\right)^{-\bb_2} \frac{n^{-\bb_1(\al-2)}}{m_n^{-\bb_2(\al-2)}}
 \end{equation*}
 \begin{equation*}\small
 \leq e^{\al} \frac{n^{-\bb_1(\al-2)}}{m_n^{-\bb_2(\al-2)}}
 \int_{n}^{\infty}\int_{0}^{\infty}
 (1+\sv{t})^{-\bb_1}(1+\sv{s})^{-\bb_2(\al-1)}\left(1+\frac{\sv{t}+1}{n}\right)^{-\bb_1(\al-1)}\left(1+\frac{\sv{s}+1}{m_n}\right)^{-\bb_2}\dd t \dd s
 \end{equation*}
 \begin{equation*}
 \leq e^{\al} \frac{n^{-\bb_1(\al-2)}}{m_n^{-\bb_2(\al-2)}}
 \int_{n}^{\infty}t^{-\bb_1}\left(1+\frac{t}{n}\right)^{-\bb_1(\al-1)}  \dd t
 \int_{0}^{\infty} s^{-\bb_2(\al-1)}\left(1+\frac{s}{m_n}\right)^{-\bb_2} \dd s
 \end{equation*}
 \begin{equation*}
 = e^{\al} \frac{n^{-\bb_1(\al-2)}}{m_n^{-\bb_2(\al-2)}}
 \int_{1}^{\infty}n^{-\bb_1}t^{-\bb_1}\left(1+\frac{nt}{n}\right)^{-\bb_1(\al-1)}  \dd nt
 \int_{0}^{\infty} m_n^{-\bb_2(\al-1)}s^{-\bb_2(\al-1)}\left(1+\frac{m_ns}{m_n}\right)^{-\bb_2} \dd m_ns
 \end{equation*}
  \begin{equation*}
  = e^{\al}{n^{1-\bb_1(\al-1)}m_n^{1-\bb_2}}
  \int_{1}^{\infty}t^{-\bb_1}\left(1+t\right)^{-\bb_1(\al-1)}  \dd t
  \int_{0}^{\infty} s^{-\bb_2(\al-1)}\left(1+s\right)^{-\bb_2} \dd s
  \end{equation*}
  \begin{equation}\label{3.4Z2}
  = e^{\al}m_n^{1-\bb_2}
  \int_{1}^{\infty}t^{-\bb_1}\left(1+t\right)^{-\bb_1(\al-1)}  \dd t
  \int_{0}^{\infty} s^{-\bb_2(\al-1)}\left(1+s\right)^{-\bb_2} \dd s.
  \end{equation}
 In the case under consideration we have $\bb_1>\frac{1}{\al}$ and $1<\bb_2<\frac{1}{\al-1}$. Thus, the integrals in \eqref{3.4Z2} are finite and $m_n^{1-\bb_2}\rightarrow0$. This implies $|Z_2|\rightarrow0,$ as $n\rightarrow\infty$.

 We continue with $Z_3$:
 \begin{equation*}
 |Z_3|\leq \sum_{i=0}^{n-1}\sum_{j=m_n}^{\infty}|q_{i,j,n}|
 \end{equation*}
 \begin{equation*}
 \leq
 e^\al \sum_{i=0}^{n-1}\sum_{j=m_n}^{\infty}
 (1+i)^{-\bb_1(\al-1)}(1+j)^{-\bb_2}\left(1+\frac{i+1}{n}\right)^{-\bb_1}\left(1+\frac{j+1}{m_n}\right)^{-\bb_2(\al-1)}
 \end{equation*}
 \begin{equation*}
 \leq
 e^\al
 \sum_{i=0}^{n-1}(1+i)^{-\bb_1(\al-1)}
 \sum_{j=m_n}^{\infty}
 (1+j)^{-\bb_2}\left(1+\frac{j+1}{m_n}\right)^{-\bb_2(\al-1)}
 \end{equation*}
  \begin{equation*}
  =
  e^\al
  \left(1+ \sum_{i=1}^{n-1}(1+i)^{-\bb_1(\al-1)} \right)
  \sum_{j=m_n}^{\infty}
  (1+j)^{-\bb_2}\left(1+\frac{j+1}{m_n}\right)^{-\bb_2(\al-1)}
  \end{equation*}
   \begin{equation*}
   =
   e^\al
   \left(1+ \int_{1}^{n}(1+\sv{t})^{-\bb_1(\al-1)} \dd t\right)
   \int_{m_n}^{\infty}
   (1+\sv{s})^{-\bb_2}\left(1+\frac{\sv{s}+1}{m_n}\right)^{-\bb_2(\al-1)} \dd s
   \end{equation*}
 \begin{equation*}
 \leq
 e^\al
 \left(1+ \int_{1}^{n}t^{-\bb_1(\al-1)} \dd t\right)
 \int_{m_n}^{\infty}
 s^{-\bb_2}\left(1+\frac{s}{m_n}\right)^{-\bb_2(\al-1)} \dd s
 \end{equation*}
  \begin{equation*}
  =
  e^\al
  \left(1+ \ln n \right)
  m_n^{1-\bb_2} \int_{1}^{\infty}
  s^{-\bb_2}\left(1+s\right)^{-\bb_2(\al-1)} \dd s.
  \end{equation*}
 The integral is finite, since $\bb_2>\frac{1}{\al}$. We have $\left(1+ \ln n \right)
 m_n^{1-\bb_2}=\left(1+ \ln n \right) n^{-\bb_1\frac{\bb_2-1}{\bb_2}}\left( \frac{m_n^{-\bb_2}}{n^{-\bb_1}} \right)^{\frac{\bb_2-1}{\bb_2}} \rightarrow0$ as $n\rightarrow\infty$, so this inequality implies $Z_3\rightarrow0$.

 We proceed with $Z_1$:
 \begin{equation*}
 Z_1=\sum_{i=0}^{n-1}\sum_{j=0}^{m_n-1}q_{i,j,n}=
 \sum_{i=0}^{\infty}\sum_{j=0}^{\infty}q_{i,j,n}\ind{[i<n,j<m_n]}.
 \end{equation*}
 For fixed values of $i,j$ we have
 \begin{equation*}
 q_{i,j,n}\ind{[i<n,j<m_n]}\rightarrow
 \frac
 {w_{(\infty,j)}w_{(i,\infty)}(1+j)^{-\bb_2}(1+i)^{-\bb_1}}
 {\left(
 	w_{(\infty,j)}^2c^{-2}(1+j)^{-2\bb_2}+
 	w_{(i,\infty)}^2(1+i)^{-2\bb_1}
 	\right)^{\frac{2-\al}{2}}} ,\ n\rightarrow\infty.
 \end{equation*}
 Moreover
 \begin{equation*}
 \abs{q_{i,j,n}}\ind{[i<n,j<m_n]}
 \end{equation*}
 \begin{equation*}
 = \frac{
 	w_{(i+n,j)}w_{(i,j+m_n)}(1+i)^{-\bb_1}(1+j)^{-\bb_2}(1+\frac{i+1}{n})^{-\bb_1}(1+\frac{j+1}{m_n})^{-\bb_2}\ind{[i<n,j<m_n]}
 }{
 \left(
 w_{(i+n,j)}^2(1+\frac{i+1}{n})^{-2\bb_1}(1+j)^{-2\bb_2}\frac{n^{-2\bb_1}}{m_n^{-2\bb_2}}
 +
 w_{(i,j+m_n)}^2(1+i)^{-2\bb_1}(1+\frac{j+1}{m_n})^{-2\bb_2}
 \right)^{\frac{2-\al}{2}}}
 \end{equation*}
  \begin{equation*}
 \leq \frac{
 	e^2(1+i)^{-\bb_1}(1+j)^{-\bb_2}
 }{
 \left(
 d^22^{-2\bb_1}(1+j)^{-2\bb_2}\frac{n^{-2\bb_1}}{m_n^{-2\bb_2}}
 +
 d^22^{-2\bb_2}(1+i)^{-2\bb_1}
 \right)^{\frac{2-\al}{2}}}.
 \end{equation*}

 The sequence $m_n$ satisfies $\frac{m_n^{-\bb_2}}{n^{-\bb_1}}\rightarrow c \in(0;\infty)$ as $n\rightarrow\infty$ so there exists $N\in\NN$ such that $\frac{n^{-\bb_1}}{m_n^{-\bb_2}}\geq\frac{c^{-1}}{2}>0$ for $n>N$. Thus, for $n>N$ we have
 \begin{equation*}
 \abs{q_{i,j,n}}\ind{[i<n,j<m_n]}\leq
 \frac{
 	e^2(1+i)^{-\bb_1}(1+j)^{-\bb_2}
 }{
 \left(
 d^22^{-2\bb_1}(1+j)^{-2\bb_2}\frac{c^{-1}}{2}
 +
 d^22^{-2\bb_2}(1+i)^{-2\bb_1}
 \right)^{\frac{2-\al}{2}}}
 \end{equation*}
 \begin{equation*}
 \leq K
 \frac{
 	(1+i)^{-\bb_1}(1+j)^{-\bb_2}
 }{
 \left(
 (1+j)^{-2\bb_2}
 +
 (1+i)^{-2\bb_1}
 \right)^{\frac{2-\al}{2}}},
\end{equation*}
 with $K=e^2\left( d^2\min(2^{-2\bb_1-1}c^{-1}, 2^{-2\bb_2}) \right)^{\frac{\al-2}{2}}$.

 We will now show that
 \begin{equation}\label{3.4finalsum}
 \sum_{i=0}^{\infty}\sum_{j=0}^{\infty}
 \frac
 {(1+j)^{-\bb_2}(1+i)^{-\bb_1}}
 {\left(
 	(1+j)^{-2\bb_2}+
 	(1+i)^{-2\bb_1}
 	\right)^{\frac{2-\al}{2}}}<\infty.
 \end{equation}

 We will use notation
 \begin{equation*}
 a_{i,j}=\frac
 {(1+j)^{-\bb_2}(1+i)^{-\bb_1}}
 {\left(
 	(1+j)^{-2\bb_2}+
 	(1+i)^{-2\bb_1}
 	\right)^{\frac{2-\al}{2}}}.
 \end{equation*}

 We split the sum in \eqref{3.4finalsum} as
 \begin{equation*}
 \sum_{i=0}^{\infty}\sum_{j=0}^{\infty}a_{i,j}=
  \sum_{i=0}^{\infty}a_{i,0}+
   \sum_{j=1}^{\infty}a_{0,j}+
   \sum_{i=1}^{\infty}\sum_{j=1}^{\infty}a_{i,j}.
 \end{equation*}

 Now
 \begin{equation*}
 \sum_{i=0}^{\infty}a_{i,0}\leq \sum_{i=0}^{\infty}(1+i)^{-\bb_1}<\infty,
 \end{equation*}
 \begin{equation*}
 \sum_{j=1}^{\infty}a_{0,j}\leq \sum_{j=1}^{\infty}(1+j)^{-\bb_2}<\infty,
 \end{equation*}
 and it remains to show that
 \begin{equation*}
 \sum_{i=1}^{\infty}\sum_{j=1}^{\infty}a_{i,j}<\infty.
 \end{equation*}

 Since
 \begin{equation*}\begin{gathered}
 \sum_{i=1}^{\infty}\sum_{j=1}^{\infty}a_{i,j}=
 \int_{1}^{\infty}\int_{1}^{\infty}\frac
 {
 	(1+\sv{s})^{-\bb_2}
 	(1+\sv{t})^{-\bb_1}
 }
 {\left(
 	(1+\sv{s})^{-2\bb_2}+
 	(1+\sv{t})^{-2\bb_1}
 	\right)^{\frac{2-\al}{2}}}    \dd t \dd s   \\ \leq
 \int_{1}^{\infty}\int_{1}^{\infty}\frac
 {
 	s^{-\bb_2}
 	t^{-\bb_1}
 }
 {\left(
 	(1+s)^{-2\bb_2}+
 	(1+t)^{-2\bb_1}
 	\right)^{\frac{2-\al}{2}}}    \dd t \dd s,
 \end{gathered}\end{equation*}
 applying the inequality $(1+x)^{-1}\geq 2^{-1}x^{-1},\ x\geq 1$ we obtain
 \begin{equation*}\begin{gathered}
 \sum_{i=1}^{\infty}\sum_{j=1}^{\infty}a_{i,j}\leq
 \int_{1}^{\infty}\int_{1}^{\infty}\frac
 {
 	s^{-\bb_2}
 	t^{-\bb_1}
 }
 {\left(
 	2^{-2\bb_2}s^{-2\bb_2}+
 	2^{-2\bb_1}t^{-2\bb_1}
 	\right)^{\frac{2-\al}{2}}}    \dd t \dd s \\
 \leq
 \int_{1}^{\infty}\int_{1}^{\infty}\min(2^{-2\bb_1},2^{-2\bb_2})^{\frac{\al-2}{2}}\frac
 {
 	s^{-\bb_2}
 	t^{-\bb_1}
 }
 {\left(
 	s^{-2\bb_2}+
 	t^{-2\bb_1}
 	\right)^{\frac{2-\al}{2}}}    \dd t \dd s.
 \end{gathered}\end{equation*}

 {The change of variables $x=t^{-\bb_1},\ y=s^{-\bb_2}$ and use of polar coordinates yield}
 \begin{equation}\label{intI}  \begin{gathered}
 \int_{1}^{\infty}\int_{1}^{\infty}\frac
 {
 	s^{-\bb_2}
 	t^{-\bb_1}
 }
 {\left(
 	s^{-2\bb_2}+
 	t^{-2\bb_1}
 	\right)^{\frac{2-\al}{2}}}    \dd t \dd s=
 \frac{1}{\bb_1\bb_2}
 \int_{0}^{1}\int_{0}^{1}\frac
 {
 	x^{-\frac{1}{\bb_1}}
 	y^{-\frac{1}{\bb_2}}
 }
 {\left(
 	x^2+ y^2
 	\right)^{\frac{2-\al}{2}}}    \dd x \dd y
 \\ \leq
 \frac{1}{\bb_1\bb_2}\int_{0}^{1}\int_{0}^{\frac{\pi}{2}}\frac
 {
 	r^{1-\frac{1}{\bb_1}-\frac{1}{\bb_2}} \cos^{-\frac{1}{\bb_1}}(\phi)
 	\sin^{-\frac{1}{\bb_2}}(\phi)
 }
 {\left(
 	r^2
 	\right)^{\frac{2-\al}{2}}}    \dd \phi \dd r  \\=
 \frac{1}{\bb_1\bb_2}
 \int_{0}^{1} r^{\al-1-\frac{1}{\bb_1}-\frac{1}{\bb_2}}  \dd r
 \int_{0}^{\frac{\pi}{2}}
 \cos^{-\frac{1}{\bb_1}}(\phi)
 \sin^{-\frac{1}{\bb_2}}(\phi)
 \dd \phi .
 \end{gathered}\end{equation}
 The integrals on the right hand side are finite since $\al-1-\frac{1}{\bb_1}-\frac{1}{\bb_2}>-1$ and $-\frac{1}{\bb_i}>-1,\ i=1,2$.

 Thus, the function $q_{i,j,n}\ind{[i<n,j<m_n]}$ converges pointwise and is dominated by an integrable function. The Dominated Convergence Theorem together with \eqref{3.4rhosplit} implies
 \begin{equation}\label{rezmnc}\begin{gathered}
 \frac{ \rho{(n,-m)}}{n^{-\bb_1}m_n^{-\bb_2(\al-1)}}\rightarrow
 \sum_{i=0}^{\infty}\sum_{j=0}^{\infty}
 \frac
 {w_{(\infty,j)}w_{(i,\infty)}(1+j)^{-\bb_2}(1+i)^{-\bb_1}}
 {\left(
 	w_{(\infty,j)}^2c^{-2}(1+j)^{-2\bb_2}+
 	w_{(i,\infty)}^2(1+i)^{-2\bb_1}
 	\right)^{\frac{2-\al}{2}}} ,\ n\rightarrow\infty.
 \end{gathered}\end{equation}
 \medskip

 $\bullet$ Let us now assume that $m_n$ is a sequence such that $h_{n}=\frac{m_n^{-\bb_2}}{n^{-\bb_1}}\rightarrow 0$.

 Fix $\varepsilon>0$ and assume $n,m_n\geq\frac{2}{\varepsilon}$. Then the following inequalities hold:
 \begin{equation*}
 \frac{\varepsilon}{2}\leq
 \varepsilon-\frac{1}{n}=
 \frac{n\varepsilon-1}{n}\leq
 \frac{\sv{n\varepsilon}}{n}\leq
 \frac{n\varepsilon}{n}=
 \varepsilon,
 \end{equation*}
 \begin{equation*}
 \frac{\varepsilon}{2}\leq
 \frac{\sv{m_n\varepsilon}}{m_n}\leq
 \varepsilon.
 \end{equation*}

 Denote
 \begin{equation}\label{Wzymejimas}\begin{gathered}
 W_{i,j,n}=\frac
 {w_{(i+n,j)}w_{(i,j+m_n)}(1+i+n)^{-\bb_1}(1+j)^{-\bb_2}(1+i)^{-\bb_1}(1+j+m_n)^{-\bb_2}}
 {\left(
 	w^2_{(i+n,j)}(1+i+n)^{-2\bb_1}(1+j)^{-2\bb_2}+
 	w^2_{(i,j+m_n)}(1+i)^{-2\bb_1}(1+j+m_n)^{-2\bb_2} \right)^{\frac{2-\al}{2}}}
 \end{gathered}\end{equation}
 and split
 \begin{equation}\label{rhoSkaid}\begin{gathered}
 \rho({n,-m_n})=\sum_{i=0}^{\infty}\sum_{j=0}^{\infty}W_{i,j,n}=D_1+D_2+D_3+D_4,
 \end{gathered}\end{equation}
 where
 \begin{equation}\label{rhoSkaidDalys}\begin{gathered}
 D_1= \sum_{i=\sv{n\varepsilon}}^{\infty}\sum_{j=\sv{m_n\varepsilon}}^{\infty}W_{i,j,n},\ \ D_2= \sum_{i=0}^{\sv{n\varepsilon}-1}\sum_{j=\sv{m_n\varepsilon}}^{\infty}W_{i,j,n},\\
 D_3= \sum_{i=\sv{n\varepsilon}}^{\infty}\sum_{j=0}^{\sv{m_n\varepsilon}-1}W_{i,j,n},\ \  D_4= \sum_{i=0}^{\sv{n\varepsilon}-1}\sum_{j=0}^{\sv{m_n\varepsilon}-1}W_{i,j,n}.
 \end{gathered}\end{equation}
The absolute value of $D_1$ can be bounded from above in the following way:
 \begin{equation*}\begin{gathered}
 |D_1|\leq
 \frac{e^2}{d^{2-\al}} \sum_{i=\sv{n\varepsilon}}^{\infty}\sum_{j=\sv{m_n\varepsilon}}^{\infty}\frac
 {(1+i+n)^{-\bb_1}(1+j)^{-\bb_2}(1+i)^{-\bb_1}(1+j+m_n)^{-\bb_2}}
 {\left(
 	(1+i+n)^{-2\bb_1}(1+j)^{-2\bb_2}+
 	(1+i)^{-2\bb_1}(1+j+m_n)^{-2\bb_2} \right)^{\frac{2-\al}{2}}}\\ \leq
 \frac{e^2}{d^{2-\al}}\sum_{i=\sv{n\varepsilon}}^{\infty}\sum_{j=\sv{m_n\varepsilon}}^{\infty}
 (1+i)^{-\bb_1}(1+j+m_n)^{-\bb_2}(1+j)^{-\bb_2(\al-1)}(1+i+n)^{-\bb_1(\al-1)}\\
 \leq
 \frac{e^2}{d^{2-\al}}\sum_{i=\sv{n\varepsilon}}^{\infty}\sum_{j=\sv{m_n\varepsilon}}^{\infty}
 (1+i)^{-\bb_1\al}(1+j)^{-\bb_2\al}\\
 =
 \frac{e^2}{d^{2-\al}}
 \int_{\sv{n\varepsilon}}^{\infty}\int_{\sv{m_n\varepsilon}}^{\infty}
 (1+\sv{t})^{-\bb_1\al}(1+\sv{s})^{-\bb_2\al}\dd s \dd t
 \\
 \leq
 \frac{e^2}{d^{2-\al}}
 \int_{\sv{n\varepsilon}}^{\infty}\int_{\sv{m_n\varepsilon}}^{\infty}
 t^{-\bb_1\al}s^{-\bb_2\al}\dd s \dd t
 \end{gathered}\end{equation*}
 \begin{equation*}\begin{gathered}
 =\frac{e^2}{d^{2-\al}} n^{1-\bb_1\al}m_n^{1-\bb_2\al}
 \int_{\frac{\sv{n\varepsilon}}{n}}^{\infty}\int_{\frac{\sv{m_n\varepsilon}}{m_n}}^{\infty}
 t^{-\bb_1\al}s^{-\bb_2\al}
 \dd s \dd t
 \end{gathered}\end{equation*}
 \begin{equation*}\begin{gathered}
\leq \frac{e^2}{d^{2-\al}} n^{1-\bb_1\al}m_n^{1-\bb_2\al}
\int_{\frac{\varepsilon}{2}}^{\infty}
t^{-\bb_1\al}\dd t
\int_{\frac{\varepsilon}{2}}^{\infty}
s^{-\bb_2\al}
\dd s .
 \end{gathered}\end{equation*}
 The integrals on the right hand side are finite. Thus,
 \begin{equation}\label{D1greitis}
 \frac{D_1}{n^{1-\bb_1\al}m_n^{1-\bb_2\al}}\leq C_1\in\RR.
 \end{equation}

 We proceed with $D_2$. Let us split
 \begin{equation*}
 D_2= \sum_{i=0}^{\sv{n\varepsilon}-1}\sum_{j=\sv{m_n\varepsilon}}^{\infty}W_{i,j,n}=\tilde{D}_2+\sum_{j=\sv{m_n\varepsilon}}^{\infty}W_{0,j,n},
 \end{equation*}
 where
 \begin{equation*}
 \tilde{D}_2=\sum_{i=1}^{\sv{n\varepsilon}-1}\sum_{j=\sv{m_n\varepsilon}}^{\infty}W_{i,j,n}.
 \end{equation*}
 We estimate the second term
 \begin{equation}\label{D2krastas}\begin{gathered}
 \abs{\sum_{j=\sv{m_n\varepsilon}}^{\infty}W_{0,j,n}}\leq
 \frac{e^2}{d^{2-\al}}
 \sum_{j=\sv{m_n\varepsilon}}^{\infty}
 \frac
 {(1+n)^{-\bb_1}(1+j)^{-\bb_2}(1+j+m_n)^{-\bb_2}}
 {\left(
 	(1+n)^{-2\bb_1}(1+j)^{-2\bb_2}+
 	(1+j+m_n)^{-2\bb_2} \right)^{\frac{2-\al}{2}}}\\ \leq
 \frac{e^2}{d^{2-\al}}
 n^{-\bb_1}\sum_{j=\sv{m_n\varepsilon}}^{\infty}
 {(1+j)^{-\bb_2}(1+j+m_n)^{-\bb_2(\al-1)}}\\
 \leq
 \frac{e^2}{d^{2-\al}}
 n^{-\bb_1}\sum_{j=\sv{m_n\varepsilon}}^{\infty}
 {(1+j)^{-\bb_2\al}}\\
 =
 \frac{e^2}{d^{2-\al}}
 n^{-\bb_1}
 \int_{\sv{m_n\varepsilon}}^{\infty}
 (1+\sv{s})^{-\bb_2\al}\dd s\leq
 \frac{e^2}{d^{2-\al}}
 n^{-\bb_1}
 \int_{\sv{m_n\varepsilon}}^{\infty}
 s^{-\bb_2\al}\dd s\\
 =
  \frac{e^2}{d^{2-\al}}
  n^{-\bb_1}m_n ^{1-\bb_2\al}
  \int_{\frac{\sv{m_n\varepsilon}}{m_n}}^{\infty}
  s^{-\bb_2\al}\dd s
  \leq
  \frac{e^2}{d^{2-\al}}n^{-\bb_1}m_n ^{1-\bb_2\al}
 \int_{\frac{\varepsilon}{2}}^{\infty}
 {s^{-\bb_2}}\dd s,
 \end{gathered}\end{equation}
 while $\tilde{D}_2$ can be estimated as follows:
 \begin{equation*}\begin{gathered}
 \abs{\tilde{D}_2}\leq \sum_{i=1}^{\sv{n\varepsilon}-1}\sum_{j=\sv{m_n\varepsilon}}^{\infty}\abs{W_{i,j,n}}\\
 \leq
 \frac{e^2}{d^{2-\al}} \sum_{i=1}^{\sv{n\varepsilon}-1}\sum_{j=\sv{m_n\varepsilon}}^{\infty}
 \frac
 {(1+i+n)^{-\bb_1}(1+j)^{-\bb_2}(1+i)^{-\bb_1}(1+j+m_n)^{-\bb_2}}
 {\left(
 	(1+i+n)^{-2\bb_1}(1+j)^{-2\bb_2}+
 	(1+i)^{-\bb_1}(1+j+m_n)^{-2\bb_2} \right)^{\frac{2-\al}{2}}}  \\ \leq
 \frac{e^2}{d^{2-\al}}  \sum_{i=1}^{\sv{n\varepsilon}-1}\sum_{j=\sv{m_n\varepsilon}}^{\infty}
 {(1+i)^{-\bb_1(\al-1)}(1+j+m_n)^{-\bb_2(\al-1)}(1+i+n)^{-\bb_1}(1+j)^{-\bb_2}} \\
  \leq
 \frac{e^2}{d^{2-\al}}
  \sum_{i=1}^{\sv{n\varepsilon}-1}(1+i)^{-\bb_1(\al-1)}(1+i+n)^{-\bb_1}
  \sum_{j=\sv{m_n\varepsilon}}^{\infty}(1+j)^{-\bb_2\al}\\
  \leq
  \frac{e^2}{d^{2-\al}} n^{-\bb_1}m_n^{1-\bb_2\al}
  \sum_{i=1}^{\sv{n\varepsilon}-1}(1+i)^{-\bb_1(\al-1)}
  \int_{\frac{\varepsilon}{2}}^{\infty}
  s^{-\bb_2\al}
  \dd s\\
  \leq
  \frac{e^2}{d^{2-\al}} n^{-\bb_1}m_n^{1-\bb_2\al}
  \int_{1}^{n\varepsilon}t^{-\bb_1(\al-1)}\dt
  \int_{\frac{\varepsilon}{2}}^{\infty}
  s^{-\bb_2\al}
  \dd s\\
  =
  \frac{e^2}{d^{2-\al}} n^{-\bb_1}m_n^{1-\bb_2\al}
  \left( \ln(\varepsilon)+\ln (n) \right)
  \int_{\frac{\varepsilon}{2}}^{\infty}
  s^{-\bb_2\al}
  \dd s  .
 \end{gathered}\end{equation*}
 The obtained inequalities imply that, for some constant $C_2$,
 \begin{equation*}
 \frac{|D_2|}{n^{-\bb_1}m_n^{1-\bb_2\al}\ln n}\leq C_2.
 \end{equation*}
  In a similar way we can obtain the following inequality
 \begin{equation}\label{L43D3ivertis}\begin{gathered}
 |D_3|\leq
 \frac{e^2}{d^{2-\al}} n^{-\bb_1}m_n^{1-\bb_2\al}
 \int_{\frac{\varepsilon}{2}}^{\infty}
 t^{-\bb_1\al}
 \dd t
 \int_{0}^{\varepsilon}
 \left(1+s\right)^{-\bb_2}
 s^{-\bb_2(\al-1)}
 \dd s ,
 \end{gathered}\end{equation}
 where the second integral is finite since $\bb_2<\frac{1}{\al-1}$.
 We see that there exists a constant $K$ such that,  for $n$ satisfying $n,m_n\geq\frac{2}{\varepsilon}$,
 \begin{equation}\label{D123}
 \frac{\abs{D_1+D_2+D_3}}{n^{-\bb_1}m_n^{1-\bb_2\al}\ln (n)} \leq K.
 \end{equation}

 Let us choose $\varepsilon=1$ and split $D_4$ as
 \begin{equation*}
 D_4=
 \sum_{i=0}^{n-1}\sum_{j=1}^{m_n-1}W_{i,j,n}+
 \sum_{i=0}^{n-1}W_{i,0,n}=
 \tilde{D}_4+\sum_{i=0}^{n-1}W_{i,0,n}.
 \end{equation*}
 We will show that
 \begin{equation}\label{D4kr}
 \frac{\sum_{i=0}^{n-1}W_{i,0,n}}{n^{-\frac{\bb_1}{\bb_2}}m_n^{1-\bb_2\al}}\rightarrow 0,\ n\rightarrow\infty.
 \end{equation}

 \begin{equation}\label{D4krastas}\begin{gathered}
 \frac{1}{n^{-\frac{\bb_1}{\bb_2}}m_n^{1-\bb_2\al}}\abs{\sum_{i=0}^{n-1}W_{i,0,n}}\\
 \leq
 {n^{\frac{\bb_1}{\bb_2}}m_n^{\bb_2\al-1}}
 \frac{e^2}{d^{2-\al}}
 \sum_{i=0}^{n-1}
 \frac
 {(1+i+n)^{-\bb_1}(1+i)^{-\bb_1}(1+m_n)^{-\bb_2}}
 {\left(
 	(1+i+n)^{-2\bb_1}+
 	(1+i)^{-\bb_1}(1+m_n)^{-2\bb_2} \right)^{\frac{2-\al}{2}}}
 \end{gathered}\end{equation}
  \begin{equation*}\begin{gathered}
  \leq {n^{\frac{\bb_1}{\bb_2}}m_n^{\bb_2\al-1}}
 \frac{e^2}{d^{2-\al}} \sum_{i=0}^{n-1}
 {(1+i+n)^{-\bb_1(\al-1)}(1+i)^{-\bb_1}(1+m_n)^{-\bb_2}}\\ \leq
 \frac{e^2}{d^{2-\al}}
 m_n^{\bb_2\al-1-\bb_2}n^{\frac{\bb_1}{\bb_2}-\bb_1(\al-1)}
 \sum_{i=0}^{\infty}
 (1+i)^{-\bb_1}\\=
  \frac{e^2}{d^{2-\al}}
  \left( m_nn^{-\frac{\bb_1}{\bb_2}} \right)^{\bb_2(\al-1)-1}
  \sum_{i=0}^{\infty}
  (1+i)^{-\bb_1}
 .
 \end{gathered}\end{equation*}
 Since $\bb_1=\frac{1}{\al-1}$ and $\al<2$, the sum is finite. In the case under consideration we have $\bb_2(\al-1)-1<0$, so $$\left( m_nn^{-\frac{\bb_1}{\bb_2}} \right)^{\bb_2(\al-1)-1}\rightarrow0,\ n\rightarrow\infty,$$
 which proves \eqref{D4kr}.

 It will be convenient to introduce notation $g_{n}= m_nn^{-\frac{\bb_1}{\bb_2}}$. We have $g_n\rightarrow\infty$ as $n\rightarrow\infty$. Then
 \begin{equation*}\footnotesize\begin{gathered}
 \tilde{D}_4= \sum_{i=0}^{n-1}\sum_{j=1}^{m_n-1}
 \frac
 {w_{(i+n,j)}w_{(i,j+m_n)}(1+i+n)^{-\bb_1}(1+j)^{-\bb_2}(1+i)^{-\bb_1}(1+j+m_n)^{-\bb_2}}
 {\left(
 	w^2_{(i+n,j)}(1+i+n)^{-2\bb_1}(1+j)^{-2\bb_2}+
 	w^2_{(i,j+m_n)}(1+i)^{-2\bb_1}(1+j+m_n)^{-2\bb_2} \right)^{\frac{2-\al}{2}}}\\=
 \sum_{i=0}^{n-1}\int_{1}^{m_n}
 \frac
 {w_{(i+n,\sv{s})}w_{(i,\sv{s}+m_n)}(1+i+n)^{-\bb_1}(1+\sv{s})^{-\bb_2}(1+i)^{-\bb_1}(1+\sv{s}+m_n)^{-\bb_2}}
 {\left(
 	w^2_{(i+n,\sv{s})}(1+i+n)^{-2\bb_1}(1+\sv{s})^{-2\bb_2}+
 	w^2_{(i,\sv{s}+m_n)}(1+i)^{-2\bb_1}(1+\sv{s}+m_n)^{-2\bb_2} \right)^{\frac{2-\al}{2}}}\dd s\\=
 \sum_{i=0}^{n-1}\int\limits_{\frac{1}{g_n}}^{n^{\frac{\bb_1}{\bb_2}}}
 \frac
 {g_nw_{(i+n,\sv{gs})}w_{(i,\sv{gs}+m_n)}
 	\left(1+i+n\right)^{-\bb_1}
 	\left(1+\sv{g_ns}\right)^{-\bb_2}
 	\left(1+i\right)^{-\bb_1}
 	\left(1+\sv{g_ns}+m_n\right)^{-\bb_2}}
 {\left(
 	w^2_{(i+n,\sv{g_ns})}
 	\left(1+i+n\right)^{-2\bb_1}
 	\left(1+\sv{g_ns}\right)^{-2\bb_2}+
 	w^2_{(i,\sv{g_ns}+m_n)}
 	\left(1+i\right)^{-2\bb_1}
 	\left(1+\sv{g_ns}+m_n\right)^{-2\bb_2}
 	\right)^{\frac{2-\al}{2}}}\dd s \\=
 \sum_{i=0}^{n-1}\int\limits_{\frac{1}{g_n}}^{n^{\frac{\bb_1}{\bb_2}}}
 \frac
 {g_n^{1-\bb_2(\al-1)}n^{-\bb_1(\al-1)}m_n^{-\bb_2}w_{(i+n,\sv{gs})}w_{(i,\sv{gs}+m_n)}
 	\left(1+\frac{i+1}{n}\right)^{-\bb_1}
 	\left(\frac{1+\sv{g_ns}}{g_n}\right)^{-\bb_2}
 	\left(1+i\right)^{-\bb_1}
 	\left(1+\frac{\sv{gs}+1}{m_n}\right)^{-\bb_2}}
 {\left(
 	w^2_{(i+n,\sv{g_ns})}
 	\left(1+\frac{i+1}{n}\right)^{-2\bb_1}
 	\left(\frac{1+\sv{g_ns}}{g_n}\right)^{-2\bb_2}+
 	w^2_{(i,\sv{g_ns}+m_n)}
 	\left(1+i\right)^{-2\bb_1}
 	\left(1+\frac{\sv{g_ns}+1}{m_n}\right)^{-2\bb_2}
 	\right)^{\frac{2-\al}{2}}}\dd s
  \end{gathered}\end{equation*}
 \begin{equation*}
 =
 m_n^{1-\bb_2\al}n^{-\frac{\bb_1}{\bb_2}}
 \sum_{i=0}^{n-1}\int_{\frac{1}{g_n}}^{n^{\frac{\bb_1}{\bb_2}}}
 u_{i,s,n}\dd s
  =
  m_n^{1-\bb_2\al}n^{-\frac{\bb_1}{\bb_2}}
  \sum_{i=0}^{\infty}\int_{0}^{\infty}
  u_{i,s,n}\ind{\left[i\leq n-1, \frac{1}{g_n}\leq s \leq n^{\frac{\bb_1}{\bb_2}}\right]}\dd s
 ,
 \end{equation*}

 where
 \begin{equation*}
 u_{i,s,n}=\frac
 {w_{(i+n,\sv{gs})}w_{(i,\sv{gs}+m_n)}
 	\left(1+\frac{i+1}{n}\right)^{-\bb_1}
 	\left(\frac{1+\sv{g_ns}}{g_n}\right)^{-\bb_2}
 	\left(1+i\right)^{-\bb_1}
 	\left(1+\frac{\sv{gs}+1}{m_n}\right)^{-\bb_2}}
 {\left(
 	w^2_{(i+n,\sv{g_ns})}
 	\left(1+\frac{i+1}{n}\right)^{-2\bb_1}
 	\left(\frac{1+\sv{g_ns}}{g_n}\right)^{-2\bb_2}+
 	w^2_{(i,\sv{g_ns}+m_n)}
 	\left(1+i\right)^{-2\bb_1}
 	\left(1+\frac{\sv{g_ns}+1}{m_n}\right)^{-2\bb_2}
 	\right)^{\frac{2-\al}{2}}}.
 \end{equation*}

 For fixed $i,s$ we have
 \begin{equation*}
 u_{i,s,n}\ind{\left[i\leq n-1, \frac{1}{g_n}\leq s \leq n^{\frac{\bb_1}{\bb_2}}\right]}\rightarrow\frac
 {w_{(i,\infty)}\left(1+i\right)^{-\bb_1}
 	s^{-\bb_2} }
 {\left(
 	w_{(i,\infty)}^2
 	\left(1+i\right)^{-2\bb_1}+
 	s^{-2\bb_2} 	
 	\right)^{\frac{2-\al}{2}}},
 \end{equation*}
 as $n\rightarrow\infty$. Moreover,
 \begin{equation*}\begin{gathered}
 u_{i,s,n}\ind{\left[i\leq n-1, \frac{1}{g_n}\leq s \leq n^{\frac{\bb_1}{\bb_2}}\right]} \\ \leq
  \frac{e^2}{d^{2-\al}}\frac
 {
 	\left(1+\frac{i+1}{n}\right)^{-\bb_1}
 	\left(\frac{1+\sv{gs}}{g}\right)^{-\bb_2}
 	\left(1+i\right)^{-\bb_1}
 	\left(1+\frac{\sv{gs}+1}{m_n}\right)^{-\bb_2}}
 {\left(
 	\left(1+\frac{i+1}{n}\right)^{-2\bb_1}
 	\left(\frac{1+\sv{gs}}{g}\right)^{-2\bb_2}+
 	\left(1+i\right)^{-2\bb_1}
 	\left(1+\frac{\sv{gs}+1}{m_n}\right)^{-2\bb_2}
 	\right)^{\frac{2-\al}{2}}}
 \end{gathered}\end{equation*}
 \begin{equation*}\begin{gathered}
 \leq
  \frac{e^2}{d^{2-\al}}\frac
 {
 	\left(\frac{1+\sv{gs}}{g}\right)^{-\bb_2}
 	\left(1+i\right)^{-\bb_1}}
 {\left(
 	2^{-2\bb_1}
 	\left(\frac{1+\sv{gs}}{g}\right)^{-2\bb_2}+
 	2^{-2\bb_2}\left(1+i\right)^{-2\bb_1}
 	\right)^{\frac{2-\al}{2}}}\\\leq
  \frac{e^2}{d^{2-\al}}\frac
 {
 	s^{-\bb_2}
 	\left(1+i\right)^{-\bb_1}}
 {\left(
 	2^{-2\bb_1}2^{-2\bb_2}s^{-2\bb_2}
 	+
 	2^{-2\bb_2}\left(1+i\right)^{-2\bb_1}
 	\right)^{\frac{2-\al}{2}}}\\\leq
  \frac{e^2}{d^{2-\al}}2^{\left( \bb_1+\bb_2 \right)\left( 2-\al \right)}
 \frac
 {
 	s^{-\bb_2}
 	\left(1+i\right)^{-\bb_1}}
 {\left(
 	s^{-2\bb_2}
 	+
 	\left(1+i\right)^{-2\bb_1}
 	\right)^{\frac{2-\al}{2}}}.
 \end{gathered}\end{equation*}

 The {dominating function on the right hand side is integrable}. The change of variables $x=\left(1+i\right)^{-\frac{\bb_1}{\bb_2}}s $ gives us
 \begin{equation*}\begin{gathered}
 \sum_{i=0}^{\infty}\int\limits_{0}^{\infty}
 \frac
 {
 	s^{-\bb_2}
 	\left(1+i\right)^{-\bb_1}}
 {\left(
 	s^{-2\bb_2}
 	+
 	\left(1+i\right)^{-2\bb_1}
 	\right)^{\frac{2-\al}{2}}}\dd s=
 \sum_{i=0}^{\infty}\int\limits_{0}^{\infty}
 \frac
 {
 	s^{-\bb_2}
 	\left(1+i\right)^{-\bb_1(\al-1)}}
 {\left(
 	\left(1+i\right)^{2\bb_1}s^{-2\bb_2}
 	+
 	1
 	\right)^{\frac{2-\al}{2}}}\dd s\\=
 \sum_{i=0}^{\infty}\int\limits_{0}^{\infty}
 \frac
 {
 	\left( \left(1+i\right)^{-\frac{\bb_1}{\bb_2}}s \right)^{-\bb_2}
 	\left(1+i\right)^{-\bb_1\al}}
 {\left(
 	\left( \left(1+i\right)^{-\frac{\bb_1}{\bb_2}}s \right)^{-2\bb_2}
 	+
 	1
 	\right)^{\frac{2-\al}{2}}}\dd s.
 \end{gathered}\end{equation*}
 \begin{equation*}\begin{gathered}
 =
 \sum_{i=0}^{\infty}\left(1+i\right)^{-\bb_1\al+\frac{\bb_1}{\bb_2}}
 \int\limits_{0}^{\infty}
 \frac
 {
 	x^{-\bb_2}
 }
 {\left(
 	x^{-2\bb_2}
 	+
 	1
 	\right)^{\frac{2-\al}{2}}}\dd x<\infty,
 \end{gathered}\end{equation*}
 since $-\bb_1\al+\frac{\bb_1}{\bb_2}+1=\bb_1\left( -\al+\frac{1}{\bb_2}+\frac{1}{\bb_1}\right)=\bb_1\left( \frac{1}{\bb_2}-1\right)<0$.

 The Dominated Convergence Theorem implies
 \begin{equation*}\begin{gathered}
 \frac{\tilde{D}_4}{n^{-\frac{\bb_1}{\bb_2}}m_n^{1-\bb_2\al}}\rightarrow
 \sum_{i=0}^{\infty}\int_{0}^{\infty}
 \frac
 {w_{(i,\infty)}\left(1+i\right)^{-\bb_1}
 	s^{-\bb_2} }
 {\left(
 	w_{(i,\infty)}^2
 	\left(1+i\right)^{-2\bb_1}+
 	s^{-2\bb_2} 	
 	\right)^{\frac{2-\al}{2}}}\dd s,\ n\rightarrow\infty.
 \end{gathered}\end{equation*}

Together with \eqref{D4kr} this implies
\begin{equation*}
\frac{{D}_4}{n^{-\frac{\bb_1}{\bb_2}}m_n^{1-\bb_2\al}}\rightarrow
\sum_{i=0}^{\infty}\int_{0}^{\infty}
\frac
{w_{(i,\infty)}\left(1+i\right)^{-\bb_1}
	s^{-\bb_2} }
{\left(
	w_{(i,\infty)}^2
	\left(1+i\right)^{-2\bb_1}+
	s^{-2\bb_2} 	
	\right)^{\frac{2-\al}{2}}}\dd s,\ n\rightarrow\infty.
\end{equation*}

Let us now show that
\begin{equation*}
\frac{D_1+D_2+D_3}{n^{-\frac{\bb_1}{\bb_2}}m_n^{1-\bb_2\al}}\rightarrow 0,\ n\rightarrow\infty.
\end{equation*}
This fact easily follows from \eqref{D123}:
\begin{equation*}
\frac{\abs{D_1+D_2+D_3}}{n^{-\frac{\bb_1}{\bb_2}}m_n^{1-\bb_2\al}} =
\frac{\abs{D_1+D_2+D_3}}{n^{-\bb_1}m_n^{1-\bb_2\al}\ln (n)}\frac{n^{-\bb_1}m_n^{1-\bb_2\al}\ln (n)}{n^{-\frac{\bb_1}{\bb_2}}m_n^{1-\bb_2\al}}\leq
K\frac{n^{-\bb_1}m_n^{1-\bb_2\al}\ln (n)}{n^{-\frac{\bb_1}{\bb_2}}m_n^{1-\bb_2\al}}
\end{equation*}
\begin{equation*}
=
K{n^{\frac{\bb_1}{\bb_2}-\bb_1}\ln (n)}=
K{n^{\frac{\bb_1}{\bb_2}(1-\bb_2)}\ln (n)}\rightarrow 0,\ n\rightarrow\infty.
\end{equation*}

In conclusion, we have
 \begin{equation*}
 \lim_{n\rightarrow\infty}\frac{\rho(n,-m_n)}{ m_n^{1-\bb_2\al}n^{-\frac{\bb_1}{\bb_2}}}=
 \lim_{n\rightarrow\infty}\frac{D_4}{ m_n^{1-\bb_2\al}n^{-\frac{\bb_1}{\bb_2}}}=
\sum_{i=0}^{\infty}\int_{0}^{\infty}
\frac
{w_{(i,\infty)}\left(1+i\right)^{-\bb_1}
	s^{-\bb_2} }
{\left(
	w_{(i,\infty)}^2
	\left(1+i\right)^{-2\bb_1}+
	s^{-2\bb_2} 	
	\right)^{\frac{2-\al}{2}}}\dd s.
 \end{equation*}

 $\bullet$ Let us now investigate the behavior of $D_4$ under the assumption $\frac{m_n^{-\bb_2}}{n^{-\bb_1}}\rightarrow \infty$ as $n\rightarrow\infty$.

 For convenience of writing let us denote $f_{n}=nm_n^{-\frac{\bb_2}{\bb_1}}$. Then $f_n\rightarrow \infty,\ n\rightarrow\infty.$ We have
 \begin{equation*}
 D_4=
 \sum_{i=0}^{\sv{n\varepsilon}-1}\sum_{j=0}^{\sv{m_n\varepsilon}-1}W_{i,j,n}.
 \end{equation*}
 For fixed $i$ we have
 \begin{equation}\label{fixedI}\begin{gathered}
 \abs{\sum_{j=0}^{\sv{m_n\varepsilon}-1}W_{i,j,n}}\leq
 \frac{e^2}{d^{2-\al}}\sum_{j=0}^{\sv{m_n\varepsilon}-1}
 \frac {(1+i+n)^{-\bb_1}(1+j)^{-\bb_2}(1+i)^{-\bb_1}(1+j+m_n)^{-\bb_2}}
 {\left(
 	(1+i+n)^{-2\bb_1}(1+j)^{-2\bb_2}+
 	(1+i)^{-2\bb_1}(1+j+m_n)^{-2\bb_2} \right)^{\frac{2-\al}{2}}}\\ \leq
 \frac{e^2}{d^{2-\al}}\sum_{j=0}^{\sv{m_n\varepsilon}-1}
 (1+i+n)^{-\bb_1}(1+j)^{-\bb_2}(1+i)^{-\bb_1(\al-1)}(1+j+m_n)^{-\bb_2(\al-1)}\\ \leq
 \frac{e^2}{d^{2-\al}}n^{-\bb_1}m_n^{-\bb_2(\al-1)}
 \sum_{j=0}^{\infty}
 (1+j)^{-\bb_2}.
 \end{gathered}\end{equation}
 Since $\bb_2>1$, the sum is finite. For any $a\in\NN$ we can split
 \begin{equation*}
 D_4=\tilde{F}_4(a)+ \tilde{D}_4(a),
 \end{equation*}
 where
 \begin{equation*}
 \tilde{F}_4(a)=\sum_{i=0}^{a-1}\sum_{j=0}^{\sv{m_n\varepsilon}-1}W_{i,j,n} , \quad
 \tilde{D}_4(a)=\sum_{i=a}^{\sv{n\varepsilon}-1}\sum_{j=0}^{\sv{m_n\varepsilon}-1}W_{i,j,n}.
 \end{equation*}
 Inequality \eqref{fixedI} gives us
 \begin{equation}\label{F4bound}\begin{gathered}
 {\abs{\tilde{F}_4(a)}}
 \leq a \frac{e^2}{d^{2-\al}}
 {n^{-\bb_1}m_n^{-\bb_2(\al-1)}}
 \sum_{j=0}^{\infty}
 (1+j)^{-\bb_2}
 \end{gathered}\end{equation}

 We continue with $\tilde{D}_4(a)$:
 \begin{equation*}\footnotesize\begin{gathered}
 \tilde{D}_4(a)= \sum_{i=a}^{\sv{n\varepsilon}-1}\sum_{j=0}^{\sv{m_n\varepsilon}-1}
 \frac
 {w_{(i+n,j)}w_{(i,j+m_n)}(1+i+n)^{-\bb_1}(1+j)^{-\bb_2}(1+i)^{-\bb_1}(1+j+m_n)^{-\bb_2}}
 {\left(
 	w^2_{(i+n,j)}(1+i+n)^{-2\bb_1}(1+j)^{-2\bb_2}+
 	w^2_{(i,j+m_n)}(1+i)^{-2\bb_1}(1+j+m_n)^{-2\bb_2} \right)^{\frac{2-\al}{2}}}\\=
 \int\limits_{a}^{\sv{n\varepsilon}}\sum_{j=0}^{\sv{m_n\varepsilon}-1}
 \frac
 {w_{(\sv{t}+n,j)}w_{(\sv{t},j+m_n)}(1+\sv{t}+n)^{-\bb_1}(1+j)^{-\bb_2}(1+\sv{t})^{-\bb_1}(1+j+m_n)^{-\bb_2}}
 {\left(
 	w^2_{(\sv{t}+n,j)}(1+\sv{t}+n)^{-2\bb_1}(1+j)^{-2\bb_2}+
 	w^2_{(\sv{t},j+m_n)}(1+\sv{t})^{-2\bb_1}(1+j+m_n)^{-2\bb_2} \right)^{\frac{2-\al}{2}}}  \dd t
 \\
 =
 \int\limits_{\frac{a}{f_n}}^{\frac{\sv{n\varepsilon}}{n}m_n^{\frac{\bb_2}{\bb_1}}}\sum_{j=0}^{\sv{m_n\varepsilon}-1}
 \frac
 {f_n w_{(\sv{f_nt}+n,j)}w_{(\sv{f_nt},j+m_n)}(1+\sv{f_nt}+n)^{-\bb_1}(1+j)^{-\bb_2}(1+\sv{f_nt})^{-\bb_1}(1+j+m_n)^{-\bb_2}}
 {\left(
 	w^2_{(\sv{f_nt}+n,j)}(1+\sv{f_nt}+n)^{-2\bb_1}(1+j)^{-2\bb_2}+
 	w^2_{(\sv{f_nt},j+m_n)}(1+\sv{f_nt})^{-2\bb_1}(1+j+m_n)^{-2\bb_2} \right)^{\frac{2-\al}{2}}}  \dd t
 \end{gathered}\end{equation*}
 \begin{equation*}\scriptsize\begin{gathered}
 =
 \int\limits_{\frac{a}{f_n}}^{\frac{\sv{n\varepsilon}}{n}m_n^{\frac{\bb_2}{\bb_1}}}\sum_{j=0}^{\sv{m_n\varepsilon}-1}
 \frac
 {f_n^{1-\bb_1}n^{-\bb_1(\al-1)}m_n^{-\bb_2}
 	w_{(\sv{f_nt}+n,j)}
 	w_{(\sv{f_nt},j+m_n)}
 	\left(1+\frac{\sv{f_nt}+1}{n}\right)^{-\bb_1}
 	\left(1+j\right)^{-\bb_2}
 	\left(\frac{1+\sv{f_nt}}{f_n}\right)^{-\bb_1}
 	\left(1+\frac{j+1}{m_n}\right)^{-\bb_2}}
 {\left(
 	w^2_{(\sv{f_nt}+n,j)}
 	\left(1+\frac{\sv{f_nt}+1}{n}\right)^{-2\bb_1}
 	\left(1+j\right)^{-2\bb_2}+
 	w^2_{(\sv{f_nt},j+m_n)}
 	\left(\frac{1+\sv{f_nt}}{f_n}\right)^{-2\bb_1}
 	\left(1+\frac{j+1}{m_n}\right)^{-2\bb_2}
 	\right)^{\frac{2-\al}{2}}}  \dd t
 \end{gathered}\end{equation*}
  \begin{equation*}\begin{gathered}
  =
  n^{1-\bb_1\al}m_n^{-\frac{\bb_2}{\bb_1}}
  \int\limits_{\frac{a}{f_n}}^{\frac{\sv{n\varepsilon}}{n}m_n^{\frac{\bb_2}{\bb_1}}}\sum_{j=0}^{\sv{m_n\varepsilon}-1}
  q_{t,j,n} \dd t.
  \end{gathered}\end{equation*}

If $\frac{a}{f_n}<t<\frac{\sv{n\varepsilon}}{n}m_n^{\frac{\bb_2}{\bb_1}}$, $0\leq j\leq \sv{m_n\varepsilon}-1$ and $n\geq\frac{2}{\varepsilon}$ the following inequalities hold:
 \begin{equation*}
 \left( 1+\frac{1}{a} \right)^{-\bb_1}t^{-\bb_1}\leq
 \left(\frac{1+\sv{f_nt}}{f_n}\right)^{-\bb_1}\leq t^{-\bb_1},
 \end{equation*}
 \begin{equation*}
 \left( 1+\varepsilon \right)^{-\bb_2}\leq
 \left(1+\frac{j+1}{m_n}\right)^{-\bb_2}\leq
 1,
 \end{equation*}
 \begin{equation*}
 \left( 1+\frac{3}{2}\varepsilon \right)^{-\bb_1}\leq \left(1+\frac{\sv{f_nt}+1}{n}\right)^{-\bb_1} \leq 1.
 \end{equation*}
 Furthermore, since $w_{(i,j)}$ uniformly converges to $w_{(\infty,j)}$, as $i \rightarrow\infty$, and $w_{(i,j)}\rightarrow 1,\ i,j\rightarrow\infty$, for all $a,n,m_n$ large enough the following inequalities hold
 \begin{equation*}
 w_{(\infty,j)}(1-\varepsilon)\leq C(\sv{f_nt}+n,j)\leq w_{(\infty,j)}(1+\varepsilon),
 \end{equation*}
 \begin{equation*}
 1-\varepsilon \leq C(\sv{f_nt},j+m_n)\leq 1+\varepsilon.
 \end{equation*}

 Previous inequalities imply the following estimate from above:
 \begin{equation*}\small\begin{gathered}
 q_{t,j,n} \leq
 \frac
 {
 	w_{(\infty,j)}\left(1+\varepsilon \right)
 	\left( 1+\varepsilon \right)
 	\left(1+j\right)^{-\bb_2}
 	t^{-\bb_1}}
 {\left(
 	w_{(\infty,j)}^2
 	\left(1-\varepsilon \right)^2
 	\left( 1+\frac{3}{2}\varepsilon \right)^{-2\bb_1}
 	\left(1+j\right)^{-2\bb_2}+
 	\left( 1-\varepsilon \right)^2
 	\left( 1+\frac{1}{a} \right)^{-2\bb_1}t^{-2\bb_1}
 	\left( 1+\varepsilon \right)^{-2\bb_2}
 	\right)^{\frac{2-\al}{2}}}\\ \leq
 \frac
 {
 	\left( 1+\varepsilon \right)
 	\left( 1+\varepsilon \right)
 } {\left( \left( 1-\varepsilon \right)^2
 \left( 1+\frac{3}{2}\varepsilon \right)^{-2\bb_1}
 \left( 1+\frac{1}{a} \right)^{-2\bb_1}
 \left( 1+\varepsilon \right)^{-2\bb_2} \right)^{\frac{2-\al}{2}}}
\frac
{	w_{(\infty,j)}\left(1+j\right)^{-\bb_2}
	t^{-\bb_1}}
{\left(   	
	w_{(\infty,j)}^2\left(1+j\right)^{-2\bb_2}+  	
	t^{-2\bb_1}
	\right)^{\frac{2-\al}{2}}}\\ =
R_d
\frac
{	w_{(\infty,j)}\left(1+j\right)^{-\bb_2}
	t^{-\bb_1}}
{\left(   	
	w_{(\infty,j)}^2\left(1+j\right)^{-2\bb_2}+  	
	t^{-2\bb_1}
	\right)^{\frac{2-\al}{2}}},
\end{gathered}\end{equation*}
where
\begin{equation*}
R_d=\frac
{
	\left( 1+\varepsilon \right)^2
} {\left(
\left( 1-\varepsilon \right)^2
\left( 1+\frac{3}{2}\varepsilon \right)^{-2\bb_1}
\left( 1+\frac{1}{a} \right)^{-2\bb_1}
\left( 1+\varepsilon \right)^{-2\bb_2} \right)^{\frac{2-\al}{2}}}.
\end{equation*}
The estimate from below is as follows:
\begin{equation*}\small\begin{gathered}
q_{t,j,n}\geq
\frac
{
	w_{(\infty,j)}	\left(1-\varepsilon \right)
	\left( 1-\varepsilon \right)
	\left( 1+\frac{3}{2}\varepsilon \right)^{-\bb_1}
	\left(1+j\right)^{-\bb_2}
	\left( 1+\frac{1}{a} \right)^{-\bb_1}t^{-\bb_1}
	\left(1+\varepsilon \right)^{-\bb_2}}
{\left(
	w_{(\infty,j)}^2\left( 1+\varepsilon \right)^2
	\left(1+j\right)^{-2\bb_2}+
	\left( 1+\varepsilon \right)^2
	t^{-2\bb_1}
	\right)^{\frac{2-\al}{2}}} \\=
R_k
\frac{w_{(\infty,j)}\left(1+j\right)^{-\bb_2}t^{-\bb_1}}{\left(
	w_{(\infty,j)}^2
	\left(1+j\right)^{-2\bb_2}+
	t^{-2\bb_1}
	\right)^{\frac{2-\al}{2}}} ,
\end{gathered}\end{equation*}
where
\begin{equation*}
R_k=\frac
{
	\left(1-\varepsilon \right)^2
	\left( 1+\frac{3}{2}\varepsilon \right)^{-\bb_1}
	\left( 1+\frac{1}{a} \right)^{-\bb_1}
	\left(1+\varepsilon \right)^{-\bb_2}}
{\left( 1+\varepsilon \right)^{2-\al}} .
\end{equation*}
The multipliers $R_k$ and $R_d$ converge to 1, as $\varepsilon\rightarrow 0,\ a\rightarrow\infty$. Let us now investigate
\begin{equation*}
\mathfrak{D} =\sum_{j=0}^{\sv{m_n\varepsilon}-1}
\int\limits_{\frac{a}{f_n}}^{\frac{\sv{n\varepsilon}}{n}m_n^{\frac{\bb_2}{\bb_1}}}
\frac{w_{(\infty,j)}\left(1+j\right)^{-\bb_2}t^{-\bb_1}}{\left(
	w_{(\infty,j)}^2
	\left(1+j\right)^{-2\bb_2}+
	t^{-2\bb_1}
	\right)^{\frac{2-\al}{2}}}\dd t.
\end{equation*}
The following inequalities hold
\begin{equation}\label{tildeDvertinimas}
R_k\mathfrak{D}
\leq
\frac{\tilde{D}_4}{n^{1-\bb_1\al}m_n^{-\frac{\bb_2}{\bb_1}}}
\leq
R_d\mathfrak{D}.
\end{equation}

{The change of the variables $x=a_j^{\frac{1}{\bb_1}}t(1+j)^{-\frac{\bb_2}{\bb_1}}$ enables us to write}
\begin{equation*}\begin{gathered}
\mathfrak{D} =\sum_{j=0}^{\sv{m_n\varepsilon}-1}
\int\limits_{\frac{a}{f_n}}^{\frac{\sv{n\varepsilon}}{n}m_n^{\frac{\bb_2}{\bb_1}}}
\frac{w_{(\infty,j)}\left(1+j\right)^{-\bb_2}t^{-\bb_1}}{\left(
	w_{(\infty,j)}^2
	\left(1+j\right)^{-2\bb_2}+
	t^{-2\bb_1}
	\right)^{\frac{2-\al}{2}}}\dd t\\=
\sum_{j=0}^{\sv{m_n\varepsilon}-1}
\int\limits_{\frac{a}{f_n}}^{\frac{\sv{n\varepsilon}}{n}m_n^{\frac{\bb_2}{\bb_1}}}
\frac{w_{(\infty,j)}^{\al-1}\left(1+j\right)^{-\bb_2(\al-1)}t^{-\bb_1}}{\left(
	1+
	w_{(\infty,j)}^{-2}
	\left(1+j\right)^{2\bb_2}t^{-2\bb_1}
	\right)^{\frac{2-\al}{2}}}\dd t
\end{gathered}\end{equation*}
\begin{equation*}\begin{gathered}
=
\sum_{j=0}^{\sv{m_n\varepsilon}-1}
\int\limits_
{\frac{a}{f_n}w_{(\infty,j)}^{\frac{1}{\bb_1}}(1+j)^{-\frac{\bb_2}{\bb_1}}}
^{\frac{\sv{n\varepsilon}}{n}m_n^{\frac{\bb_2}{\bb_1}}w_{(\infty,j)}^{\frac{1}{\bb_1}}(1+j)^{-\frac{\bb_2}{\bb_1}}}
\frac{
	w_{(\infty,j)}^{\al-1}\left(1+j\right)^{-\bb_2(\al-1)}w_{(\infty,j)}(1+j)^{-\bb_2}x^{-\bb_1}
}{\left(
1+
x^{-2\bb_1}
\right)^{\frac{2-\al}{2}}}w_{(\infty,j)}^{-\frac{1}{\bb_1}}(1+j)^{\frac{\bb_2}{\bb_1}}\dd x
\\=
\sum_{j=0}^{\sv{m_n\varepsilon}-1}
\int\limits_
{\frac{a}{f_n}w_{(\infty,j)}^{\frac{1}{\bb_1}}(1+j)^{-\frac{\bb_2}{\bb_1}}}
^{\frac{\sv{n\varepsilon}}{n}m_n^{\frac{\bb_2}{\bb_1}}w_{(\infty,j)}^{\frac{1}{\bb_1}}(1+j)^{-\frac{\bb_2}{\bb_1}}}
\frac{
	w_{(\infty,j)}^{\al-\frac{1}{\bb_1}}\left(1+j\right)^{-\bb_2\al+\frac{\bb_2}{\bb_1}}x^{-\bb_1}
}{\left(
1+
x^{-2\bb_1}
\right)^{\frac{2-\al}{2}}}\dd x
\end{gathered}\end{equation*}
\begin{equation}\label{pagr}\begin{gathered}
=
\sum_{j=0}^{\sv{m_n\varepsilon}-1}
w_{(\infty,j)}\left(1+j\right)^{-\bb_2}
\int\limits_
{\frac{a}{f_n}w_{(\infty,j)}^{\frac{1}{\bb_1}}(1+j)^{-\frac{\bb_2}{\bb_1}}}
^{\frac{\sv{n\varepsilon}}{n}m_n^{\frac{\bb_2}{\bb_1}}w_{(\infty,j)}^{\frac{1}{\bb_1}}(1+j)^{-\frac{\bb_2}{\bb_1}}}
\frac{
	x^{-\bb_1(\al-1)}
}{\left(
1+
x^{2\bb_1}
\right)^{\frac{2-\al}{2}}}\dd x.
\end{gathered}\end{equation}
For a moment, let us investigate the integral
\begin{equation}\label{BjaurusisIntegralas}
\int\limits_
{\frac{a}{f_n}w_{(\infty,j)}^{\frac{1}{\bb_1}}(1+j)^{-\frac{\bb_2}{\bb_1}}}
^{\frac{\sv{n\varepsilon}}{n}m_n^{\frac{\bb_2}{\bb_1}}w_{(\infty,j)}^{\frac{1}{\bb_1}}(1+j)^{-\frac{\bb_2}{\bb_1}}}
\frac{
	x^{-1}
}{\left(
1+
x^{2\bb_1}
\right)^{\frac{2-\al}{2}}}\dd x
\end{equation}
separately.

Variable $j$ satisfies the inequality $j\leq \sv{m_n\varepsilon}-1$, so the upper limit of the integral satisfies
\begin{equation*}
\frac{\sv{n\varepsilon}}{n}m_n^{\frac{\bb_2}{\bb_1}}a^{\frac{1}{\bb_1}}_j(1+j)^{-\frac{\bb_2}{\bb_1}}=
\frac{\sv{n\varepsilon}}{n}a^{\frac{1}{\bb_1}}_j\left(\frac{1+j}{m_n}\right)^{-\frac{\bb_2}{\bb_1}}
\geq \frac{\varepsilon}{2} A^{\frac{1}{\bb_1}}
\varepsilon^{-\frac{\bb_2}{\bb_1}}=\frac{A^{\frac{1}{\bb_1}}}{2}\varepsilon^{1-\frac{\bb_2}{\bb_1}},  \end{equation*}
with $A=\inf\left( a_j, j\geq 0 \right)$. Notice that $1-\frac{\bb_2}{\bb_1}=1-\bb_2(\al-1)>0$. Denote $\delta=\frac{A^{\frac{1}{\bb_1}}}{2}\varepsilon^{1-\frac{\bb_2}{\bb_1}}$.

Let us split the integral in \eqref{BjaurusisIntegralas} as
\begin{equation*}\begin{gathered}
\int\limits_
{\frac{a}{f_n}w_{(\infty,j)}^{\frac{1}{\bb_1}}(1+j)^{-\frac{\bb_2}{\bb_1}}}
^{\frac{\sv{n\varepsilon}}{n}m_n^{\frac{\bb_2}{\bb_1}}w_{(\infty,j)}^{\frac{1}{\bb_1}}(1+j)^{-\frac{\bb_2}{\bb_1}}}
\frac{
	x^{-1}
}{\left(
1+
x^{2\bb_1}
\right)^{\frac{2-\al}{2}}}\dd x\\=
\int\limits_
{\frac{a}{f_n}w_{(\infty,j)}^{\frac{1}{\bb_1}}(1+j)^{-\frac{\bb_2}{\bb_1}}}
^{\delta}
\frac{
	x^{-1}
}{\left(
1+
x^{2\bb_1}
\right)^{\frac{2-\al}{2}}}\dd x+
\int\limits_
{\delta}
^{\frac{\sv{n\varepsilon}}{n}m_n^{\frac{\bb_2}{\bb_1}}w_{(\infty,j)}^{\frac{1}{\bb_1}}(1+j)^{-\frac{\bb_2}{\bb_1}}}
\frac{
	x^{-1}
}{\left(
1+
x^{2\bb_1}
\right)^{\frac{2-\al}{2}}}\dd x=:I_{1,j}+I_{2,j}.
\end{gathered}\end{equation*}

The second integral is bounded from above by
\begin{equation*}\begin{gathered}
\int_
{\delta}
^{\frac{\sv{n\varepsilon}}{n}m_n^{\frac{\bb_2}{\bb_1}}w_{(\infty,j)}^{\frac{1}{\bb_1}}(1+j)^{-\frac{\bb_2}{\bb_1}}}
\frac{
	x^{-1}
}{\left(
1+
x^{2\bb_1}
\right)^{\frac{2-\al}{2}}}\dd x\leq
\int_
{\delta}
^{\infty}
{
	x^{-\bb_1}
}\dd x<\infty,
\end{gathered}\end{equation*}
since $\bb_1=\frac{1}{\al-1}>1$.

We proceed with the first integral. The integration variable $x$ belongs to the interval $(0,\delta)$, so
\begin{equation}\label{PointVertNulioApl}\begin{gathered}
\frac{
	x^{-1}
}{\left(
1+
\delta^{2\bb_1}
\right)^{\frac{2-\al}{2}}}
\leq
\frac{
	x^{-1}
}{\left(
1+
x^{2\bb_1}
\right)^{\frac{2-\al}{2}}}\leq x^{-1}.
\end{gathered}\end{equation}

 Since simple integration gives
\begin{equation*}
\int\limits_
{\frac{a}{f_n}a^{\frac{1}{\bb_1}}_j(1+j)^{-\frac{\bb_2}{\bb_1}}}
^{\delta}x^{-1}\dd x= \theta_j+\ln\left( f_n \right),
\end{equation*}
where $\theta_j=\ln(\delta) - \ln\left( {a}w_{(\infty,j)}^{\frac{1}{\bb_1}}(1+j)^{-\frac{\bb_2}{\bb_1}} \right)$,
integrating \eqref{PointVertNulioApl} we obtain
\begin{equation}\label{I1jVert}
\frac{
	\theta_j+\ln(f_n)
}{\left(
1+
\delta^{2\bb_1}
\right)^{\frac{2-\al}{2}}}\leq I_{1,j}\leq \theta_j+\ln(f_n).
\end{equation}

We can split
\begin{equation}\label{Dsplit}
\mathfrak{D}=\sum_{j=0}^{\sv{m_n\varepsilon}-1}w_{(\infty,j)}\left(1+j\right)^{-\bb_2}I_{1,j}+
\sum_{j=0}^{\sv{m_n\varepsilon}-1}w_{(\infty,j)}\left(1+j\right)^{-\bb_2}I_{2,j}.
\end{equation}

The second term is bounded from above by
\begin{equation}\label{DsplitI2bound}
\sum_{j=0}^{\sv{m_n\varepsilon}-1}w_{(\infty,j)}\left(1+j\right)^{-\bb_2}I_{2,j}\leq
\sum_{j=0}^{\infty}w_{(\infty,j)}\left(1+j\right)^{-\bb_2}\int_
{\delta}
^{\infty}
{
	x^{-\bb_1}
}\dd x<\infty.
\end{equation}

Let us denote the first term as $\bar{S}_1$. Then with the help of \eqref{I1jVert} we can obtain the bounds
\begin{equation*}
\sum_{j=0}^{\sv{m_n\varepsilon}-1}w_{(\infty,j)}\left(1+j\right)^{-\bb_2}\frac{
	\theta_j+\ln(f_n)
}{\left(
1+
\delta^{2\bb_1}
\right)^{\frac{2-\al}{2}}}
\leq \bar{S}_1\leq
\sum_{j=0}^{\sv{m_n\varepsilon}-1}w_{(\infty,j)}\left(1+j\right)^{-\bb_2}\left( \theta_j+\ln(f_n) \right).
\end{equation*}

Notice that
\begin{equation*}\begin{gathered}
\sum_{j=0}^{\infty}\left|w_{(\infty,j)}\left(1+j\right)^{-\bb_2}\theta_j\right|\\\leq
\sum_{j=0}^{\infty}w_{(\infty,j)}\left(1+j\right)^{-\bb_2} \left|\ln(\delta) - \ln\left( {a}w_{(\infty,j)}^{\frac{1}{\bb_1}}\right)\right|+
\frac{\bb_2}{\bb_1}\sum_{j=0}^{\infty}w_{(\infty,j)}\left(1+j\right)^{-\bb_2}\ln\left(1+j \right)<\infty,
\end{gathered}\end{equation*}
as $\left|\ln(\delta) - \ln\left( {a}w_{(\infty,j)}^{\frac{1}{\bb_1}}\right)\right|$ is bounded and $\bb_2>1$.

We obtain
\begin{equation*}
\frac{ 1 }{\left(
	1+
	\delta^{2\bb_1}
	\right)^{\frac{2-\al}{2}}}\sum_{j=0}^{\infty}w_{(\infty,j)}\left(1+j\right)^{-\bb_2}\leq
\liminf\limits_{n\rightarrow\infty}\frac{\bar{S}}{\ln(f_n)},\quad
\limsup\limits_{n\rightarrow\infty}\frac{\bar{S}}{\ln(f_n)} \leq \sum_{j=0}^{\infty}w_{(\infty,j)}\left(1+j\right)^{-\bb_2}.
\end{equation*}
Together with \eqref{Dsplit} and \eqref{DsplitI2bound} this implies
\begin{equation*}
\frac{ 1 }{\left(
	1+
	\delta^{2\bb_1}
	\right)^{\frac{2-\al}{2}}}\sum_{j=0}^{\infty}w_{(\infty,j)}\left(1+j\right)^{-\bb_2}\leq
\liminf\limits_{n\rightarrow\infty}\frac{\mathfrak{D}}{\ln(f_n)},\quad
\limsup\limits_{n\rightarrow\infty}\frac{\mathfrak{D}}{\ln(f_n)} \leq \sum_{j=0}^{\infty}w_{(\infty,j)}\left(1+j\right)^{-\bb_2}.
\end{equation*}
Returning to \eqref{tildeDvertinimas} and keeping in mind \eqref{F4bound} we obtain

\begin{equation*}
\frac{ R_k }{\left(
	1+
	\delta^{2\bb_1}
	\right)^{\frac{2-\al}{2}}}\sum_{j=0}^{\infty}w_{(\infty,j)}\left(1+j\right)^{-\bb_2}
\leq
\liminf\limits_{n\rightarrow\infty}\frac{{D}_4}{n^{1-\bb_1\al}m_n^{-\frac{\bb_2}{\bb_1}}\ln(f_n)},
\end{equation*}
\begin{equation*}
\limsup\limits_{n\rightarrow\infty}\frac{{D}_4}{n^{1-\bb_1\al}m_n^{-\frac{\bb_2}{\bb_1}}\ln(f_n)} \leq R_d\sum_{j=0}^{\infty}w_{(\infty,j)}\left(1+j\right)^{-\bb_2}.
\end{equation*}

Now,
\begin{equation*}
\frac{\abs{D_1+D_2+D_3}}{n^{1-\bb_1\al}m_n^{-\frac{\bb_2}{\bb_1}}\ln\left( nm^{-\frac{\bb_2}{\bb_1}} \right)} =
\frac{\abs{D_1+D_2+D_3}}{n^{-\bb_1}m_n^{1-\bb_2\al}\ln (n)}\frac{n^{-\bb_1}m_n^{1-\bb_2\al}\ln (n)}{n^{1-\bb_1\al}m_n^{-\frac{\bb_2}{\bb_1}}\ln\left( nm^{-\frac{\bb_2}{\bb_1}} \right)}
\end{equation*}
\begin{equation*}
\leq
K\frac{n^{-\bb_1}m_n^{1-\bb_2\al}\ln (n)}{n^{1-\bb_1\al}m_n^{-\frac{\bb_2}{\bb_1}}\ln\left( nm^{-\frac{\bb_2}{\bb_1}} \right)}=
K\frac{m_n^{1-\bb_2}\ln (n)}{\ln\left( nm^{-\frac{\bb_2}{\bb_1}} \right)}
=
K\frac{m_n^{1-\bb_2}\left( \ln \left( nm^{-\frac{\bb_2}{\bb_1}} \right)+\ln \left( m_n^{\frac{\bb_2}{\bb_1}} \right) \right)}{\ln\left( nm^{-\frac{\bb_2}{\bb_1}} \right)}
\end{equation*}
\begin{equation*}
=
K\left(m_n^{1-\bb_2}+ \frac{m_n^{1-\bb_2}\ln \left( m_n^{\frac{\bb_2}{\bb_1}} \right) }{\ln\left( nm^{-\frac{\bb_2}{\bb_1}} \right)} \right)\rightarrow0,\ n\rightarrow\infty.
\end{equation*}

This implies
\begin{equation*}
\frac{ R_k }{\left(
	1+
	\delta^{2\bb_1}
	\right)^{\frac{2-\al}{2}}}\sum_{j=0}^{\infty}w_{(\infty,j)}\left(1+j\right)^{-\bb_2}
\leq
\liminf\limits_{n\rightarrow\infty}\frac{\rho(n,-m_n)}{n^{1-\bb_1\al}m_n^{-\frac{\bb_2}{\bb_1}}\ln(f_n)},
\end{equation*}
\begin{equation*}
\limsup\limits_{n\rightarrow\infty}\frac{\rho(n,-m_n)}{n^{1-\bb_1\al}m_n^{-\frac{\bb_2}{\bb_1}}\ln(f_n)} \leq R_d\sum_{j=0}^{\infty}w_{(\infty,j)}\left(1+j\right)^{-\bb_2}.
\end{equation*}

Passing to the limit as $\varepsilon\rightarrow0$, then as $a\rightarrow\infty$, we obtain
\begin{equation*}
\lim\limits_{n\rightarrow\infty}\frac{ \rho({n,-m_n})}{n^{1-\bb_1\al}m_n^{-\frac{\bb_2}{\bb_1}}\ln\left(nm_n^{-\frac{\bb_2}{\bb_1}}\right)}=\sum_{j=0}^{\infty}w_{(\infty,j)}\left(1+j\right)^{-\bb_2}.
\end{equation*}

\end{document}